\documentclass[pdflatex,sn-nature]{sn-jnl}

\AtBeginDocument{\setcitestyle{super,sort&compress,open={},close={}}}

\usepackage{graphicx}%
\usepackage{multirow}%
\usepackage{amsmath,amssymb,amsfonts}%
\usepackage{amsthm}%
\usepackage{mathrsfs}%
\usepackage[title]{appendix}%
\usepackage{xcolor}%
\usepackage{textcomp}%
\usepackage{manyfoot}%
\usepackage{booktabs}%
\usepackage{algorithm}%
\usepackage{algorithmicx}%
\usepackage{algpseudocode}%
\usepackage{listings}%

\usepackage{graphicx}
\usepackage{subcaption}
\usepackage{array}
\usepackage{comment}
\usepackage{tabularx}

\usepackage{natbib}

\usepackage{booktabs}
\usepackage{siunitx}
\theoremstyle{thmstyleone}%

\theoremstyle{thmstyletwo}%
\newtheorem{remark}{Remark}%

\theoremstyle{thmstylethree}%

\begin{document}

\title[Article Title]{Geometry-aware LegONet for PDE Learning on Arbitrary Domains}


\author[1]{\fnm{Jiahao} \sur{Zhang}}
\author[1]{\fnm{Yueqi} \sur{Wang}}

\author[1,2]{\fnm{Guang} \sur{Lin}}


\affil[1]{\orgdiv{Department of Mathematics}, \orgname{Purdue University},
\orgaddress{\city{West Lafayette}, \state{IN}, \postcode{47907}, \country{USA}}}

\affil[2]{\orgdiv{School of Mechanical Engineering}, \orgname{Purdue University},
\orgaddress{\city{West Lafayette}, \state{IN}, \postcode{47907}, \country{USA}}}


\abstract{

Learned PDE solvers often entangle governing operators with the geometry, boundary conditions, and discretization used for training. This limits reuse when the same physics is posed on new domains, and it also makes physical-law discovery geometry-dependent. We introduce Geometry-aware LegONet (gLegONet), a boundary-manifold extension of Lego-like operator learning. Physical mechanisms are pretrained once as modular variational blocks on an ambient spectral domain. For a target geometry, sampled boundary constraints define an affine admissible manifold. Its mass-orthonormal tangent coordinates are used to evolve the dynamics and evaluate candidate law-discovery features directly. Changing the domain therefore changes only an algebraic coordinate interface, not the learned operator blocks. This converts arbitrary-domain PDE learning from geometry-specific retraining or soft penalty enforcement into boundary-guaranteed assembly of reusable mechanisms. In forward simulations and sparse identification tests on unseen domains, the method maintains boundary residuals near the algebraic tolerance and yields predictive governing laws from short-time observations.

}




\maketitle

Learned solvers for partial differential equations have advanced rapidly, including physics-informed networks, neural operators and learned simulators
\cite{sirignano2018dgm,raissi2019pinn,lu2021deeponet,li2020fno,
kovachki2023neuraloperator,pfaff2020meshgraphnets}. Their use in scientific computing, however, is still limited by geometry. A physical mechanism such as diffusion, transport or reaction can appear on many domains, but many learned solvers are trained on a fixed family of geometries. When the domain changes from a regular grid to a disk, channel or mixed-boundary region, the model may require new training data, architectural changes or an additional boundary treatment. This coupling between learned physics and geometric realization makes it difficult to build reusable learned solvers for related physical problems.

Classical numerical methods make this separation explicit. The PDE
specifies the interior operator, while the mesh, trial space, quadrature and boundary constraints determine how that operator is realized on a particular domain \citep{leveque2007fd,brenner2008fem,boyd2001chebyshev}. Neural solvers often mix these roles inside one trained model. Physics-informed neural networks can handle irregular domains through collocation, and some trial functions can impose selected boundary conditions exactly
\citep{raissi2019pinn,karniadakis2021piml,sukumar2022exact}. However, the
optimization is generally repeated for each target problem, and exact trial
constructions become increasingly specialized for complex, mixed or
derivative-valued boundary conditions.

Neural operators address this problem by amortizing inference across families of inputs. Foundational DeepONet, Fourier neural operator and U-shaped neural operator architectures were primarily developed on fixed domains or structured representations
\citep{lu2021deeponet,li2020fno,rahman2022uno,kovachki2023neuraloperator}.
Recent works have made neural operators more geometry-aware. Geo-FNO, DAFNO and GINO develop deformations, domain indicators or graph-to-Fourier interfaces
\citep{li2023fourier,liu2023domain,li2023geometry}.  GNOT and Transolver use
attention on irregular discretizations, while PI-GANO combines a geometry encoder with physics-informed operator training \citep{hao2023gnot,wu2024transolver,zhong2025physics}.  Boundary-focused
methods such as BENO further condition the solution operators on complex geometry
\citep{wang2024beno,mousavi2026imposing}.  These methods greatly
expand geometric flexibility. However, geometry and boundary data are still inputs to an end-to-end learned model, usually trained over a distribution of domains or boundary conditions. Our goal is different: to keep the learned physical mechanisms fixed and realize them on each new domain through a deterministic boundary-adapted Galerkin interface.

Our previous Lego-like variational operator learning framework \cite{zhang2026legonet} addressed one
part of this problem: \emph{physics modularity}. Rather than learning a
monolithic PDE solver, it decomposes an equation into reusable variational
operator blocks, such as diffusion, transport and reaction, and composes them
by structure-aware time integration
\citep{zhang2026legonet,strang1968splitting,hairer2006geometric}. That
formulation makes physical mechanisms reusable, but its most accurate
implementation still relied on boundary-compatible bases tailored to particular
baseplates. In other words, the physics became modular, but the geometry
remained built into the representation.

This work supplies the missing geometry interface. The key idea is to reuse a frozen ambient mechanism library on a new domain by changing coordinates, rather than by retraining the library or correcting the boundary after each step. Let \(\mathcal B_\theta^Q=\{\mathbf F_i^\theta\}_{i=1}^{N_{\rm lib}}\) be a library pretrained once on a regular domain \(Q\). We ask whether this library can be reused on new embedded target domains while enforcing new linear trace constraints algebraically.  A tempting alternative is to evolve an unconstrained ambient state and project it back to the admissible space after each time step. This projection--evolution splitting is generally not equivalent to integrating the \(M_\Omega\)-Galerkin vector field on the admissible tangent space, unless the projection metric, tangent space and time integrator are
specifically matched \citep{hairer2006geometric}.  We instead build the sampled boundary constraints into the coordinates before time integration begins.

For a target domain \(\Omega\subset Q\), consider a PDE decomposed into physical mechanisms as
\[
    \partial_t u(\mathbf x,t)
    =
    \sum_{i=1}^{N_{\rm blk}}
    c_i L_i^\Omega(u)(\mathbf x,t),
    \qquad
    \mathbf x\in\Omega ,
\]
where \(L_i^\Omega\) is the \(i\)th target-domain mechanism and \(c_i\) is its equation-specific scalar coefficient.  Sampling the boundary operator in the ambient coefficient space gives
\[
    C\mathbf a(t)=\mathbf d(t),
    \qquad
    \mathbf a(t)
    =
    \mathbf a_{\rm bc}(t)+N_\Omega\mathbf z(t).
\]
Here \(\mathbf a_{\rm bc}(t)\) is an affine lift of the prescribed boundary data. The columns of \(N_\Omega\) span the homogeneous boundary-compatible directions and are orthonormal in the \(M_\Omega\) inner product. Thus the state \(\mathbf a(t)\) satisfies the sampled boundary constraints by construction. For each mechanism \(L_i^\Omega\), the frozen ambient block is realized on \(\Omega\) as a reduced response \(\mathbf q_i^\Omega(\mathbf z,t)\in\mathbb R^r\).  The boundary-adapted dynamics are therefore
\[
    \dot{\mathbf z}(t)
    =
    \sum_{i=1}^{N_{\rm blk}}
    c_i\,\mathbf q_i^\Omega(\mathbf z(t),t)
    -
    N_\Omega^\top M_\Omega
    \dot{\mathbf a}_{\rm bc}(t).
\]
For static boundary data, \(\dot{\mathbf a}_{\rm bc}(t)=\mathbf0\). Thus this equation keeps the same physical coefficients \(c_i\), but replaces each mechanism by its boundary-adapted reduced response. Geometry and boundary conditions change only the algebraic realization layer
\((C,\mathbf a_{\rm bc},N_\Omega,M_\Omega)\). The neural parameters of the ambient mechanism library remain fixed.

This formulation clarifies the role of geometry transfer. Boundary conditions are not imposed by adding a penalty to the loss; they are built into the coordinates. The neural blocks are not retrained for each new domain; they are realized through a new algebraic interface. It is also distinct from post-step state projection, because the reduced dynamics are derived and integrated directly in the boundary-adapted coordinates. As a result, the sampled boundary constraints are satisfied to the tolerance of the null-space construction throughout the rollout. The main remaining errors come from the expressiveness of the boundary-adapted space, the accuracy of the frozen block responses, and the time discretization.

These boundary-adapted coordinates also provide a natural space for physical-law identification. Traditional sparse discovery methods such as SINDy and PDE-FIND identify governing terms from candidate libraries \citep{brunton2016discovering,rudy2017data,messenger2021weak}. In our setting, the candidate library is not rebuilt from raw spatial derivatives on each new geometry. Instead, frozen mechanism blocks and analytic terms are first realized through the same boundary-adapted Galerkin interface and then regressed in the reduced coordinate \(\mathbf z\).  Once their coefficients are identified, the discovered law already defines a boundary-constrained rollout model on the target domain.

\begin{figure}[H]
    \centering
    \includegraphics[width=\linewidth]{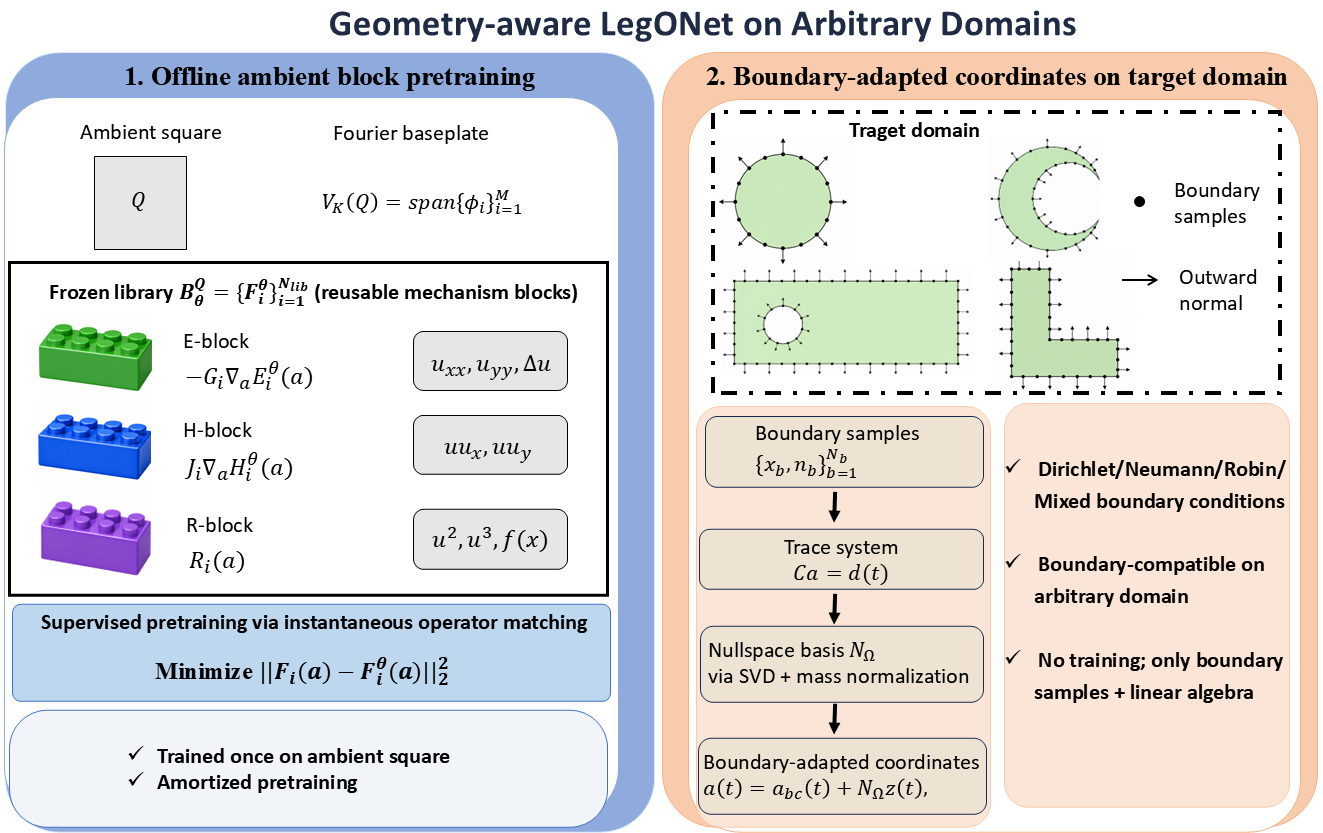}
    \hfill
    \includegraphics[width=\linewidth]{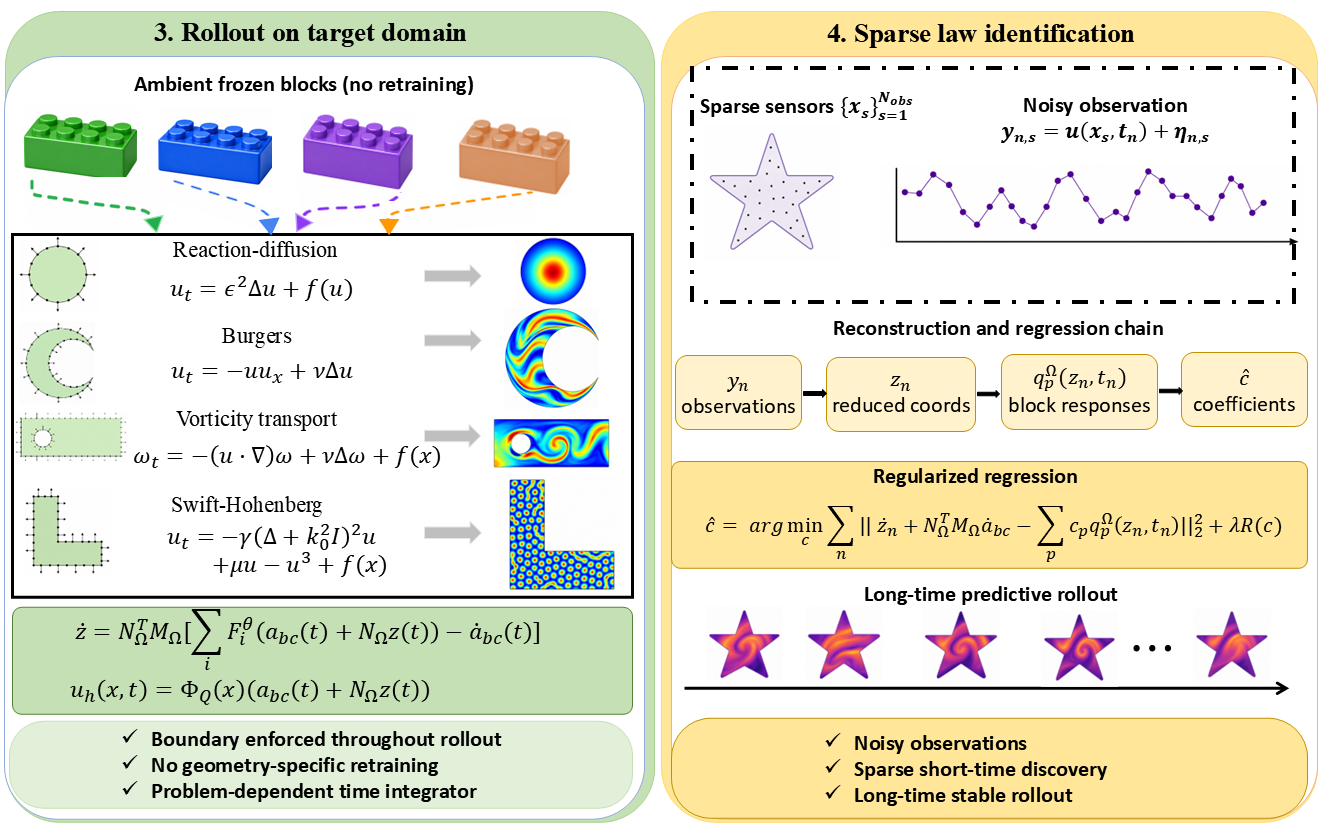}
\caption{
Workflow of Geometry-aware LegONet.
\textbf{1}, Single-mechanism blocks are pretrained once on the ambient square
\(Q\) and frozen as \(\mathcal B_\theta^Q\).
\textbf{2}, Boundary samples define \(C\mathbf a=\mathbf d(t)\); null-space and
mass normalization yield admissible coordinates
\(\mathbf a=\mathbf a_{\rm bc}+N_\Omega\mathbf z\).
\textbf{3}, The frozen blocks are assembled through the
\(\Omega\)-Galerkin dynamics and rolled out on new domains without neural
retraining.
\textbf{4}, The same coordinates map sparse short-time observations to reduced
states, enabling mechanism identification and direct long-time prediction.
}
    \label{fig:geometry_aware_workflow}
\end{figure}

Fig.~\ref{fig:geometry_aware_workflow} summarizes the resulting pretrain--realize--rollout pipeline.  We evaluate this idea on unseen geometries with diverse boundary conditions, nonlinear transport, obstacle flow, higher-order operator composition and short-time inverse discovery.  Across these tests, the boundary constraints are controlled by the algebraic coordinate

\noindent construction, while the same pretrained mechanisms are reused without geometry-specific neural retraining. Geometry-aware LegONet therefore connects modular neural mechanism learning with the classical
principle that interior physics and boundary-constrained realization should be separated.
\section*{Results}\label{sec:Results}

\subsection*{Reusable mechanism blocks are trained once and frozen}

All experiments reuse ambient libraries pretrained on the square
\(Q=[-1,1]^2\) before any target geometry or boundary condition is introduced.
For each Fourier cutoff \(K\), we train one resolution-matched library
\(\mathcal B_\theta^Q=
\{\mathbf F_\Delta^\theta,\mathbf F_x^\theta,\mathbf F_y^\theta\}\),
containing diffusion $\Delta u$ and the two directional quadratic-transport mechanisms $uu_x,uu_y$.
Once trained, the neural parameters are frozen. All subsequent target problems---regardless of geometry, boundary type or rollout equation---use the same pretrained blocks.
Unless stated otherwise, all experiments use \(K=22\).  Held-out block errors
remain low and stable over the tested resolutions
(Extended Data Figs.~\ref{fig:ambient_block_error_bins} and
\ref{fig:ambient_block_k_sensitivity}).

\subsection*{Frozen blocks transfer across diverse geometries and boundary conditions}

\begin{figure}[H]
\centering
\scriptsize
\setlength{\tabcolsep}{2.6pt}
\renewcommand{\arraystretch}{1.04}
\setlength{\abovecaptionskip}{2pt}
\setlength{\belowcaptionskip}{0pt}

\begin{minipage}[t]{0.995\textwidth}
\centering
\resizebox{\linewidth}{!}{
\begin{tabular}{llllcc}
\toprule
Test & Domain & Equation & Boundary condition & \(N_b\) & \(r\) \\
\midrule
MMS-I
& Rosette
& \(u_t=\varepsilon^2\Delta u+\mu u-\beta u^3+f^\Omega\)
& \(u=0\)
& \(420\)
& \(207\) \\

MMS-II
& Crescent
& \(u_t=\varepsilon^2\Delta u+\mu u-\beta u^3+f^\Omega\)
& \(\partial_{\mathbf n} u=0\)
& \(620\)
& \(220\) \\

MMS-III
& Bunny
& \(u_t=D\Delta u+r_g u(1-u/K_c)+f^\Omega\)
& \(\partial_{\mathbf n} u+\kappa u=g_R\)
& \(420\)
& \(843\) \\

MMS-IV
& Annular star
& \(\left\{
\begin{aligned}
u_t&=D_u\Delta u+a-u+u^2v+f_u^\Omega,\\
v_t&=D_v\Delta v+b-u^2v+f_v^\Omega
\end{aligned}
\right.\)
& \(\left\{
\begin{aligned}
&u=v=0 \quad \text{on }\Gamma_{\rm out},\\
&\partial_{\mathbf n} u=\partial_{\mathbf n} v=0 \quad \text{on }\Gamma_{\rm in}
\end{aligned}
\right.\)
& \(720\)
& \(832\) \\

MMS-V
& Pinwheel shell
& \(u_t=\nu\Delta u-u u_x-u u_y+f^\Omega\)
& \(u=0\)
& \(420\)
& \(644\) \\
\bottomrule
\end{tabular}}
\vspace{0.3mm}

{\scriptsize \textbf{(a)} Five manufactured-solution benchmarks.}
\end{minipage}

\vspace{0.8mm}

\noindent
\begin{minipage}[t]{0.496\textwidth}
    \centering
    {\scriptsize \textbf{(b)} MMS-I}\par
    \vspace{0.25mm}
    \begin{minipage}[c]{0.225\linewidth}
        \centering
        \includegraphics[
            height=0.121\textheight,
            keepaspectratio,
            trim=965 210 0.0 115,
            clip
        ]{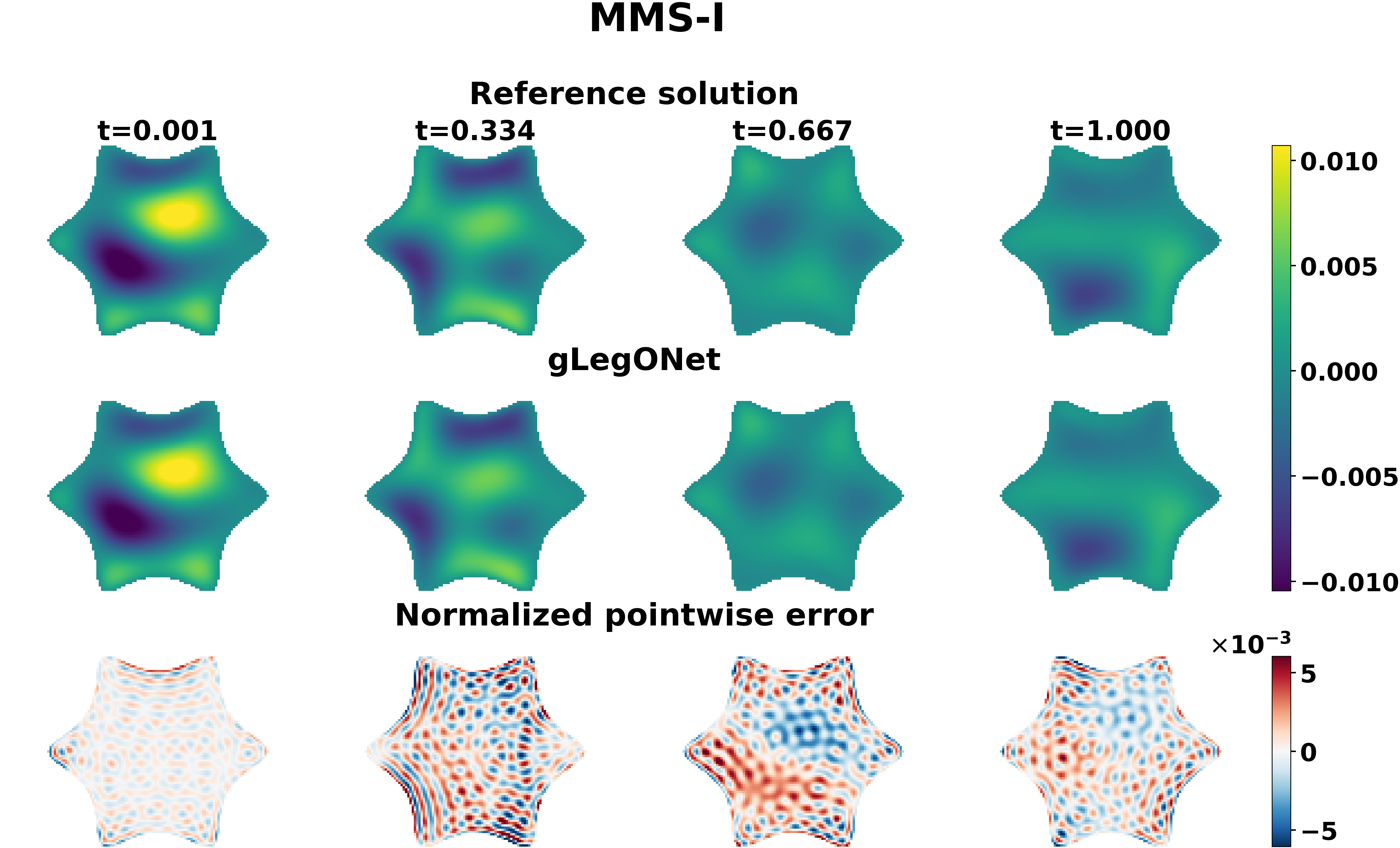}
    \end{minipage}%
    \begin{minipage}[c]{0.765\linewidth}
        \centering
        \includegraphics[
            width=\linewidth,
            height=0.125\textheight,
            keepaspectratio,
            trim=4 6 4 22,
            clip
        ]{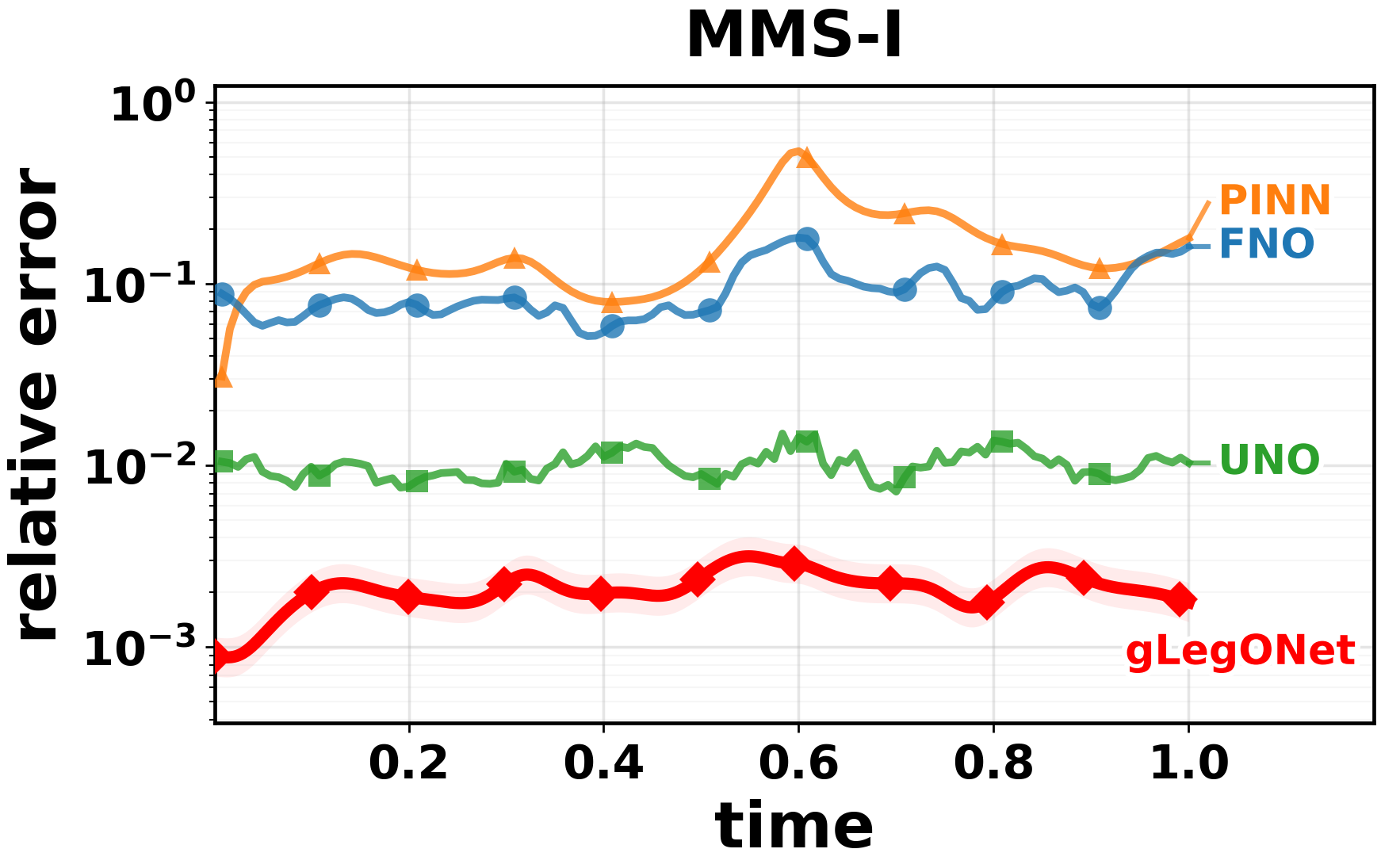}
    \end{minipage}
\end{minipage}%
\hspace{0.004\textwidth}%
\begin{minipage}[t]{0.496\textwidth}
    \centering
    {\scriptsize \textbf{(c)} MMS-II}\par
    \vspace{0.25mm}
    \begin{minipage}[c]{0.225\linewidth}
        \centering
        \includegraphics[
            height=0.125\textheight,
            keepaspectratio,
            trim=965 210 2 115,
            clip
        ]{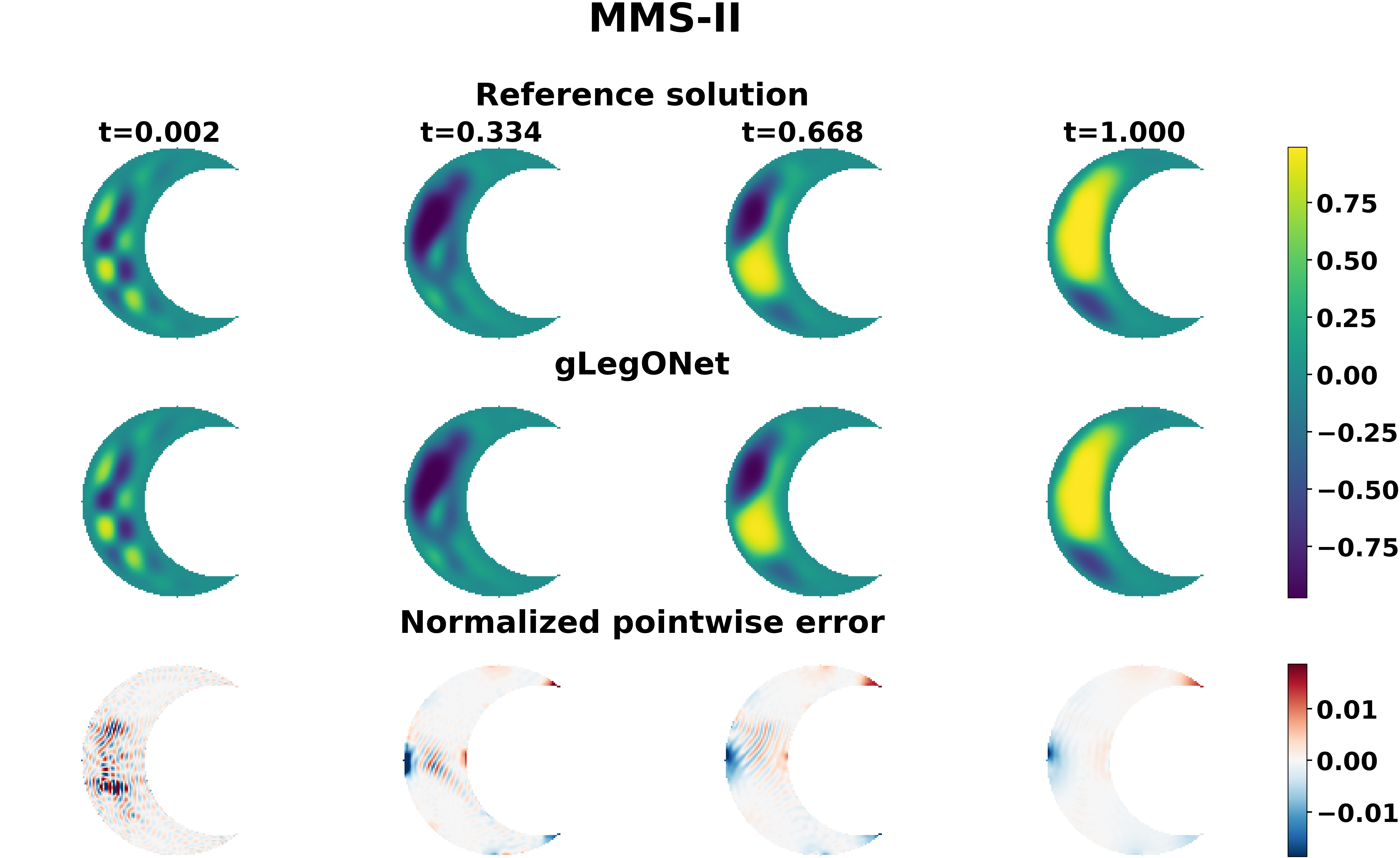}
    \end{minipage}%
    \begin{minipage}[c]{0.765\linewidth}
        \centering
        \includegraphics[
            width=\linewidth,
            height=0.125\textheight,
            keepaspectratio,
            trim=4 6 4 22,
            clip
        ]{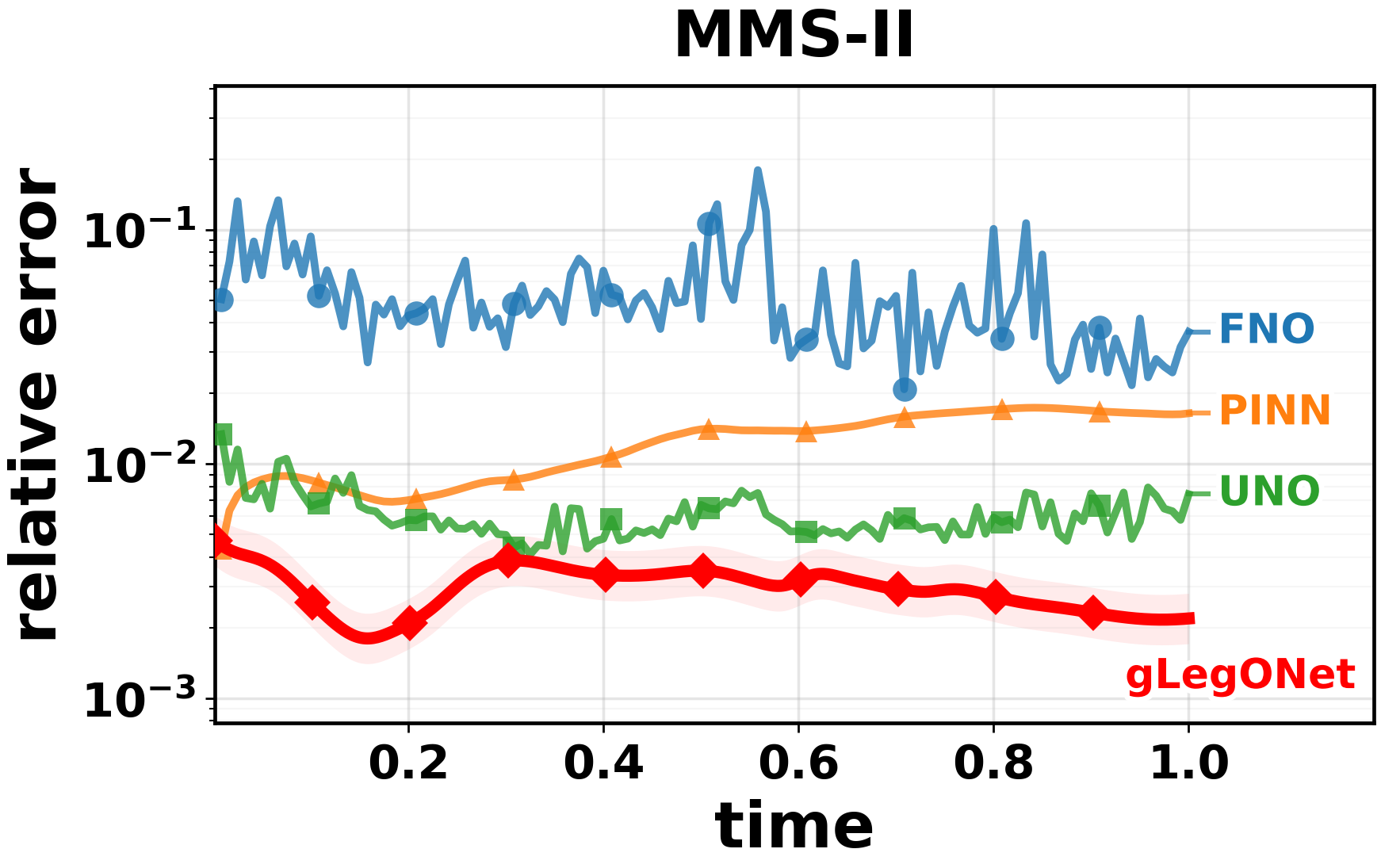}
    \end{minipage}
\end{minipage}

\vspace{0.9mm}

\noindent
\begin{minipage}[t]{0.496\textwidth}
    \centering
    {\scriptsize \textbf{(d)} MMS-III}\par
    \vspace{0.25mm}
    \begin{minipage}[c]{0.225\linewidth}
        \centering
        \includegraphics[
            height=0.125\textheight,
            keepaspectratio,
            trim=965 210 2 115,
            clip
        ]{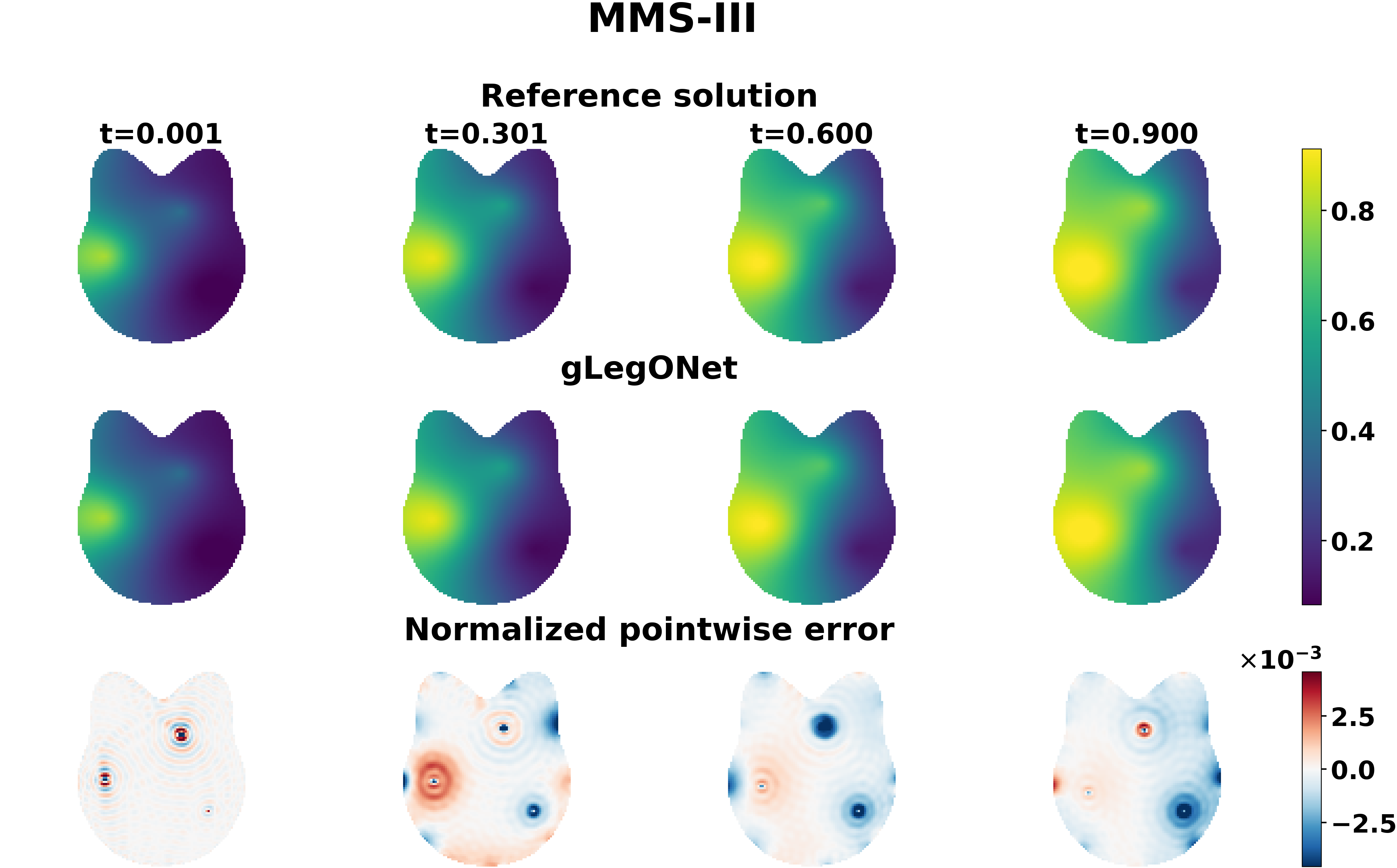}
    \end{minipage}%
    \begin{minipage}[c]{0.765\linewidth}
        \centering
        \includegraphics[
            width=\linewidth,
            height=0.125\textheight,
            keepaspectratio,
            trim=4 6 4 22,
            clip
        ]{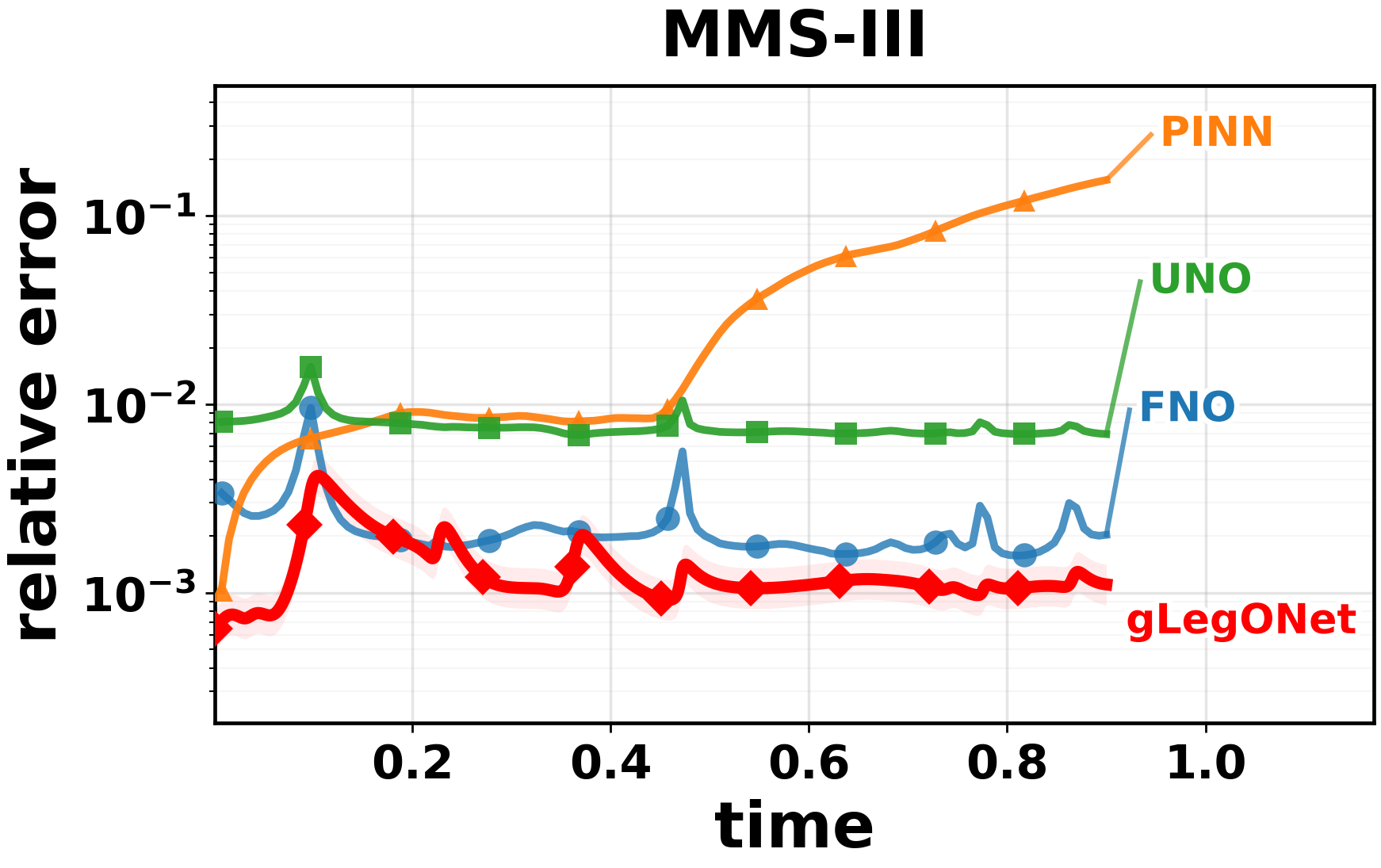}
    \end{minipage}
\end{minipage}%
\hspace{0.004\textwidth}%
\begin{minipage}[t]{0.496\textwidth}
    \centering
    {\scriptsize \textbf{(e)} MMS-V}\par
    \vspace{0.25mm}
    \begin{minipage}[c]{0.225\linewidth}
        \centering
        \includegraphics[
            height=0.121\textheight,
            keepaspectratio,
            trim=965 210 2 115,
            clip
        ]{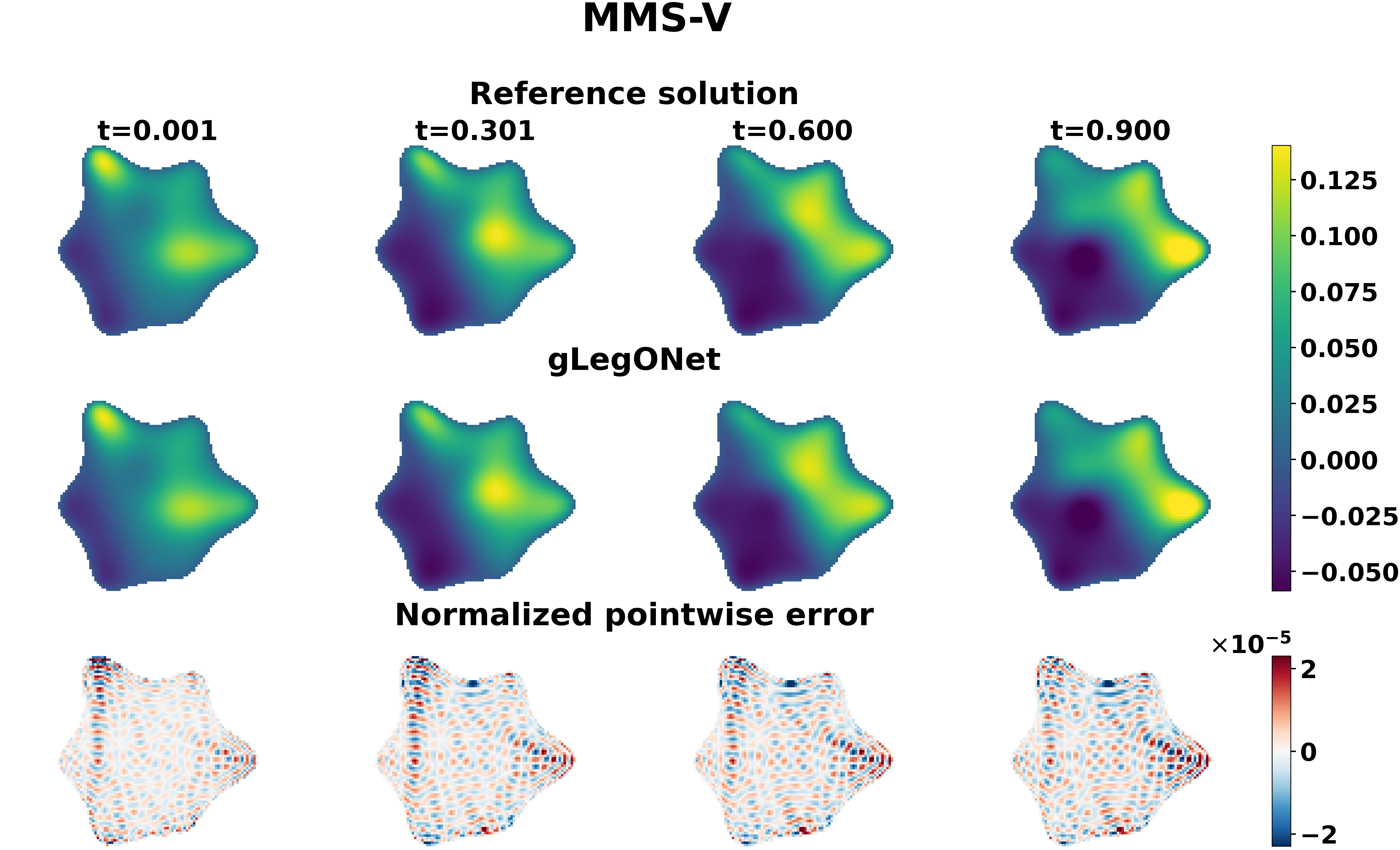}
    \end{minipage}%
    \begin{minipage}[c]{0.765\linewidth}
        \centering
        \includegraphics[
            width=\linewidth,
            height=0.125\textheight,
            keepaspectratio,
            trim=4 6 4 22,
            clip
        ]{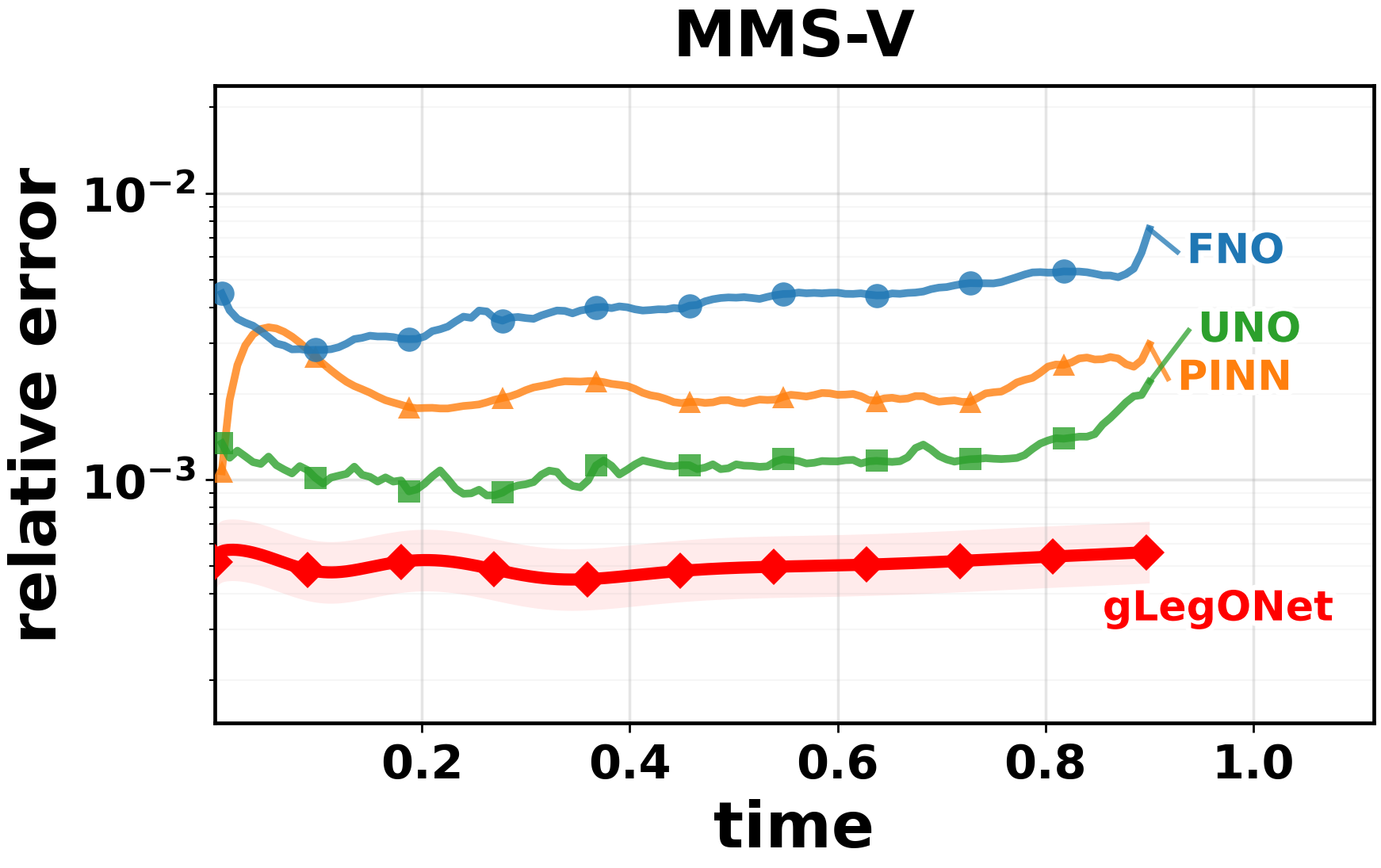}
    \end{minipage}
\end{minipage}

\vspace{0.9mm}

\begin{minipage}[t]{0.62\textwidth}
\centering
{\scriptsize \textbf{(f)} MMS-IV}\par
\vspace{0.25mm}
\begin{minipage}[c]{0.145\linewidth}
    \centering
    \includegraphics[
        height=0.125\textheight,
        keepaspectratio,
        trim=965 210 2 115,
        clip
    ]{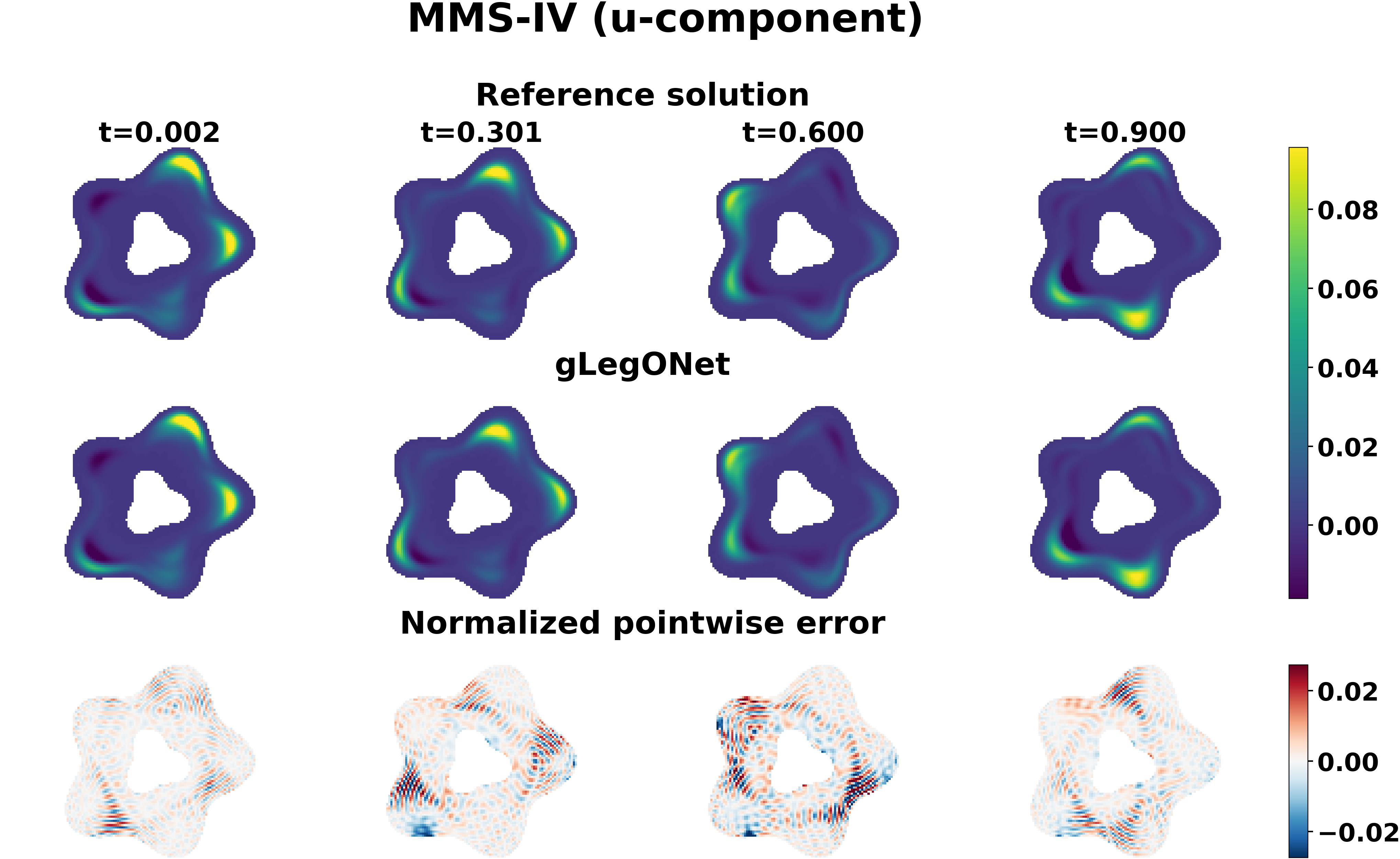}
\end{minipage}%
\begin{minipage}[c]{0.145\linewidth}
    \centering
    \includegraphics[
        height=0.125\textheight,
        keepaspectratio,
        trim=965 210 2 115,
        clip
    ]{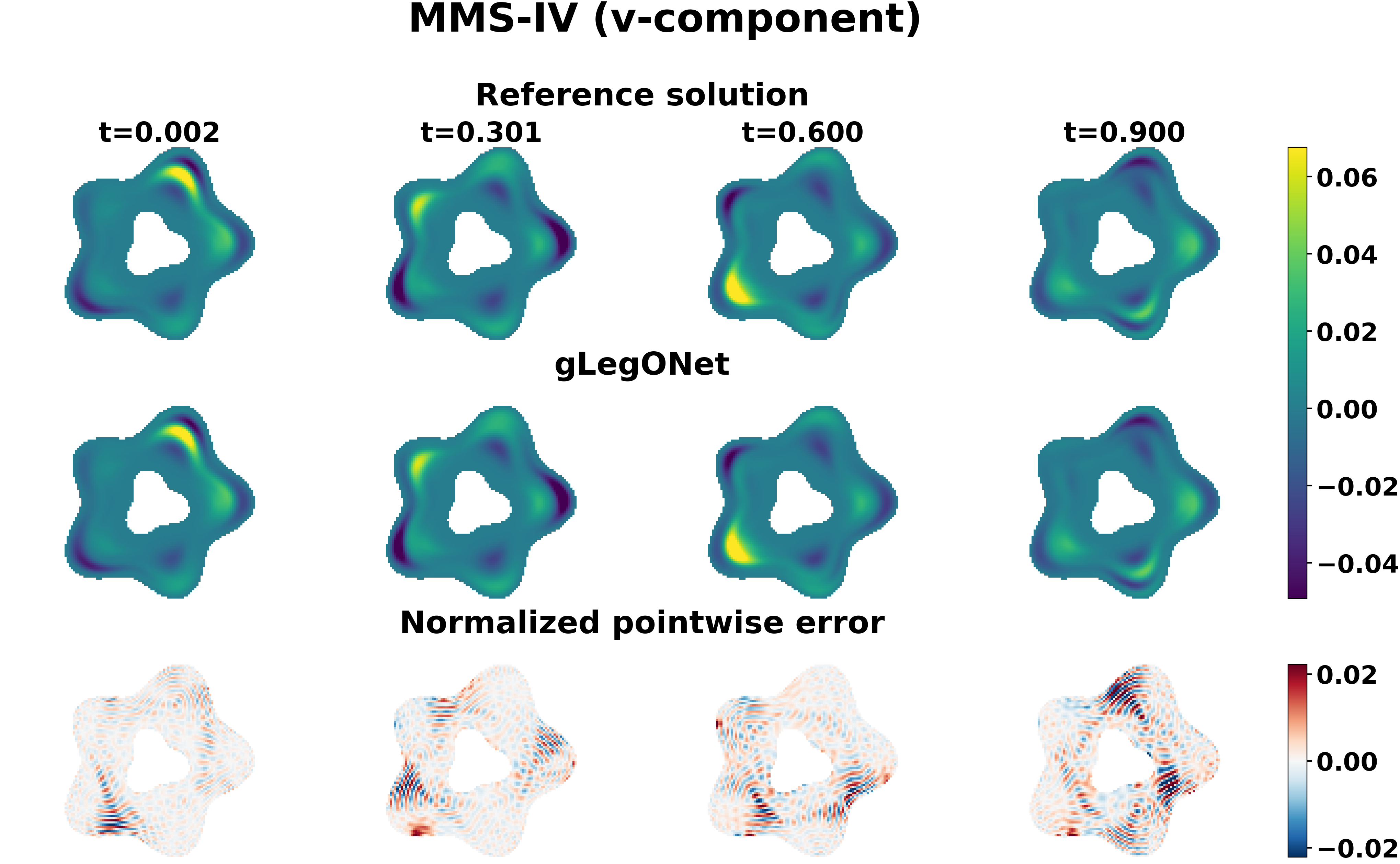}
\end{minipage}%
\begin{minipage}[c]{0.690\linewidth}
    \centering
    \includegraphics[
        width=\linewidth,
        height=0.125\textheight,
        keepaspectratio,
        trim=4 6 4 22,
        clip
    ]{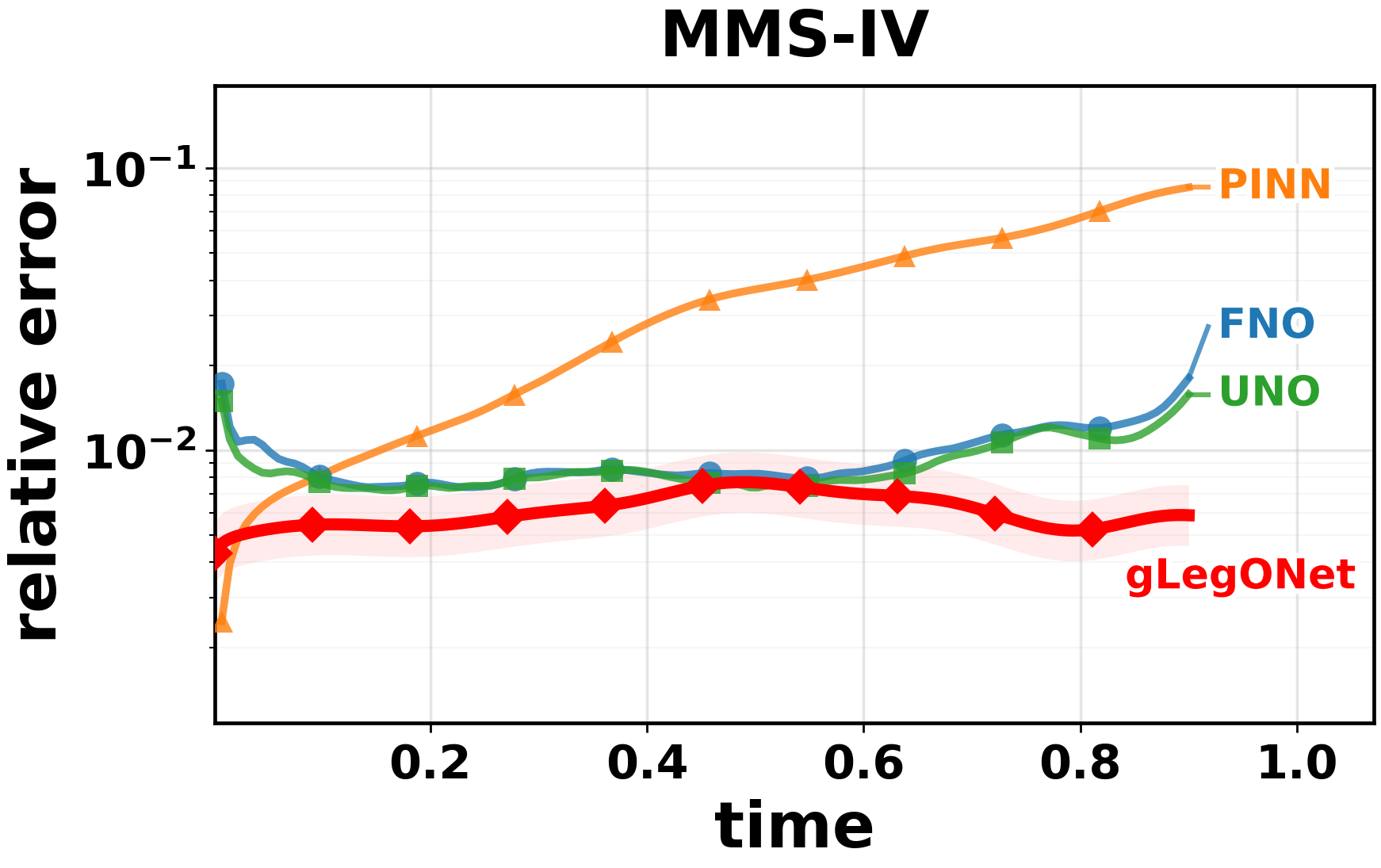}
\end{minipage}
\end{minipage}

\vspace{0.5mm}
\caption{
Manufactured-solution benchmarks on unseen embedded domains.  Panel (a)
summarizes the PDE, boundary condition, number of boundary samples \(N_b\), and
reduced dimension \(r\).  Panels (b)--(f) show rollout snapshots and relative \(L^2(\Omega)\) error
histories compared with PINN, FNO and UNO; in each snapshot block, the first row
is the reference solution and the second row is the gLegONet prediction.
MMS-IV shows separate final-time snapshots for the \(u\)- and \(v\)-components.
Full rollout panels and numerical tables are provided in Supplementary
Information.
}
\label{fig:mms}
\end{figure}

We first ask whether the same frozen library can be realized on unseen domains
and boundary operators without geometry-specific retraining.  Five
manufactured-solution benchmarks provide exact trajectories and therefore
reflect boundary realization and reduced-space rollout: Dirichlet
reaction--diffusion on a rosette, Neumann reaction--diffusion on a crescent,
nonhomogeneous Robin dynamics on a bunny-shaped domain, a coupled annular-star
system with mixed Dirichlet--Neumann constraints, and Burgers transport on a
pinwheel domain (Fig.~\ref{fig:mms}).  All five tests use the same \(K=22\)
library; only \(C\), \(\mathbf a_{\rm bc}\), \(M_\Omega\) and \(N_\Omega\) are
rebuilt.

PINN is optimized separately for each target, whereas FNO and UNO are trained
on supervised trajectories from nearby manufactured families.  Geometry-aware
LegONet uses neither target-specific neural training nor target-family
trajectories.  Across the five tests, it achieves \(10^{-4}\)--\(10^{-3}\)-level
final errors with runtimes of \(2.12\)--\(19.5\) s, compared with
hundreds to thousands of seconds for the baselines
(Extended Data Table~\ref{tab:mms_baseline_comparison}).  Homogeneous traces
remain near the algebraic tolerance, and nonhomogeneous Robin data are imposed
through an affine lift.  Boundary leakage is therefore suppressed below the
dominant interior-approximation, retained-rank and time-integration errors.
Full benchmark and baseline specifications are provided in Supplementary
Information.

\subsection*{Geometry-aware LegONet preserves constrained dynamics}

We next ask whether gLegONet preserves the structure of a nonlinear constrained system, rather than only matching field values.  The benchmark is the volume-constrained Allen--Cahn
equation on the disk
\(\Omega_{\rm disk}=\{\mathbf x:x^2+y^2\le0.4^2\}\), with \(\varepsilon^2=5\times10^{-3}\) and

\noindent final time \(T=4.0\):
\[
\left\{
\begin{aligned}
    u_t&=\varepsilon^2\Delta u+u-u^3-\lambda(t),
    &&\mathbf x\in\Omega_{\rm disk},\\
    \partial_{\mathbf n}u&=0,
    &&\mathbf x\in\partial\Omega_{\rm disk},
\end{aligned}
\right.
\qquad
\lambda(t)=
\frac{1}{|\Omega_{\rm disk}|}
\int_{\Omega_{\rm disk}}(u-u^3)\,d\mathbf x .
\]
Here \(\lambda(t)\) preserves the spatial mean, making this
a joint test of field accuracy, boundary enforcement and global conservation.
We use \(N_b=1600\) Neumann boundary samples and retain \(r=412\) boundary-adapted coordinates.
An independent \(P_1\) finite-element solution is used as reference.  We report
\(D_{\rm mass}^{\max}=\max_t|\bar u_h(t)-\bar u_h(0)|\), where
\(\bar u_h\) is the quadrature-weighted spatial mean.

\begin{figure}[H]
\centering
\scriptsize
\setlength{\abovecaptionskip}{2pt}
\setlength{\belowcaptionskip}{0pt}

\begin{minipage}[t]{0.96\textwidth}
    \centering
    \includegraphics[
        width=\linewidth,
        height=0.245\textheight,
        keepaspectratio
    ]{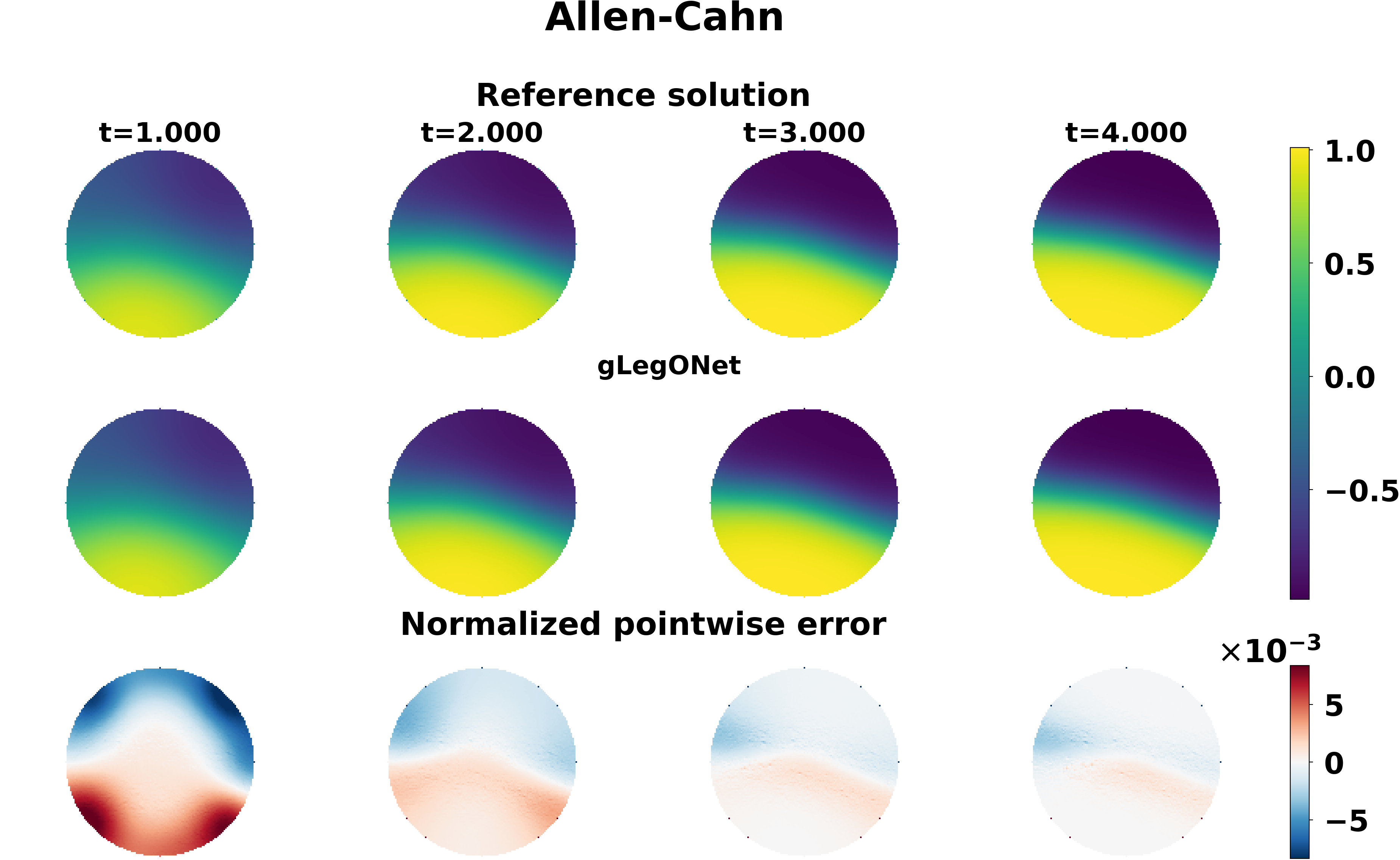}

    \vspace{0.4mm}
    {\scriptsize
    \textbf{(a)} FEM reference, gLegONet prediction and normalized
    pointwise error.}
\end{minipage}

\vspace{1.5mm}

\begin{minipage}[t]{0.48\textwidth}
    \centering
    \includegraphics[
        width=\linewidth,
        height=0.245\textheight,
        keepaspectratio
    ]{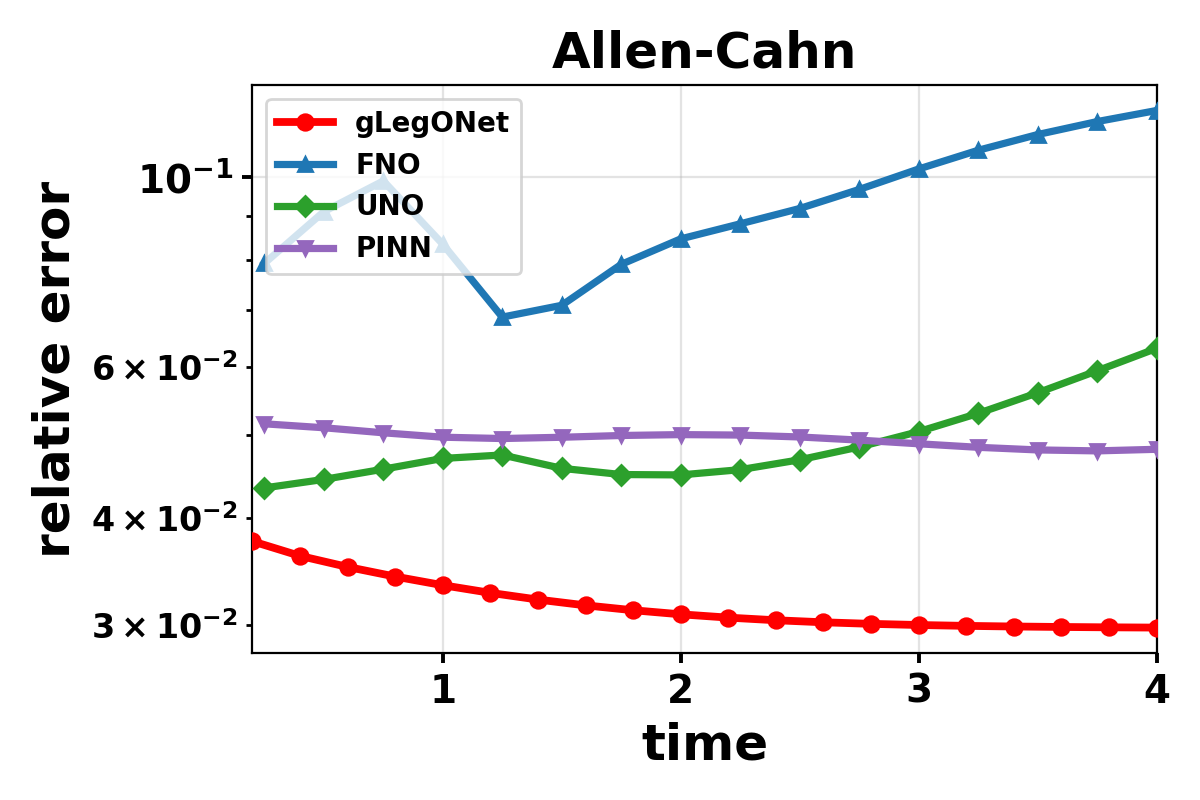}
    {\scriptsize
    \textbf{(b)} Relative \(L^2(\Omega_{\rm disk})\) error history.}
\end{minipage}
\hfill
\begin{minipage}[t]{0.48\textwidth}
    \centering
    \includegraphics[
        width=\linewidth,
        height=0.245\textheight,
        keepaspectratio
    ]{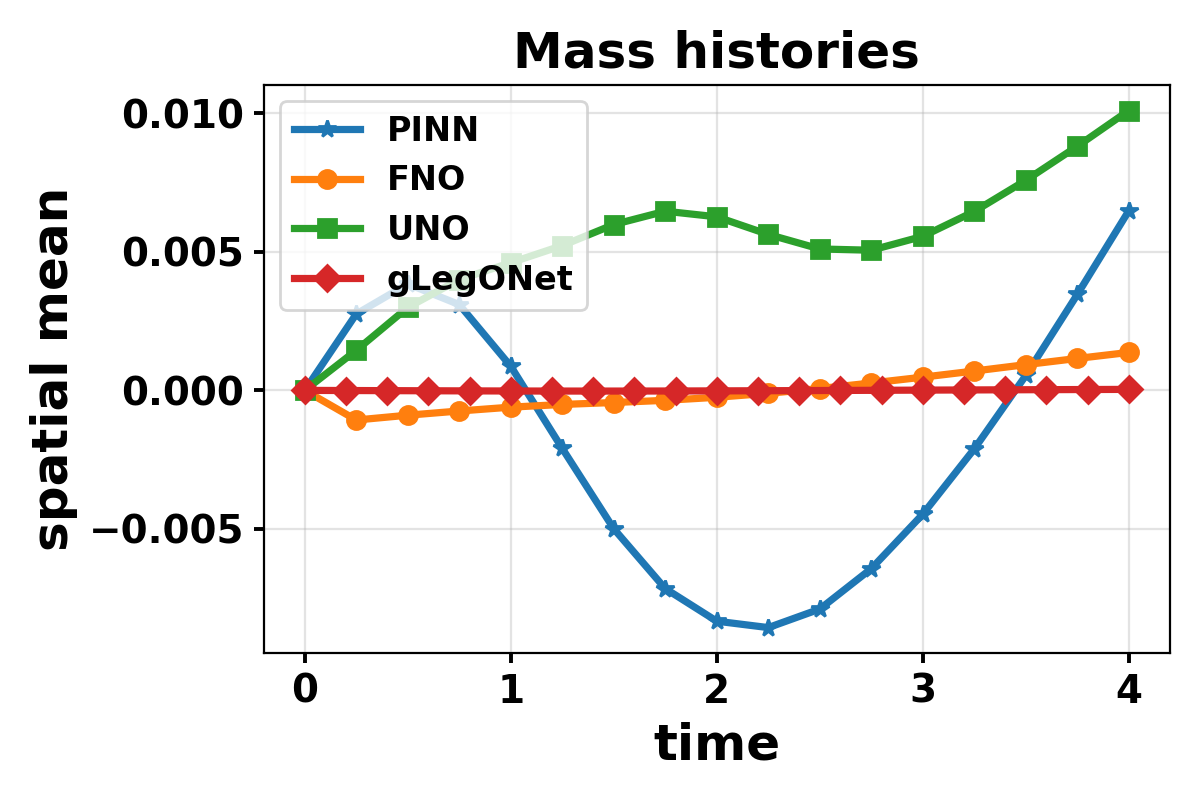}

    \vspace{0.3mm}
    {\scriptsize
    \textbf{(c)} Spatial mean history.}
\end{minipage}

\vspace{1.5mm}

\begin{minipage}[t]{0.48\textwidth}
    \centering
    \includegraphics[
        width=\linewidth,
        height=0.245\textheight,
        keepaspectratio
    ]{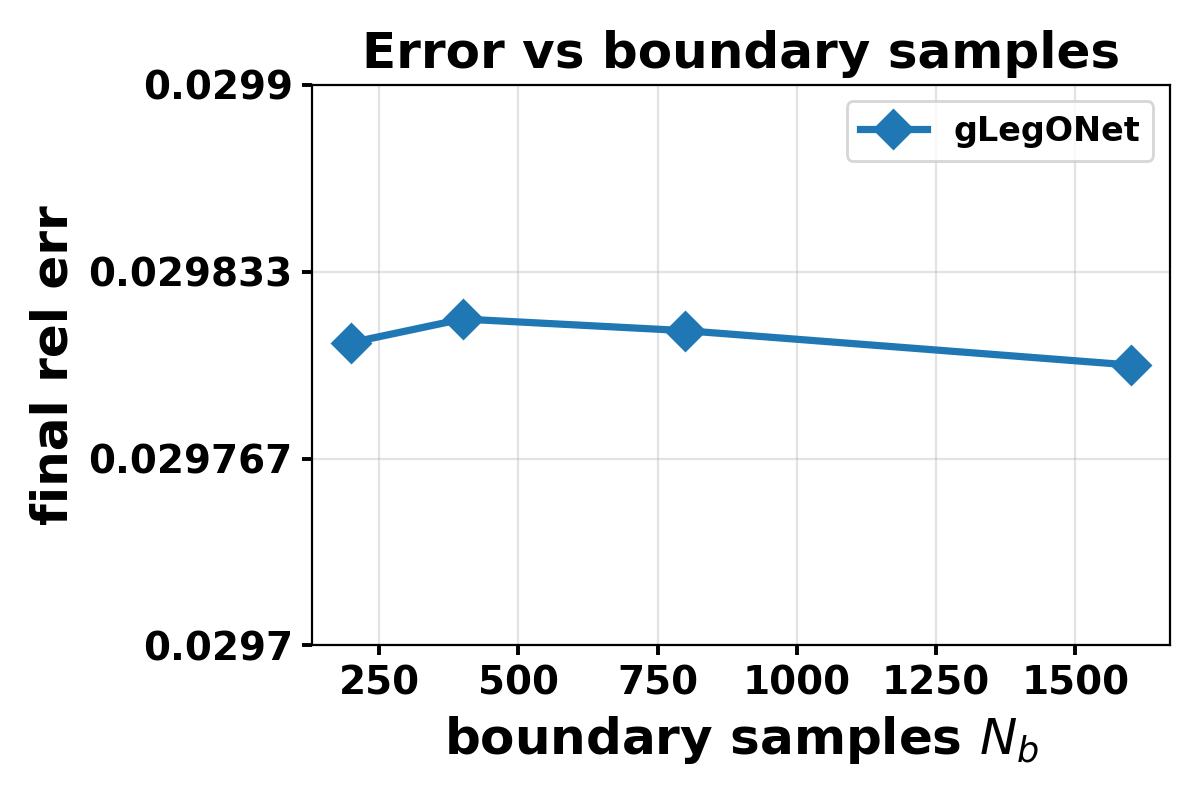}
    {\scriptsize
    \textbf{(d)} Sensitivity to the number of Neumann boundary samples
    \(N_b\).}
\end{minipage}
\hfill
\begin{minipage}[t]{0.48\textwidth}
    \centering
    \includegraphics[
        width=\linewidth,
        height=0.245\textheight,
        keepaspectratio
    ]{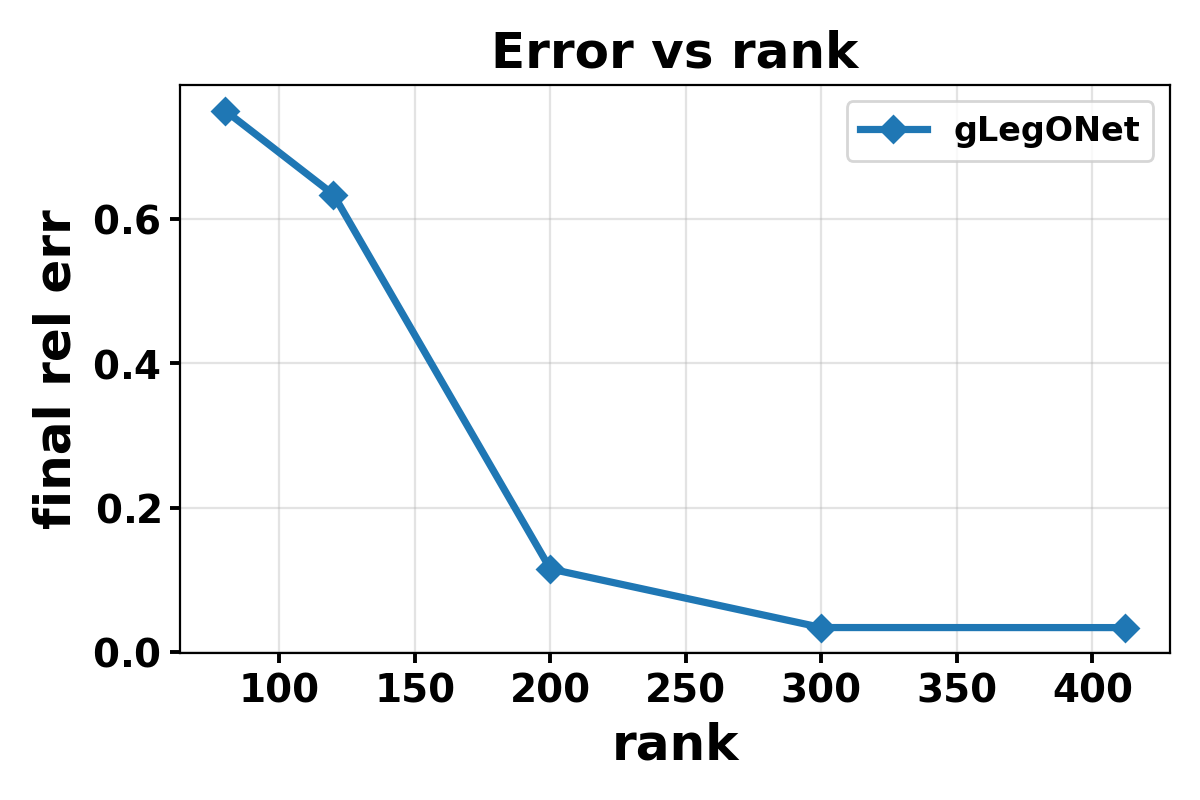}

    \vspace{0.3mm}
    {\scriptsize
    \textbf{(e)} Sensitivity to the retained reduced rank \(r\).}
\end{minipage}

\caption{
Boundary-adapted gLegONet rollout and ablations for volume-constrained
Allen--Cahn dynamics on the disk.  Panel (a) compares the finite-element
reference with the gLegONet prediction.  Panels (b,c) report the relative field error and
spatial mean, respectively.  Panels (d,e) vary the number of boundary samples
and the retained reduced rank.  The error is weakly sensitive to \(N_b\) once
the Neumann trace is resolved, but decreases sharply with \(r\) before
saturating, showing that reduced-space capacity rather than boundary sampling
is the dominant limitation.
}
\label{fig:allen_cahn_rollout_ablation}
\end{figure}

Fig.~\ref{fig:allen_cahn_rollout_ablation}(a,b) shows that gLegONet reproduces
the FEM two-phase state with final relative error \(2.98\times10^{-2}\).  The transferred coordinates also preserve the constrained structure:
the maximum mass drift is \(3.63\times10^{-5}\)
(Fig.~\ref{fig:allen_cahn_rollout_ablation}c), and the dense Neumann residual is
\(3.47\times10^{-9}\)
(Extended Data Table~\ref{tab:downstream_ac_burgers_baselines}).  PINN is the
closest baseline in final field error, but its error is \(1.6\) times larger
and its mass drift is \(2.4\times10^2\) times larger; FNO and UNO are less
accurate.  Panels (d,e) show that increasing \(r\) reduces the error until
saturation, whereas additional boundary samples have little effect once the
smooth Neumann trace is resolved.  Reduced-space capacity, rather than boundary
sampling, is therefore the dominant limitation.  Numerical settings are provided in
Supplementary Information.

\subsection*{Transferred reduced bases match geometry-specific bases}

We next test whether a reduced basis obtained by boundary transfer can match a basis constructed directly for the target geometry in a nonlinear transport problem.  On the embedded square \(\Omega=[-h,h]^2\), with \(h=0.55\), \(\nu=10^{-2}\) and final time \(T=0.5\), we solve
\[
\left\{
\begin{aligned}
    u_t+u u_x+v u_y&=\nu\Delta u,\\
    v_t+u v_x+v v_y&=\nu\Delta v,
\end{aligned}
\right.
\qquad
\mathbf x\in\Omega,
\qquad
u=v=0\ \text{on }\partial\Omega .
\]
This benchmark requires diffusion, directional transport and cross-advection
to be composed in one Dirichlet-compatible coordinate space.  The same boundary-adapted matrix \(N_\Omega\) is used for both velocity components. It is constructed from \(N_b=1000\) boundary samples and retains \(r=833\) coordinates. The frozen diffusion and transport blocks are reused through this interface, while cross-advection is assembled by analytic weak projection. An independent \(201\times201\) upwind finite-volume solution provides the reference.

\begin{figure}[H]
\centering
\scriptsize
\setlength{\abovecaptionskip}{2pt}
\setlength{\belowcaptionskip}{0pt}

\begin{minipage}[t]{0.49\textwidth}
    \vspace{0pt}
    \centering
    \includegraphics[
        width=\linewidth,
        height=0.219\textheight,
        keepaspectratio
    ]{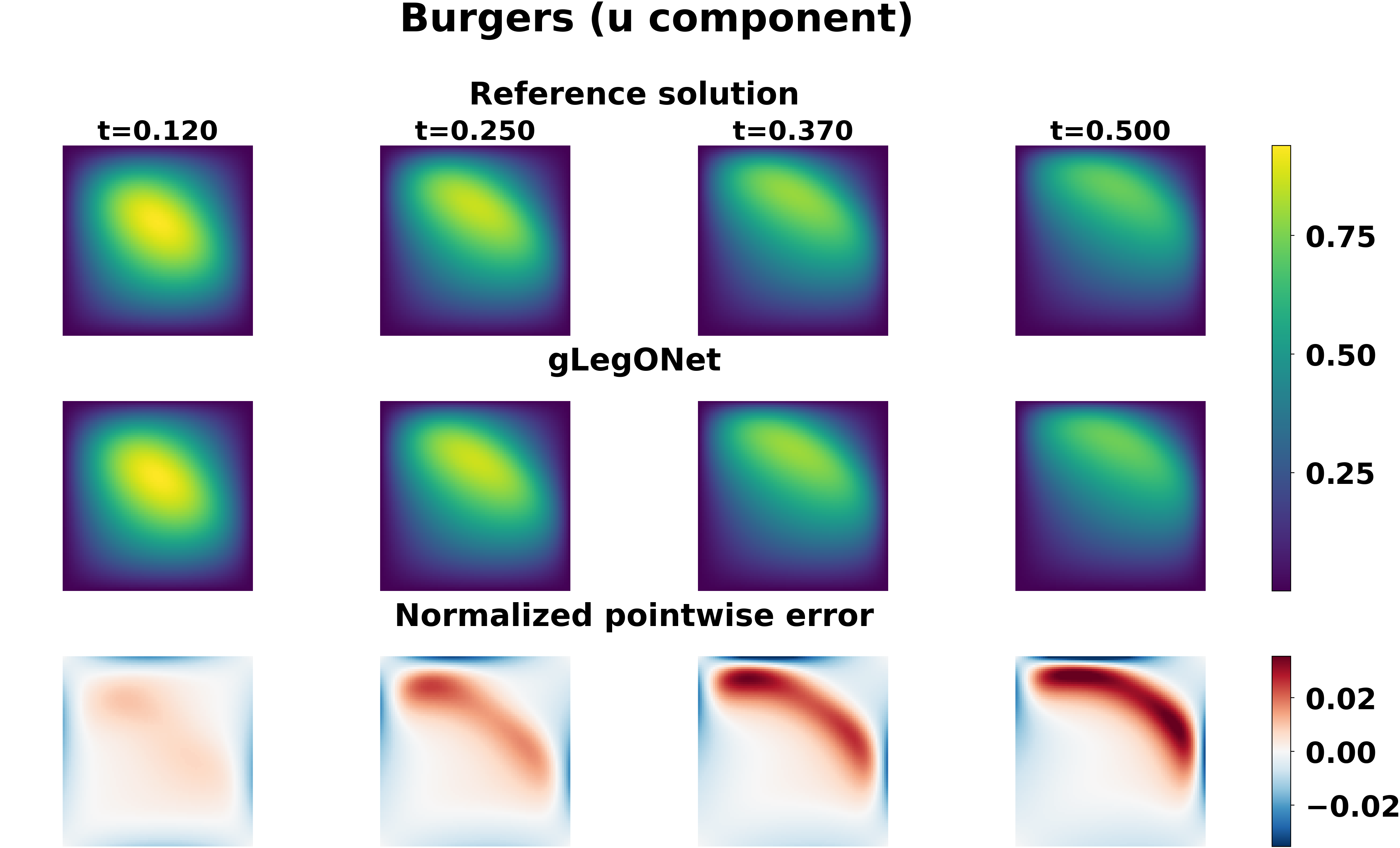}

    \vspace{0.3mm}
    {\scriptsize
    \textbf{(a)} Reference, gLegONet prediction and normalized
    pointwise error for \(u\).}
\end{minipage}
\hfill
\begin{minipage}[t]{0.49\textwidth}
    \vspace{0pt}
    \centering
    \includegraphics[
        width=\linewidth,
        height=0.219\textheight,
        keepaspectratio
    ]{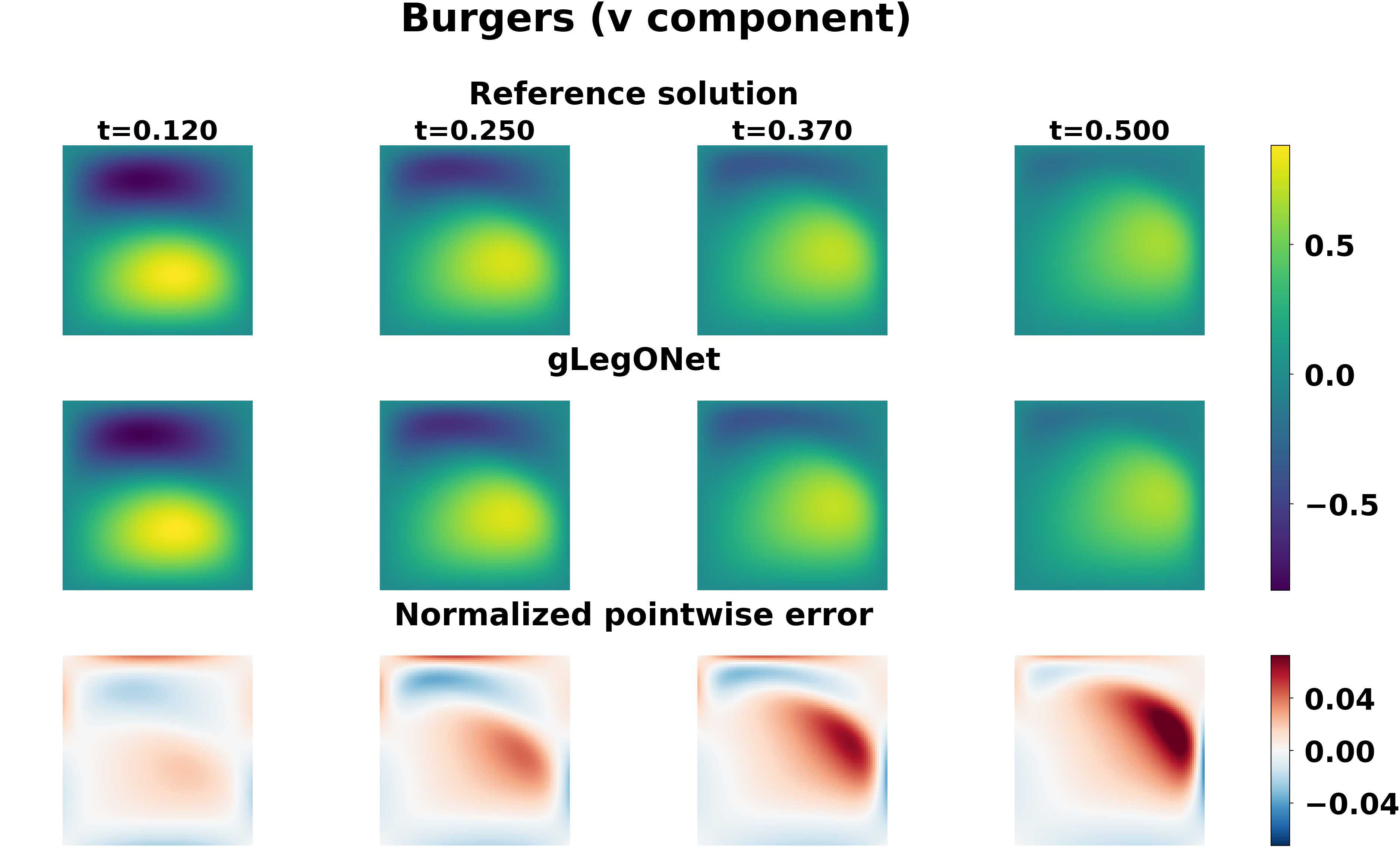}

    \vspace{0.3mm}
    {\scriptsize
    \textbf{(b)} Reference, gLegONet prediction and normalized
    pointwise error for \(v\).}
\end{minipage}

\vspace{1.5mm}

\begin{minipage}[t]{0.49\textwidth}
    \vspace{0pt}
    \centering
    \includegraphics[
        width=\linewidth,
        height=0.215\textheight,
        keepaspectratio
    ]{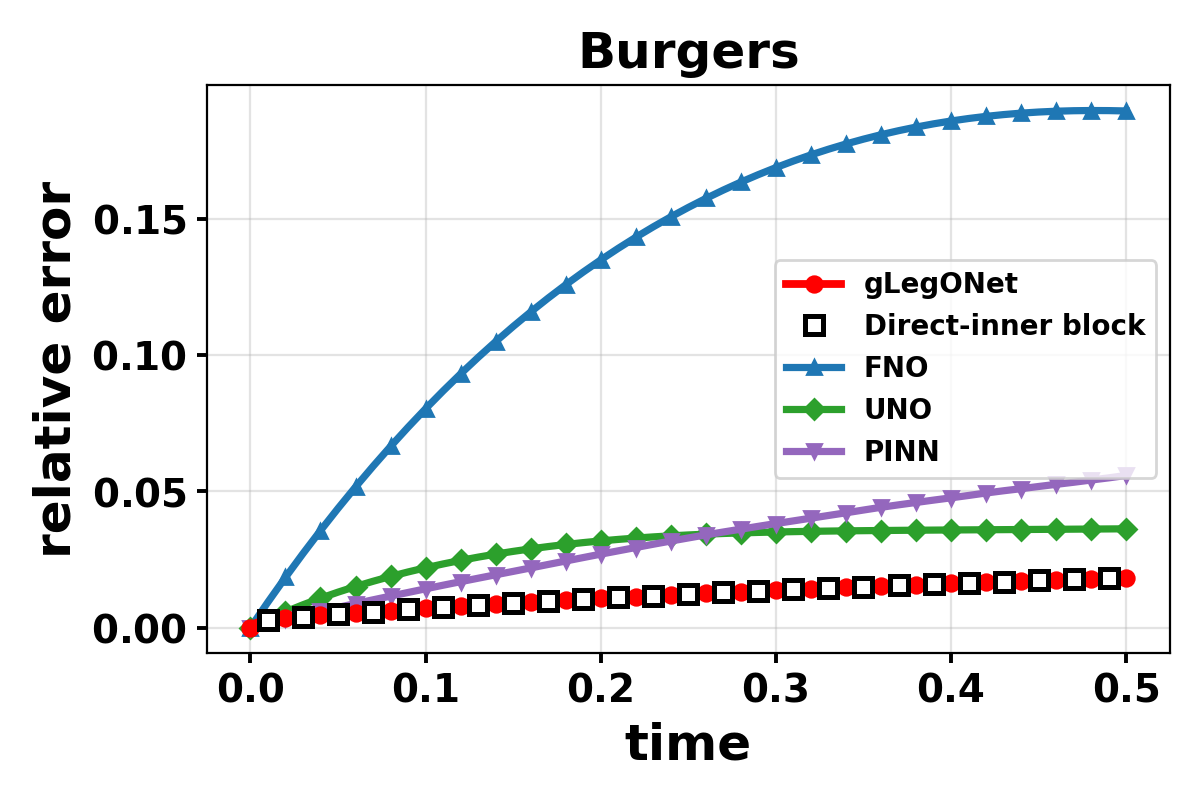}
    {\scriptsize
    \textbf{(c)} Relative \(L^2(\Omega)\) error history.}
\end{minipage}
\hfill
\begin{minipage}[t]{0.49\textwidth}
    \vspace{0pt}
    \centering
    \includegraphics[
        width=\linewidth,
        height=0.215\textheight,
        keepaspectratio
    ]{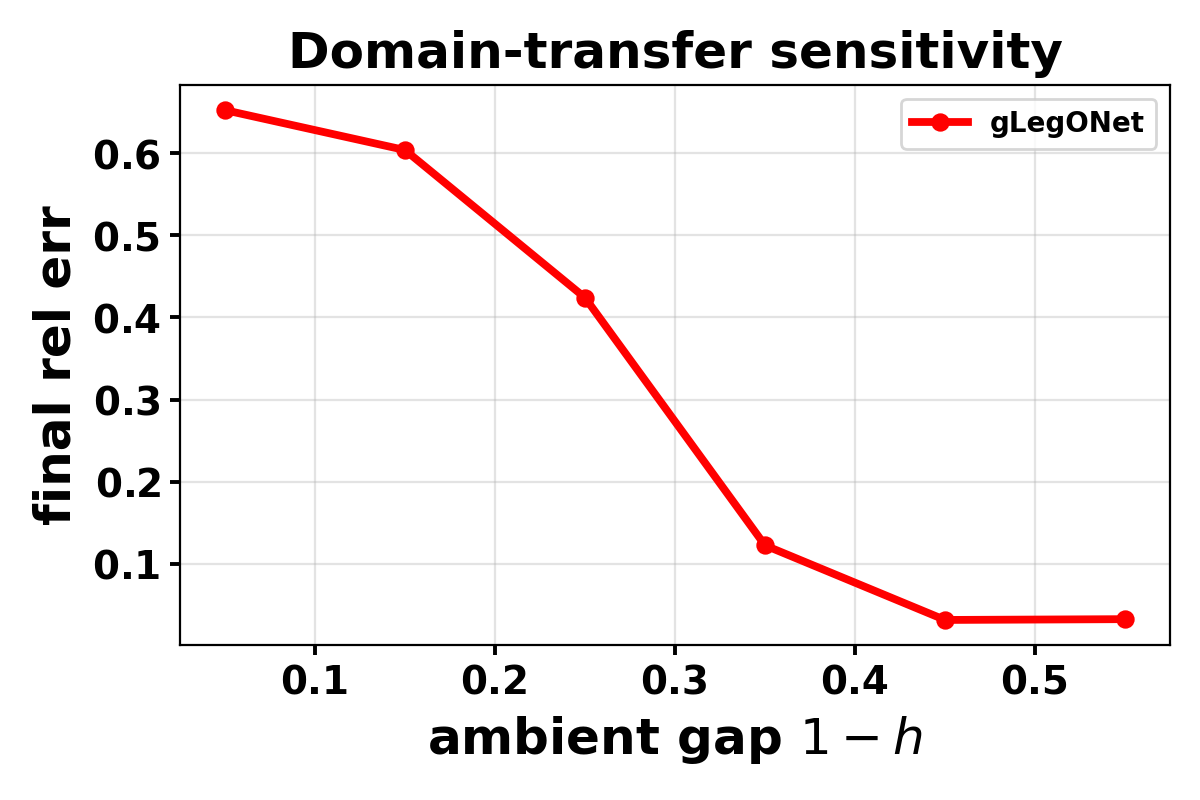}
    {\scriptsize
    \textbf{(d)} Sensitivity to the ambient-boundary gap \(1-h\).}
\end{minipage}

\vspace{1.5mm}

\begin{minipage}[t]{0.49\textwidth}
    \vspace{0pt}
    \centering
    \includegraphics[
        width=\linewidth,
        height=0.215\textheight,
        keepaspectratio
    ]{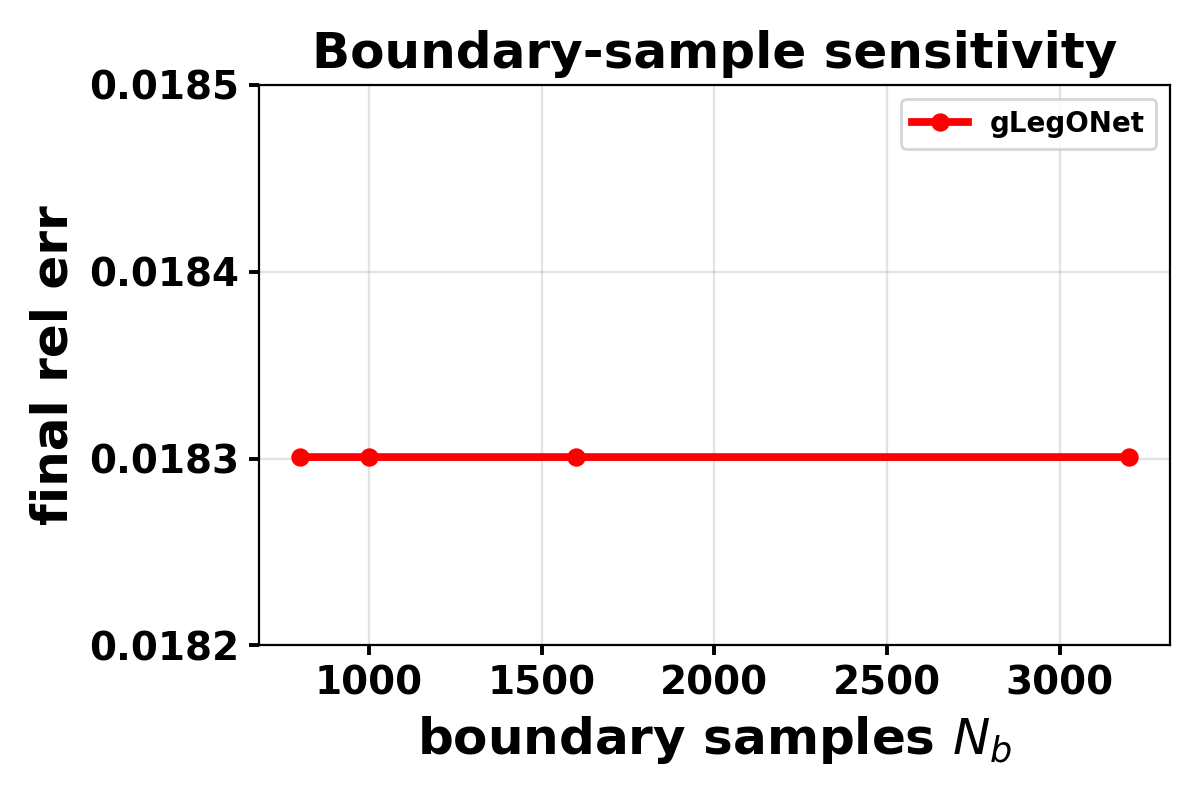}
    {\scriptsize
    \textbf{(e)} Sensitivity to the number of Dirichlet boundary samples
    \(N_b\).}
\end{minipage}
\hfill
\begin{minipage}[t]{0.49\textwidth}
    \vspace{0pt}
    \centering
    \includegraphics[
        width=\linewidth,
        height=0.215\textheight,
        keepaspectratio
    ]{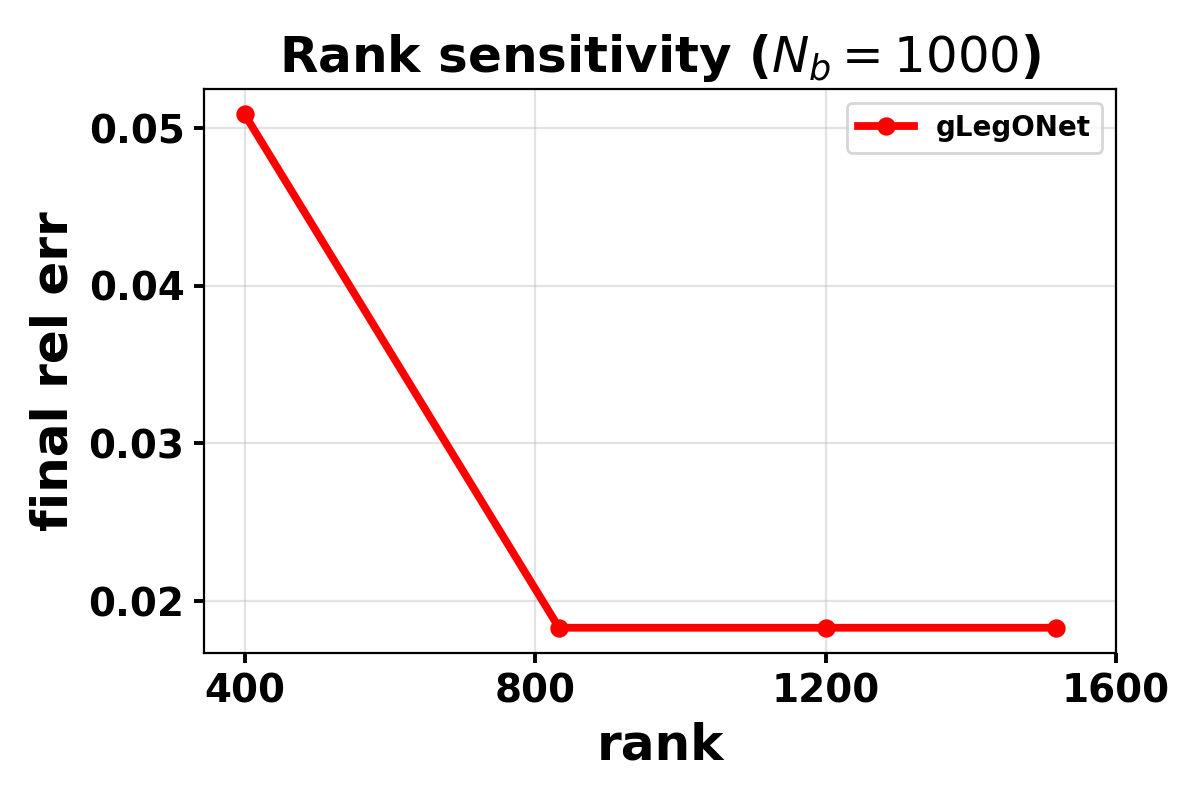}
    {\scriptsize
    \textbf{(f)} Sensitivity to the retained reduced rank \(r\).}
\end{minipage}

\caption{
gLegONet prediction for nonlinear vector Burgers dynamics.  Panels (a,b) compare
the reference and gLegONet rollouts for the \(u\)- and \(v\)-components.
Panel (c) compares the relative error histories with the neural and
geometry-specific baselines.  Panel (d) varies the ambient-boundary gap \(1-h\), while panels (e,f) vary the number of boundary samples and the retained
rank.  The error rises as the target domain approaches the ambient boundary,
changes little once the Dirichlet trace is sufficiently sampled, and decreases
with \(r\) until saturation.  The limiting factors are therefore ambient
resolution and reduced-space capacity rather than boundary leakage.
}
\label{fig:burgers_rollout_ablation}
\end{figure}

A direct inner-square block solver from our previous work \citep{zhang2026legonet} serves as a geometry-specific control, because its basis is built for the target boundary. gLegONet reaches the same accuracy without that tailored basis (Extended Data Table~\ref{tab:downstream_ac_burgers_baselines}): the blocks are reused as frozen ambient mechanisms, and only the boundary-adapted realization changes. As shown in Fig.~\ref{fig:burgers_rollout_ablation}(a--c), both velocity components track the reference, and the transferred and direct solvers have nearly overlapping error histories below the neural baselines. The ablations in Fig.~\ref{fig:burgers_rollout_ablation}(d--f) clarify the
expected operating regime: the embedded square should remain sufficiently
inside the ambient box, boundary-sample refinement has little effect once the
constraints are resolved, and increasing \(r\) improves accuracy until the
reduced space saturates. Thus the transfer limit is ambient-resolution and reduced-space capacity, not boundary leakage.  Numerical
settings are provided in Supplementary Information.

\subsection*{Frozen blocks compose into complex operator assemblies}

\begin{figure}[H]
\centering
\scriptsize
\setlength{\tabcolsep}{3pt}
\renewcommand{\arraystretch}{1.0}
\setlength{\abovecaptionskip}{2pt}
\setlength{\belowcaptionskip}{0pt}

\begin{minipage}[t]{0.95\textwidth}
    \centering
    \includegraphics[
        width=\linewidth,
        trim=0 0 17 0,
        clip
    ] {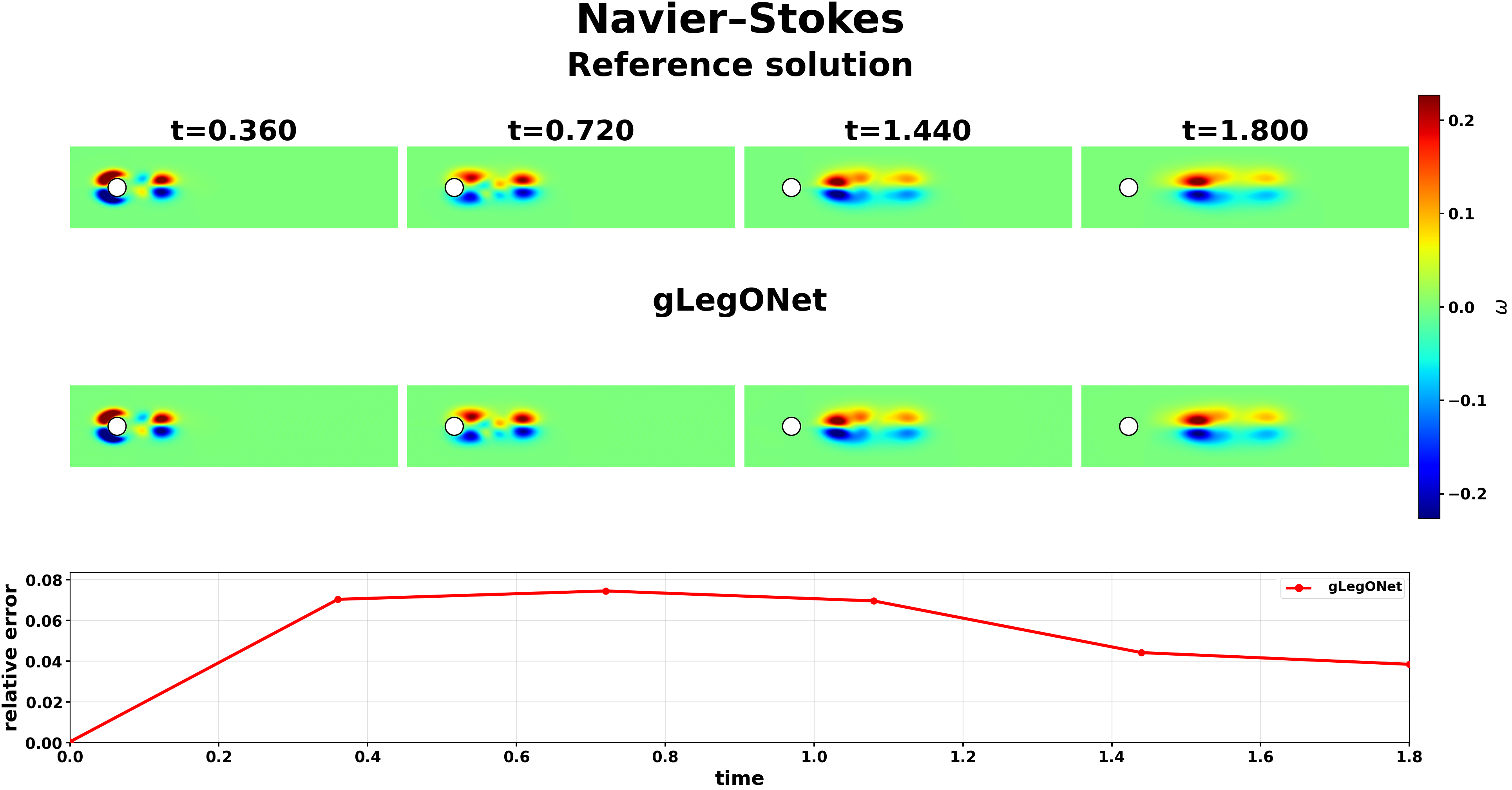}  
    \vspace{0.5mm}
    {\scriptsize \textbf{(a)} Navier--Stokes cylinder-wake vorticity rollout:
    masked-grid reference, gLegONet prediction, and
    relative-error history.}
\end{minipage}
\begin{minipage}[t]{0.93\textwidth}
    \centering
    \includegraphics[
        width=\linewidth,
        trim=0 0 15 0,
        clip
    ]{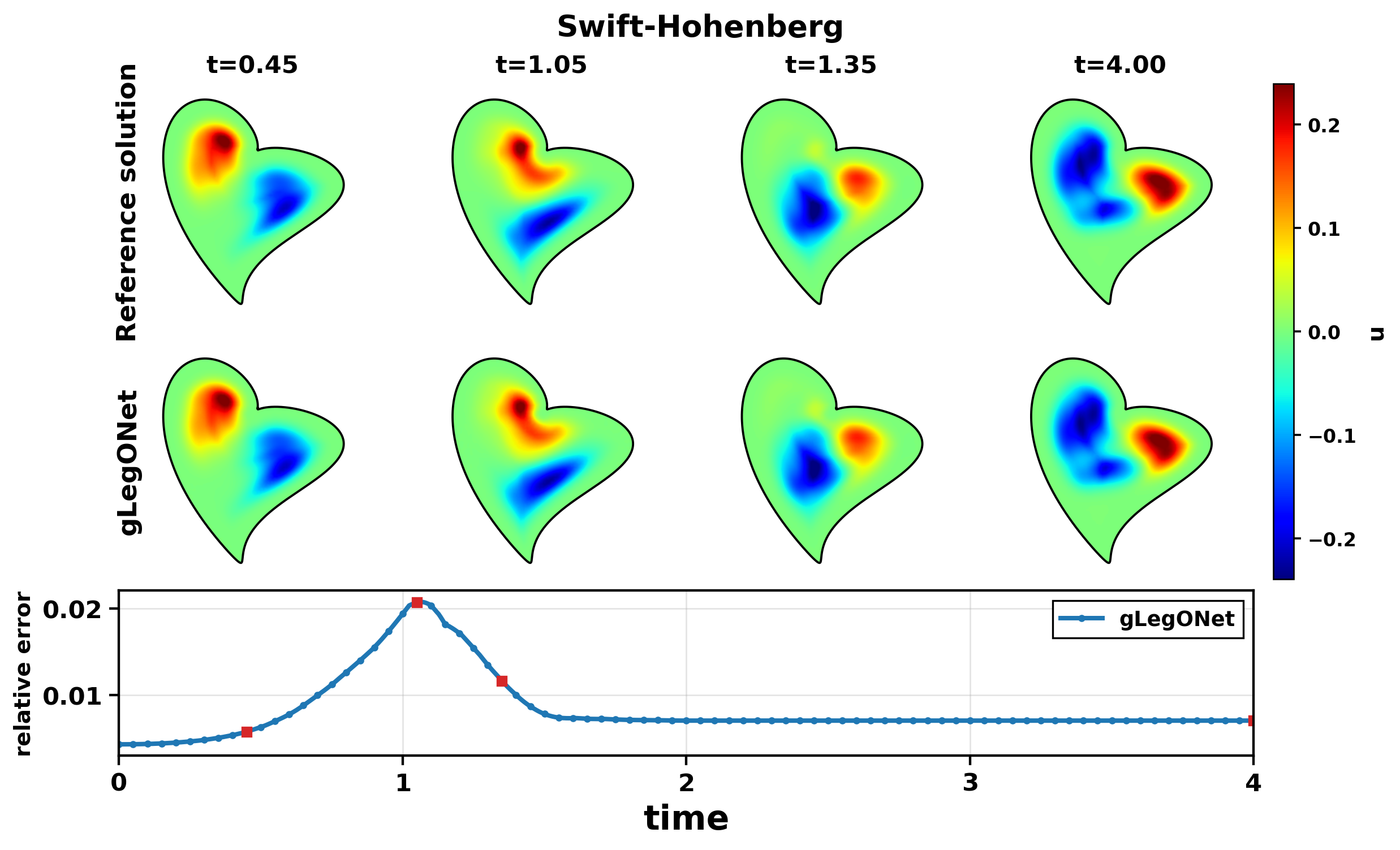}

    \vspace{0.5mm}
    {\scriptsize \textbf{(b)} Clamped Swift--Hohenberg rollout on a heart-shaped domain: manufactured reference, gLegONet prediction, and relative-error history.}
\end{minipage}

\caption{ gLegONet stress tests for boundary transfer and operator composition. Panel (a) shows vorticity transport past an embedded cylinder, where the solver must preserve homogeneous traces on both the channel wall and the obstacle while propagating a wake interaction. Panel (b) shows the clamped Swift--Hohenberg benchmark on a heart-shaped domain, where the solver simultaneously enforces \(u=0\) and \(\partial_{\mathbf n}u=0\), and composes the reusable second-order Laplace block into a fourth-order mechanism without training a fourth-order block. Detailed configurations, manufactured forcing, Fourier-resolution ablations and boundary diagnostics are given in Supplementary Information. }
\label{fig:ns_sh_rollout_main}
\end{figure}

We next test whether frozen second-order mechanisms can be assembled into more complex operator systems.  The first case is Navier–Stokes vorticity transport past an embedded cylinder:
\[
\left\{
\begin{aligned}
    \partial_t\omega+
    (\mathbf u_{\rm bg}+\nabla^\perp\psi)\cdot\nabla\omega
    &=\nu\Delta\omega+f,\\
    -\Delta\psi&=\omega,
\end{aligned}
\right.
\qquad
\omega=\psi=0\ \text{on }\partial\Omega.
\]
The embedded channel is $\Omega
    = [-0.95,0.95]\times[-0.2375,0.2375] \setminus B_{0.0523}\bigl((-0.6769,0)\bigr)$, where
\(B_R(\boldsymbol\xi_0)\) denotes the open disk of radius
\(R\) centred at \(\boldsymbol\xi_0\).
We use \(\nu=10^{-2}\) and \(T=1.8\).
This benchmark combines a multiply connected boundary,
an auxiliary Poisson solve for \(\psi\) and advective wake transport.  With \(K=45\) and
\(r=6099\), gLegONet achieves final relative error \(3.85\times10^{-2}\), and the dense boundary residual is
\(1.59\times10^{-12}\) (Fig.~\ref{fig:ns_sh_rollout_main}a).  The
\(K\)-ablation shows that boundary accuracy is retained while interior errors
decrease as the ambient representation is enriched
(Extended Data Fig.~\ref{fig:ns_smooth_left_k_ablation} and Extended Data
Table~\ref{tab:ns_smooth_left_ablation}).

The second case is the clamped Swift--Hohenberg equation on a
heart-shaped embedded domain,
\[
    \partial_t u
    =
    -\gamma(\Delta+k_0^2I)^2u+\mu u-u^3+f,
    \qquad
    u=\partial_{\mathbf n}u=0
    \ \text{on }\partial\Omega_{\rm heart}.
\]
We use \(\gamma=4.0\times10^{-3}\), \(k_0=5.3\), \(\mu=0.05\),
\(T=4.0\), \(\Delta t=10^{-3}\), and a \(K=30\) boundary-adapted space with
retained dimension \(r=795\).
No fourth-order neural block is trained.  Instead, $(\Delta+k_0^2I)^2$ is assembled by composing the pretrained Laplace block
inside the boundary-adapted space.
This controlled test prescribes a smooth manufactured path.
Fig.~\ref{fig:ns_sh_rollout_main}(b) shows that the composed fourth-order
rollout tracks the reference evolution with sub-\(10^{-2}\) final relative
error and peak error below \(2.1\times10^{-2}\).  Both clamped traces are
retained, with dense value- and normal-trace residuals below \(2\times10^{-8}\).
The Fourier-resolution ablation in Extended Data
Table~\ref{tab:heart_sh_k_ablation} confirms that increasing \(K\) primarily
improves interior resolution rather than boundary realization.

Together, these tests show that the same
interface supports multiply connected domains, auxiliary solves, multiple
boundary traces and higher-order operator composition.  Full configurations
are provided in Supplementary Information.

\subsection*{Sparse physical-law discovery on unseen domains}

Finally, we test whether gLegONet can use the same boundary-transfer layer for inverse discovery. To isolate the effect of the candidate mechanisms from the size of the reduced representation, all identification experiments use the same rank-\(80\) homogeneous Dirichlet boundary-adapted space on each unseen domain. The candidate library is
\[
    \mathcal D_\theta^\Omega
    =
    \{
        \mathbf q_\Delta^\Omega,
        \mathbf q_x^\Omega,
        \mathbf q_y^\Omega,
        \mathbf q_u^\Omega,
        \mathbf q_{u^2}^\Omega,
        \mathbf q_{u^3}^\Omega
    \},
\]

\noindent combining realized frozen blocks with analytic polynomial
maps. No neural parameter is updated during identification; only the sparse coefficients multiplying these realized mechanisms are fitted from short-time observations.

Identification is deliberately data limited. We use \(N_{\rm traj}=16\) independent short trajectories over \(T_{\rm id}=0.012\), sampled every \(\Delta t_{\rm obs}=5\times10^{-5}\), from either \(N_{\rm obs}=1000\) or \(240\) fixed random interior sensors. The recovered law is then rolled out to \(T_{\rm roll}=1.0\), more than \(80\) times longer than the identification window. The noise level \(\sigma\) denotes relative additive observation noise.

The peanut dumbbell tests
reaction--diffusion through a narrow bottleneck, with true coefficient vector
\[
    \mathbf c_{\rm peanut}^{\star}
    =
    (0.04,0,0,0.06,0,-0.35).
\]
The oblique sinusoidal channel combines nonlinear transport, reaction terms and
a thin anisotropic geometry, with
\[
    \mathbf c_{\rm channel}^{\star}
    =
    (0.02,-0.06,0.09,-0.03,0.025,-0.30).
\]
With \(N_{\rm obs}=1000\) and no observation noise, gLegONet recovers the
peanut law almost exactly, giving
\(e_{\rm coef}=6.89\times10^{-5}\) and
\(e_{\rm roll}=3.17\times10^{-5}\).  Even at
\(\sigma=10^{-2}\), the long-time rollout remains stable with final error
\(4.53\times10^{-2}\).  The channel is more demanding because the active law
contains both directional transport terms. Nevertheless, at zero noise
gLegONet gives \(e_{\rm coef}=1.25\times10^{-2}\) and
\(e_{\rm roll}=3.44\times10^{-2}\), and remains predictive under
low-to-moderate noise before degrading at the highest noise level.

The ablations show that the reduced boundary-adapted representation is
essential.  Replacing it by a rank-matched ambient Fourier basis fails: coefficient errors exceed \(10^2\) and the
rollout errors are near one on both domains.  Sparse finite differences also
fail under the same sensor budget, because scattered observations do not
provide stable derivative features. The sparse FD coefficient error on the
peanut domain reaches \(9.33\times10^{2}\), and the channel rollout becomes
unstable.  Dense finite differences are competitive only when given the full
masked grid, which is a much stronger observation setting than the sparse
sensor data used by gLegONet.  Reducing the sensor count to \(N_{\rm obs}=240\) still gives a predictive peanut rollout
(\(e_{\rm roll}=1.83\times10^{-4}\)), while the channel becomes less accurate
but remains finite.

\begin{figure}[H]
\centering
\footnotesize
\setlength{\tabcolsep}{4.3pt}
\renewcommand{\arraystretch}{1.05}
\setlength{\abovecaptionskip}{2pt}
\setlength{\belowcaptionskip}{0pt}

\begin{minipage}[t]{0.995\textwidth}
\centering
\resizebox{\linewidth}{!}{%
\begin{tabular}{llcccc}
\toprule
Setting
& Data
& Peanut \(e_{\rm coef}\)
& Peanut \(e_{\rm roll}\)
& Channel \(e_{\rm coef}\)
& Channel \(e_{\rm roll}\) \\
\midrule
\multicolumn{6}{l}{\textit{Noise robustness, gLegONet with \(N_{\rm obs}=1000\)}}\\
\midrule
gLegONet & \(\sigma=0\)
& \(6.89{\times}10^{-5}\pm5.02{\times}10^{-5}\)
& \(3.17{\times}10^{-5}\pm1.14{\times}10^{-5}\)
& \(1.25{\times}10^{-2}\pm3.74{\times}10^{-4}\)
& \(3.44{\times}10^{-2}\pm3.66{\times}10^{-4}\) \\

gLegONet & \(\sigma=10^{-4}\)
& \(3.90{\times}10^{-3}\pm3.90{\times}10^{-3}\)
& \(4.68{\times}10^{-4}\pm1.52{\times}10^{-4}\)
& \(1.31{\times}10^{-2}\pm6.55{\times}10^{-4}\)
& \(3.17{\times}10^{-2}\pm5.10{\times}10^{-3}\) \\

gLegONet & \(\sigma=10^{-3}\)
& \(3.96{\times}10^{-2}\pm3.95{\times}10^{-2}\)
& \(4.58{\times}10^{-3}\pm1.59{\times}10^{-3}\)
& \(2.23{\times}10^{-2}\pm2.39{\times}10^{-2}\)
& \(6.78{\times}10^{-2}\pm1.97{\times}10^{-2}\) \\

gLegONet & \(\sigma=10^{-2}\)
& \(4.16{\times}10^{-1}\pm4.01{\times}10^{-1}\)
& \(4.53{\times}10^{-2}\pm1.38{\times}10^{-2}\)
& \(1.72{\times}10^{-1}\pm2.35{\times}10^{-1}\)
& \(2.84{\times}10^{-1}\pm1.62{\times}10^{-1}\) \\

\midrule
\multicolumn{6}{l}{\textit{Zero-noise sparsity and basis ablation}}\\
\midrule
gLegONet & \(N_{\rm obs}=240\)
& \(6.20{\times}10^{-4}\pm4.44{\times}10^{-4}\)
& \(1.83{\times}10^{-4}\pm5.64{\times}10^{-5}\)
& \(5.04{\times}10^{-2}\pm2.66{\times}10^{-2}\)
& \(2.69{\times}10^{-1}\pm1.30{\times}10^{-1}\) \\

Ambient Fourier & \(N_{\rm obs}=1000\)
& \(6.42{\times}10^{2}\pm3.63{\times}10^{2}\)
& \(1.0174\pm6.58{\times}10^{-1}\)
& \(1.77{\times}10^{2}\pm8.11{\times}10^{1}\)
& \(9.98{\times}10^{-1}\pm5.01{\times}10^{-1}\) \\

\midrule
\multicolumn{6}{l}{\textit{Zero-noise finite-difference cost study}}\\
\midrule
FD dense & full masked grid
& \(3.12{\times}10^{-2}\pm3.26{\times}10^{-2}\)
& \(3.05{\times}10^{-4}\pm1.04{\times}10^{-4}\)
& \(8.02{\times}10^{-3}\pm4.60{\times}10^{-4}\)
& \(8.39{\times}10^{-3}\pm1.24{\times}10^{-3}\) \\

FD sparse & \(N_{\rm obs}=1000\)
& \(9.33{\times}10^{2}\pm9.77{\times}10^{2}\)
& \(7.57{\times}10^{-1}\pm3.24{\times}10^{-1}\)
& \(4.77{\times}10^{2}\pm2.67{\times}10^{2}\)
& \(\mathrm{unstable}\) \\
\bottomrule
\end{tabular}%
}

\vspace{0.5mm}
{\scriptsize \textbf{(a)} Boundary-aware law identification on unseen domains.}
\end{minipage}

\vspace{1.8mm}

\begin{minipage}[t]{0.48\textwidth}
\centering
\includegraphics[width=\linewidth]{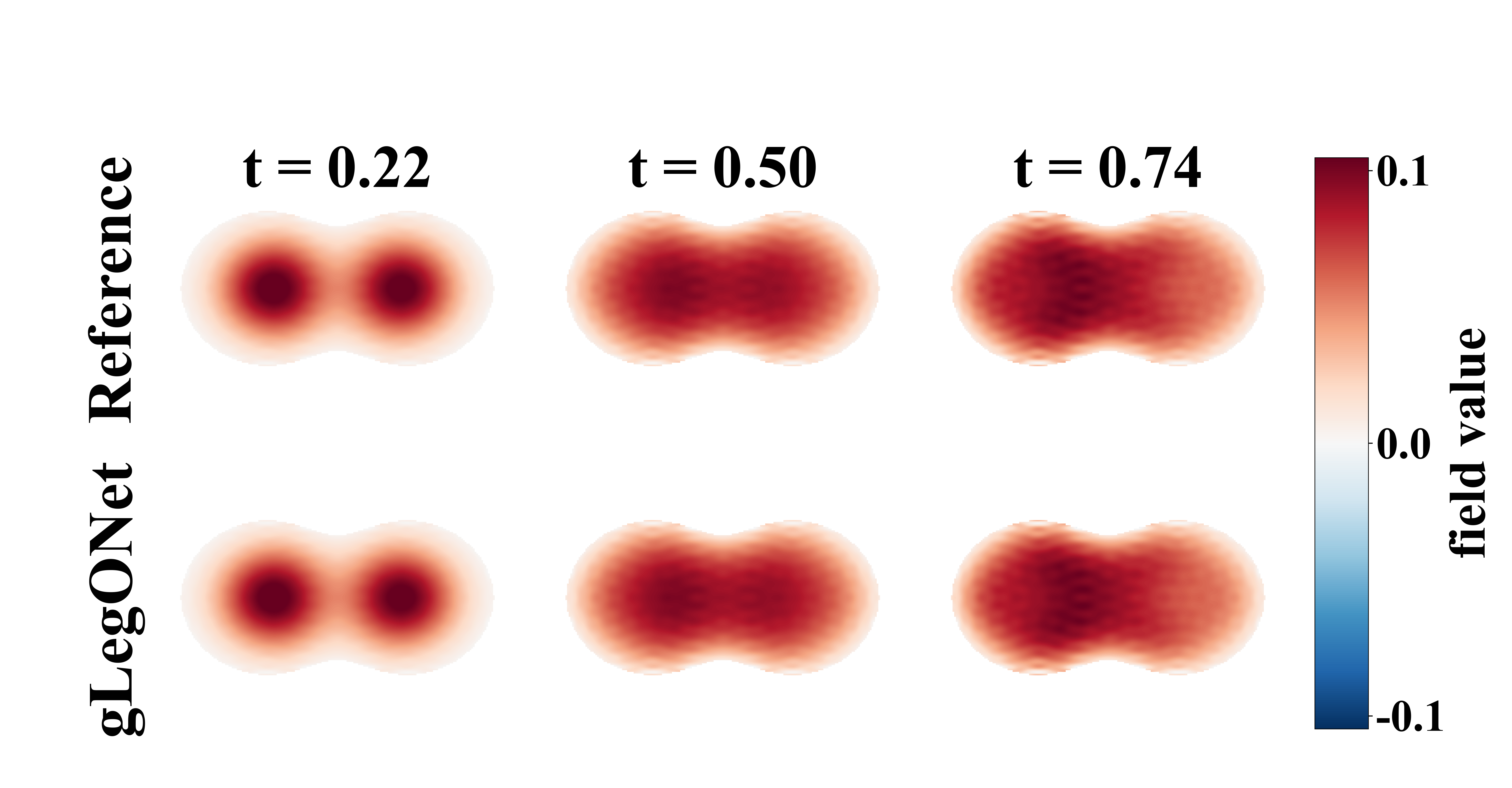}

\vspace{0.5mm}
{\scriptsize
\textbf{(b)} Peanut dumbbell: reference and identified boundary-adapted
rollout.}
\end{minipage}
\begin{minipage}[t]{0.48\textwidth}
\centering
\includegraphics[width=\linewidth]{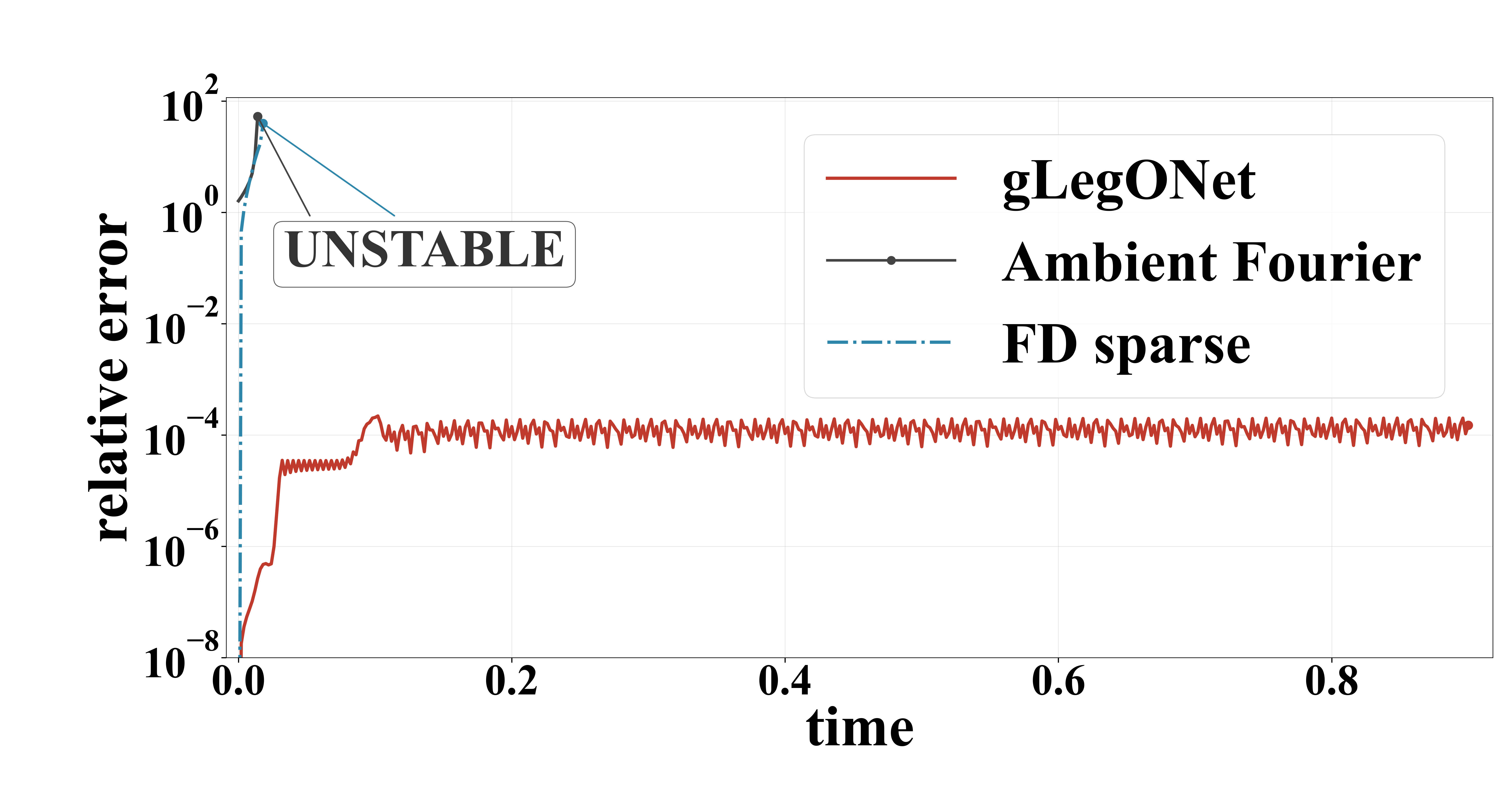}

\vspace{0.5mm}
{\scriptsize
\textbf{(c)} Peanut dumbbell: relative \(L^2(\Omega)\) error history.}
\end{minipage}

\begin{minipage}[t]{0.48\textwidth}
\centering
\includegraphics[width=\linewidth]{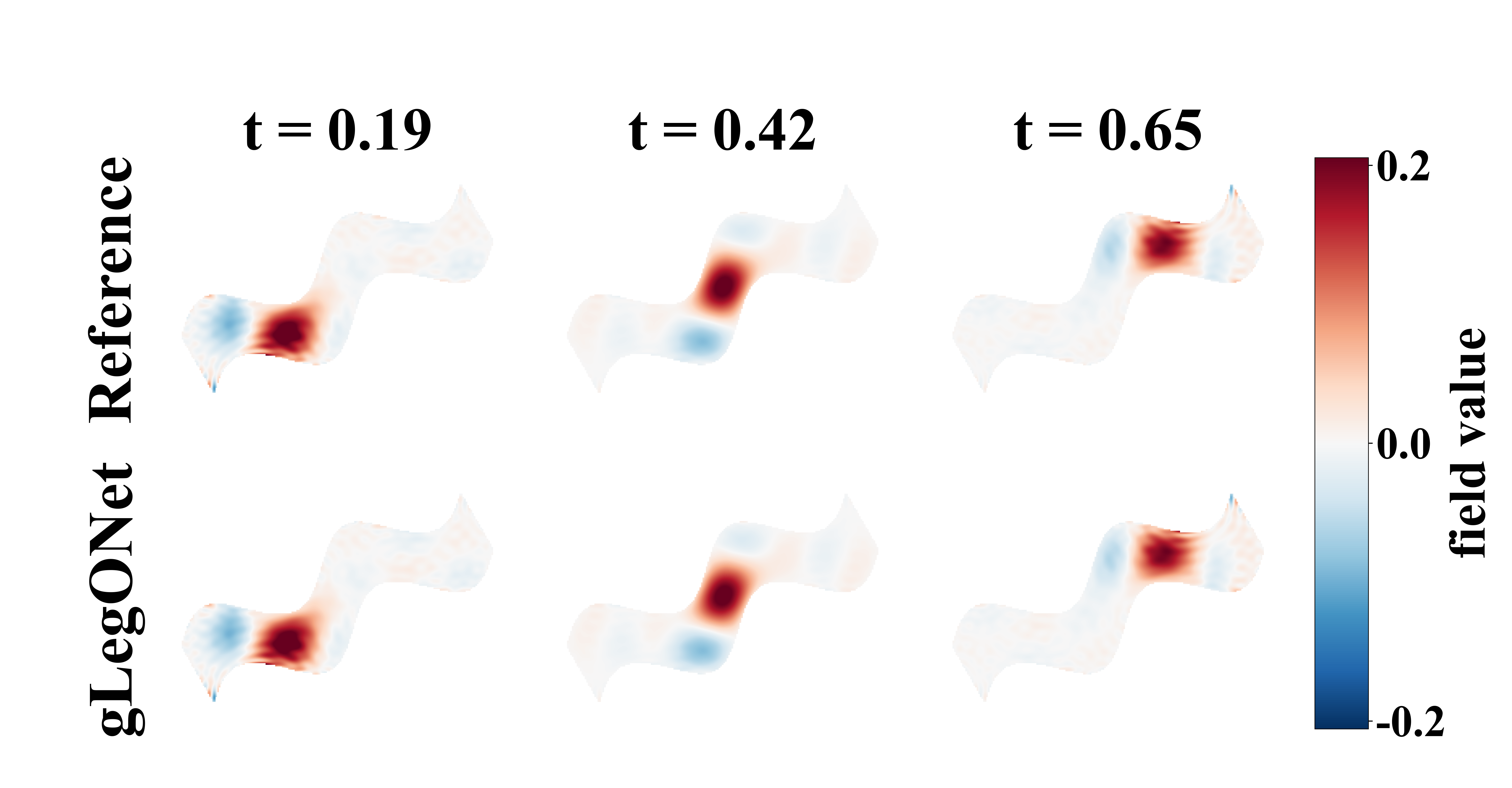}

\vspace{0.5mm}
{\scriptsize
\textbf{(d)} Oblique sinusoidal channel: reference and identified
boundary-adapted rollout.}
\end{minipage}
\begin{minipage}[t]{0.48\textwidth}
\centering
\includegraphics[width=\linewidth]{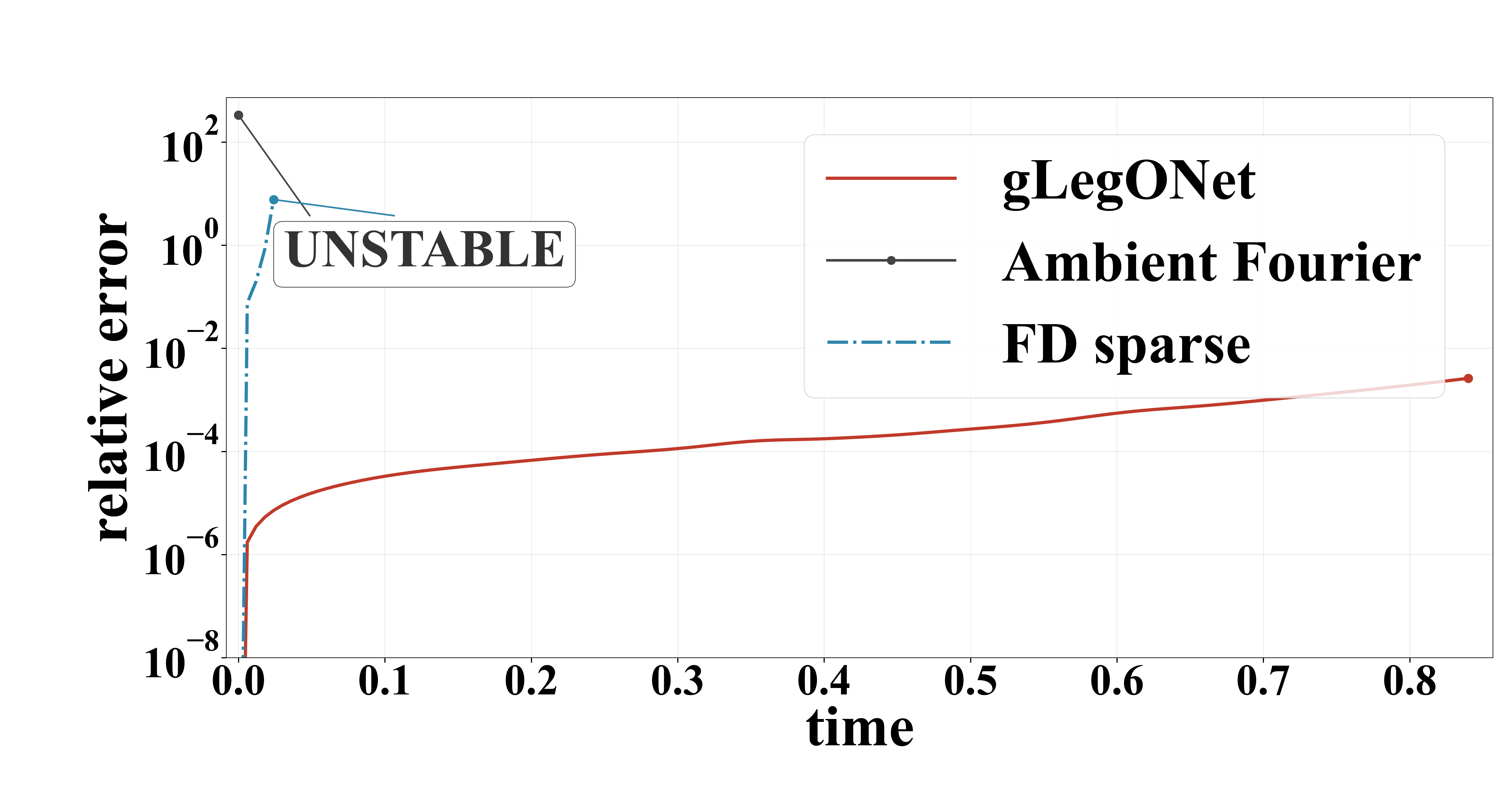}

\vspace{0.5mm}
{\scriptsize
\textbf{(e)} Oblique sinusoidal channel: relative \(L^2(\Omega)\) error
history.}
\end{minipage}
\caption{
Boundary-adapted physical-law identification on unseen domains.
Panel (a) shows active-coefficient error \(e_{\rm coef}\) and final rollout error
\(e_{\rm roll}\) under observation noise, reduced sensor counts, an ambient
Fourier representation and finite-difference baselines.  Entries with
\(\pm\) report the mean and standard deviation over three randomized seeds.
Panels (b)--(e) show representative zero-noise rollout visualizations from fixed
illustrative initial conditions chosen to make the transported structures
clearly visible; these panels are single-seed examples and are separate from
the randomized statistics in Panel (a).  Panels (b) and (c)
show the identified rollout and relative \(L^2(\Omega)\) error history on the
peanut dumbbell, respectively; Panels (d) and (e) show the
corresponding results on the oblique sinusoidal channel.  In Panels
(c) and (e), \emph{Ambient Fourier} uses an ambient Fourier representation of the same
rank without encoding the target boundary; and \emph{FD sparse} denotes the
finite-difference discovery baseline using sparse observations.  The label
\emph{UNSTABLE} indicates that the corresponding baseline diverges before
reaching the final rollout time.  Here \(e_{\rm coef}\) is the relative error
of the active mechanism coefficients and \(e_{\rm roll}\) is the final
relative \(L^2(\Omega)\) error at \(T_{\rm roll}=1.0\), using only data from
\(T_{\rm id}=0.012\).
}
\label{fig:lawid_main}
\end{figure}

Fig.~\ref{fig:lawid_main}(a) summarizes coefficient recovery and long-time
prediction under noise and sensor ablations.  Panels (b,c) and (d,e) show
representative peanut and channel rollouts, respectively.  Together, these results show that gLegONet turns sparse observations into
predictive physical laws by evaluating frozen mechanism blocks in a
geometry-compatible boundary-adapted coordinate system.

\section*{Discussion}\label{sec:Discussion}

This work introduces Geometry-aware LegONet as a geometry-transfer layer for reusable PDE learning. The central idea is to separate two roles that are often entangled in neural PDE solvers: learning interior physical mechanisms and realizing those mechanisms on a particular geometry. The physical mechanisms are represented by Lego-like variational blocks pretrained once in an ambient coefficient space. The geometry and boundary condition are represented separately through algebraic constraints that define boundary-adapted coordinates. In this formulation, changing the domain changes the coordinate interface, not the learned physics.

The results show that this separation is useful across both forward and inverse problems. In the manufactured benchmarks, one frozen ambient block library transfers across several boundary operators without target-domain neural retraining. In the Allen--Cahn disk problem, hard Neumann coordinates help preserve a constrained nonlinear evolution against an independent finite-element reference. In the vector Burgers experiment, nonlinear transport blocks transfer with little loss relative to a geometry-specific block construction when the target domain is resolved by the ambient space. The cylinder-wake and clamped Swift--Hohenberg examples further show that the same boundary interface can support Poisson coupling and higher-order operators assembled from lower-order reusable blocks. Finally, the law-identification tests show that boundary-adapted coordinates provide stable candidate features for discovering governing equations from short-time data on unseen domains.

A useful interpretation is that the method changes the unit of reuse in scientific machine learning. In many neural-operator approaches, the reusable object is a full solution map on a prescribed geometry or mesh distribution. In the original LegONet framework, the reusable object is a physical mechanism block. In the present work, these reusable mechanism blocks are connected to new geometries by an algebraic boundary interface. This suggests a route toward scientific machine-learning libraries in which physical mechanisms are pretrained once and deployed on new domains by boundary assembly, rather than by retraining a monolithic solver for every new configuration.

The framework also clarifies where errors come from. Because the boundary condition is embedded in the coordinate system, boundary leakage is largely removed as a dominant failure mode. When the boundary residual remains near numerical precision but the field error increases, the cause is not weak boundary enforcement. It is instead associated with the expressiveness of the ambient basis, the retained rank of the boundary-adapted space, the accuracy of the pretrained blocks, quadrature or time integration. This separation is difficult to obtain in monolithic neural solvers, where geometry, boundary constraints and learned dynamics are usually mixed inside one model.

The method has several practical considerations. It inherits the resolution limits of the chosen ambient trial space: if the target geometry contains fine boundary features or requires modes beyond the selected cutoff, the boundary can remain satisfied while the interior approximation error grows. The deployment cost also differs from a fully amortized black-box neural time-stepper. A new target geometry requires boundary sampling, null-space construction, reduced operator assembly and coefficient-space rollout. However, these are algebraic and numerical assembly costs rather than neural retraining costs. Once the ambient block library has been trained, new geometries can be handled by rebuilding the boundary-coordinate interface, which is precisely the intended advantage of the framework. The formulation is not restricted to two dimensions, but large-scale three-dimensional arbitrary-boundary applications will require efficient implementations of boundary sampling, quadrature, rank selection and sparse reduced algebra. Strongly nonlocal, stochastic, multiscale or geometry-dependent mechanisms may also require additional reusable blocks or hybrid coupling with classical solvers.

These considerations clarify the role of Geometry-aware LegONet. The method is most useful when the same physical mechanisms appear across many geometries and boundary conditions, and when hard boundary satisfaction, interpretability and mechanism reuse are important. It is not designed merely to be the fastest black-box predictor for one fixed geometry. Rather, it provides a reusable mechanism library together with an algebraic interface for deploying that library on new domains. In this sense, the work extends Lego-like operator learning from physics modularity to geometry modularity. The first step was to make physical mechanisms reusable; the step taken here is to make those mechanisms deployable on new domains through boundary-adapted coordinates. This turns arbitrary-domain PDE learning from a geometry-specific retraining problem into an algebraic assembly problem over reusable scientific operators.

\section*{Methods}\label{sec:Methods}

\subsection*{Ambient representation and reusable mechanism blocks}

Geometry-aware LegONet does not learn a new trajectory model for each geometry
or boundary condition.  It pretrains reusable physical mechanism blocks once on
an ambient domain, using a shared coefficient representation, and
then realizes the frozen blocks on each target domain through boundary-adapted Galerkin coordinates.  Thus geometry and boundary data
change the linear realization layer, while the neural block parameters remain
fixed.

We first place all target problems inside a fixed regular base domain.  In two
dimensions, we use the ambient square \(Q=[-1,1]^2\), the standard
tensor-product box for Fourier representation and quadrature.  A physical
domain \(\Omega_{\rm phys}\) is mapped into \(Q\) by an affine rescaling
\(T_\Omega:\Omega_{\rm phys}\rightarrow\Omega\subset Q\), and the rescaled PDE
is posed on the embedded target domain \(\Omega\).  The mapping consists
only of translation, isotropic scaling and, when needed, a rigid rotation, so it
does not change the physical mechanism.
In Supplementary Information, we specify this
physical-to-embedded rescaling before block pretraining and transfer are
described.

In all experiments, the ambient trial space is the truncated real Fourier
space on \(Q\).  Let
\(\mathcal I_K=\{(k,\ell)\in\mathbb Z^2:k^2+\ell^2\le K^2\}\) and let
\(\mathcal I_K^+=\{(k,\ell)\in\mathcal I_K:k>0\ {\rm or}\ (k=0,\ell>0)\}\),
which selects one representative from each nonzero pair
\((k,\ell)\) and \((-k,-\ell)\).  We use the normalized real Fourier basis
\begin{equation*}
    \begin{aligned}
        \mathcal V_K(Q)
    &=
    \operatorname{span}
    \Bigl(
        \{1\}
        \cup
        \{\sqrt{2}\cos(\pi(kx+\ell y)),
          \sqrt{2}\sin(\pi(kx+\ell y)):
          (k,\ell)\in\mathcal I_K^+\}
    \Bigr)\\
    &:=
    \operatorname{span}\{\phi_1,\ldots,\phi_M\}.
    \end{aligned}
\end{equation*}
The basis is orthonormal with respect to the normalized inner product
\(\langle u,v\rangle_Q=|Q|^{-1}\int_Q u(x)v(x)\,dx\).  
The choice of Fourier functions is deliberate.  They provide a complete
orthonormal coordinate system on the regular ambient box \(Q\), allow analytic
differentiation of basis functions, and give explicit coefficient projection.
Let \(P_K\) be the orthogonal projection onto \(\mathcal V_K(Q)\) under the
normalized inner product.  The key property we use
is density in \(L^2(Q)\) \citep{trefethen2000spectral,boyd2001chebyshev}: for every \(u\in L^2(Q)\),
\[
    \|u-P_Ku\|_{L^2(Q)}\rightarrow 0
    \qquad
    \text{as } K\rightarrow\infty .
\]
Consequently, if a target-domain field \(u\in L^2(\Omega)\) is represented by
an ambient extension \(\tilde u\in L^2(Q)\) with
\(\tilde u|_\Omega=u\), then the restricted Fourier approximation also
converges on the embedded domain $\Omega$.
Thus the Fourier baseplate supplies a universal ambient coefficient interface,
while the actual target geometry and boundary condition are imposed later by
boundary-adapted coordinates.  

\begin{remark}
The same realization principle is not tied to Fourier functions.  One may
replace \(\mathcal V_K(Q)\) by another finite-dimensional trial space. The essential requirement is an approximation-complete trial family: the spaces \(\mathcal V_K(Q)\) should be dense, as \(K\to\infty\), in the function space in which the target PDE is posed. Polynomial spectral bases, finite-element bases and wavelet bases are examples satisfying this requirement in their standard
settings
\citep{canuto2006spectral,brenner2008fem,cohen2003numerical}.  The Fourier
choice in this work is therefore a controlled baseplate for testing reusable
block transfer, rather than an empirical claim that all admissible bases behave
identically.
\end{remark}

We now use this ambient approximation space as the coefficient language for
mechanism learning.
A coefficient vector \(\mathbf a\in\mathbb R^M\) represents the ambient field
\[
    u_{\mathbf a}(\mathbf x)
    =
    \sum_{j=1}^{M}a_j\phi_j(\mathbf x)
    =
    \Phi_Q(\mathbf x)\mathbf a ,
\]
where \(\Phi_Q(\mathbf x)=(\phi_1(\mathbf x),\ldots,\phi_M(\mathbf x))\).  We denote the reconstruction
map by \(\mathcal R_Q:\mathbf a\mapsto u_{\mathbf a}\), and the coefficient
projection by
\(\mathcal P_Q:u\mapsto(\langle u,\phi_1\rangle_Q,\ldots,
\langle u,\phi_M\rangle_Q)\).  Thus every ambient field used by the method is
first represented in the finite coefficient space \(\mathcal V_K(Q)\).

Let \(\Omega\subset Q\) be an embedded target domain.  A target PDE is posed on
\(\Omega\) and decomposed into physical mechanisms as
\[
    \partial_t u(\mathbf x,t)
    =
    \sum_{i=1}^{N_{\rm blk}}
    c_iL_i^\Omega(u)(\mathbf x,t),
    \qquad
    \mathbf x\in\Omega .
\]
Here \(L_i^\Omega\) is the \(i\)th mechanism acting on fields defined on
\(\Omega\), and \(c_i\) is the scalar coefficient.
In our coefficient representation, a target-domain state is
represented by the restriction of an ambient Fourier field
\(u_{\mathbf a}\) to \(\Omega\).  
We define the ambient counterpart \(L_i^Q\) by using the same interior
differential expression on \(Q\), i.e.,
\(\left. L_i^Q(u)\right|_\Omega\) agrees with
\(L_i^\Omega(\left.u\right|_\Omega)\) in the interior of \(\Omega\). Boundary conditions are not part of this ambient mechanism. They enter only through the coordinate interface in the next subsection.
For \(\mathbf a\in\mathbb R^M\), the ambient mechanism \(L_i^Q\) induces the
coefficient-space vector field
\[
    \mathbf F_i(\mathbf a)
    =
    \mathcal P_Q
    \left[
        L_i^Q(\mathcal R_Q\mathbf a)
    \right]
    \in\mathbb R^M .
\]
The pretraining objective is to construct a neural block
\(\mathbf F_i^\theta:\mathbb R^M\rightarrow\mathbb R^M\) that approximates
\(\mathbf F_i\).  This block is therefore a coefficient-space representation
of one mechanism on \(Q\): for instance, a Laplace block maps the coefficients
of \(u_{\mathbf a}\) to the Fourier coefficients of \(\Delta u_{\mathbf a}\).
Note that, at this stage, the block does not encode a target geometry, a boundary
condition, or a complete trajectory model.  It only learns the reusable action
of a single physical mechanism on the shared ambient coefficient state.

Following our previous structure-preserving LegONet design \cite{zhang2026legonet}, a mechanism block can be
written in the generator form
\[
    \mathbf F_i^\theta(\mathbf a)
    =
    -G_i\nabla_{\mathbf a}E_i^\theta(\mathbf a)
    +
    J_i\nabla_{\mathbf a}H_i^\theta(\mathbf a)
    +
    \mathbf R_i(\mathbf a).
\]
Here $G_i$ is a symmetric positive-semidefinite structure operator,
\(J_i^\top=-J_i\),
\(E_i^\theta:\mathbb R^M\rightarrow\mathbb R\) is a dissipative generator,
\(H_i^\theta:\mathbb R^M\rightarrow\mathbb R\) is a conservative generator, and
\(\mathbf R_i:\mathbb R^M\rightarrow\mathbb R^M\) contains residual
effects.  This formula is a template inherited from our previous work \cite{zhang2026legonet}.  In the
present library, blocks are deliberately single-purpose: a diffusion or
relaxation mechanism is represented as an \(E\)-block, a transport or
Hamiltonian mechanism as an \(H\)-block, and a forcing, reaction or closure
mechanism as an \(R\)-block.  This modular design is what makes the library
plug-and-play: mechanisms can be selected, scaled and recombined without
retraining a monolithic solver.

The block parameters are obtained by supervised pretraining \cite{zhang2026legonet}.  For each mechanism, training data are generated
from either an analytic expression or a trusted reference discretization, and
\(\theta\) is optimized by an \(L^2\)-type supervised loss between the predicted
and reference mechanism outputs.  The exact sampling rules, reference
mechanisms and losses used for each benchmark are reported in the Supplementary Information.

After pretraining, the ambient library on \(Q\) is frozen:
\[
    \mathcal B_\theta^Q
    =
    \{\mathbf F_i^\theta\}_{i=1}^{N_{\rm lib}} .
\]
A target PDE, however, is posed only on the embedded domain
\(\Omega\subset Q\).  The next subsection explains how the boundary condition
on \(\partial\Omega\) defines an admissible coordinate interface inside the
ambient coefficient space, so that the frozen library \(\mathcal B_\theta^Q\)
can be realized on \(\Omega\) without retraining its neural parameters.

\subsection*{Boundary-adapted coordinates on target domains}

The frozen library \(\mathcal B_\theta^Q\) acts on Fourier coefficients on the
ambient square \(Q\), whereas a target problem is posed on an embedded domain
\(\Omega\subset Q\) with its own boundary operator.  Geometry-aware LegONet
resolves this mismatch by changing coordinates rather than retraining the
neural blocks.  We construct an affine coordinate space in which every
numerical state satisfies the sampled target boundary condition before any
block composition or time integration is applied.

Let
\(\{(\mathbf x_b,\mathbf n_b)\}_{b=1}^{N_b}\) be samples on
\(\partial\Omega\), where \(\mathbf n_b\) is the outward unit normal.  We
consider Dirichlet, Neumann and Robin trace operators
\(\mathcal B_Du=u\),
\(\mathcal B_Nu=\partial_{\mathbf n}u\), and
\(\mathcal B_Ru=\partial_{\mathbf n}u+\kappa u\), with prescribed data
\(\mathcal B_\star u=g_\star(\mathbf x,t)\).  Here
\(\partial_{\mathbf n}u=\mathbf n\cdot\nabla u\), and \(\kappa\) is the Robin
coefficient.

After stacking all scalar trace constraints, boundary sampling gives
\begin{equation}\label{eq:boundary-system}
    C\mathbf a=\mathbf d(t),
\end{equation}
where \(C\in\mathbb R^{N_b\times M}\),
\(\mathbf d(t)\in\mathbb R^{N_b}\).  For a row \(\ell\) generated from the sample
\(\mathbf x_b\), the entries are
\[
    C_{\ell,:}
    =
    \begin{cases}
        \Phi_Q(\mathbf x_b),
        & \star=D,\\
        \partial_{\mathbf n}\Phi_Q(\mathbf x_b),
        & \star=N,\\
        \partial_{\mathbf n}\Phi_Q(\mathbf x_b)
        +\kappa(\mathbf x_b)\Phi_Q(\mathbf x_b),
        & \star=R,
    \end{cases}
    \qquad
    d_\ell(t)=g_\star(\mathbf x_b,t),
\]
where
\(\partial_{\mathbf n}\Phi_Q(\mathbf x_b)
=\mathbf n_b^\top\nabla\Phi_Q(\mathbf x_b)\in\mathbb R^{1\times M}\).
Mixed conditions are
handled by assigning different trace types to different boundary segments;
clamped conditions are obtained by stacking both value and normal-derivative
rows.

The sampled boundary-admissible coefficient set is
\[
    \mathcal M_{\partial\Omega}(t)
    =
    \left\{
        \mathbf a\in\mathbb R^M:
        C\mathbf a=\mathbf d(t)
    \right\}.
\]
Its tangent space is \(\ker C\).  We compute a numerical basis for this
nullspace from the full singular value decomposition (SVD)
\begin{equation}\label{eq:boundary-svd}
    C=U\Sigma V^\top .
\end{equation}
Let \(p=\min(N_b,M)\),
\(V=[\mathbf v_1,\ldots,\mathbf v_M]\in\mathbb R^{M\times M}\), and order the
singular values as \(\sigma_1\ge\cdots\ge\sigma_p\ge0\).  Using the numerical
rank threshold
\(r_C=\#\{j\le p:\sigma_j>\tau_C\max(\sigma_1,1)\}\), we define
\[
    Z
    =
    [\mathbf v_{r_C+1},\ldots,\mathbf v_M]
    \in\mathbb R^{M\times m},
    \qquad
    m=M-r_C .
\]
Then \(CZ\approx\mathbf0\), so the columns of \(Z\) span the coefficient
directions left unconstrained by the sampled boundary operator
\citep{golub2013matrix}.  Unless otherwise stated, we use
\(\tau_C=10^{-10}\).

The Galerkin metric is defined on the physical domain \(\Omega\).  Let
\(\{(\mathbf x_q,w_q)\}_{q=1}^{N_{\rm q}}\) be volume quadrature points and
weights, define
\((\Phi_\Omega)_{qj}=\phi_j(\mathbf x_q)\), and set
\(W_\Omega=\operatorname{diag}(w_1,\ldots,w_{N_{\rm q}})\) and
\(M_\Omega=\Phi_\Omega^\top W_\Omega\Phi_\Omega\).
Here \(N_{\rm q}\) denotes the number of volume quadrature points, whereas
\(N_\Omega\) below denotes the boundary-adapted coordinate matrix.
The raw nullspace basis \(Z\) may contain directions with negligible
\(L^2(\Omega)\) mass or poor scaling.  We therefore apply the standard
Gram-matrix normalization
\citep{golub2013matrix,trefethen2022numerical}.  Form
\(G=Z^\top M_\Omega Z\in\mathbb R^{m\times m}\) and compute
\(G=P\Lambda P^\top\), where
\(P=[\mathbf p_1,\ldots,\mathbf p_m]\) is orthogonal and
\(\Lambda=\operatorname{diag}(\lambda_1,\ldots,\lambda_m)\), with
\(\lambda_1\ge\cdots\ge\lambda_m\ge0\).  We retain the eigenvalues satisfying
\(\lambda_j>\tau_M\max(\lambda_1,1)\).  If \(r\) directions remain, let
\(P_r=[\mathbf p_1,\ldots,\mathbf p_r]\) and
\(\Lambda_r=\operatorname{diag}(\lambda_1,\ldots,\lambda_r)\).  The
boundary-adapted coordinate matrix is
\begin{equation}\label{eq:boundary-adapted-N}
    N_\Omega
    =
    ZP_r\Lambda_r^{-1/2}
    \in\mathbb R^{M\times r}.
\end{equation}
Here \(r\) is the retained reduced-space dimension reported in the numerical
experiments.  We use \(\tau_M=10^{-11}\) unless otherwise stated.  By
construction,
\[
    CN_\Omega\approx\mathbf0,
    \qquad
    N_\Omega^\top M_\Omega N_\Omega=I_r .
\]
Thus every homogeneous state
\(\mathbf a=N_\Omega\mathbf z\), with
\(\mathbf z\in\mathbb R^r\), satisfies the sampled boundary condition to the
numerical rank tolerance and is orthonormal in the physical-domain mass inner
product.

For nonhomogeneous data, we compute an affine lift
\(\mathbf a_{\rm bc}(t)\in\mathbb R^M\) satisfying
\(C\mathbf a_{\rm bc}(t)\approx\mathbf d(t)\).  Every admissible numerical
state is then represented as
\[
    \mathbf a(t)
    =
    \mathbf a_{\rm bc}(t)+N_\Omega\mathbf z(t),
    \qquad
    \mathbf z(t)\in\mathbb R^r .
\]
For homogeneous data, \(\mathbf a_{\rm bc}(t)=\mathbf0\).  The target boundary
condition is therefore built into the state coordinates: frozen ambient blocks
are evaluated only at boundary-admissible coefficient states, and their
responses are restricted to the same tangent space in the next subsection.

This construction uses only boundary samples, volume quadrature and numerical linear algebra. No neural parameter is updated for a new domain or boundary condition.

\subsection*{Galerkin realization and rollout of frozen blocks}

The preceding subsections provide two ingredients: the frozen ambient library
\(\mathcal B_\theta^Q=\{\mathbf F_i^\theta\}_{i=1}^{N_{\rm lib}}\) and the
boundary-adapted coordinate map
\(\mathbf a(t)=\mathbf a_{\rm bc}(t)+N_\Omega\mathbf z(t)\).  We now realize
the reusable mechanisms directly on the target domain.

Recall that the target equation is decomposed as
\[
    \partial_t u(\mathbf x,t)
    =
    \sum_{i=1}^{N_{\rm blk}}
    c_i L_i^\Omega(u)(\mathbf x,t),
    \qquad
    \mathbf x\in\Omega ,
\]
where \(c_i\) is the prescribed scalar weight of the \(i\)th mechanism.

At the admissible state
\(\mathbf a_\star(\mathbf z,t)
=\mathbf a_{\rm bc}(t)+N_\Omega\mathbf z\), a default block returns the
complete ambient coefficient response
\(\mathbf F_i^\theta(\mathbf a_\star)\in\mathbb R^M\).  Let
\(\mathbf q_i^\Omega(\mathbf z,t)\in\mathbb R^r\) denote the reduced
coordinates of its \(M_\Omega\)-orthogonal projection onto
\(\operatorname{range}(N_\Omega)\).  Equivalently,
\(N_\Omega\mathbf q_i^\Omega(\mathbf z,t)\) is defined by the Galerkin
orthogonality condition
\[
    N_\Omega^\top M_\Omega
    \left[
        N_\Omega\mathbf q_i^\Omega(\mathbf z,t)
        -
        \mathbf F_i^\theta
        \left(
            \mathbf a_\star(\mathbf z,t)
        \right)
    \right]
    =
    \mathbf 0 .
\]
Using \(N_\Omega^\top M_\Omega N_\Omega=I_r\) gives
\begin{equation}\label{eq:default-reduced-response}
    \mathbf q_i^\Omega(\mathbf z,t)
    =
    N_\Omega^\top M_\Omega
    \mathbf F_i^\theta
    \left(
        \mathbf a_{\rm bc}(t)+N_\Omega\mathbf z
    \right)
    \in\mathbb R^r .
\end{equation}
Thus Eq.~\eqref{eq:default-reduced-response} is the unique Galerkin coordinate vector of the ambient block response in the
boundary-adapted tangent space.

For each mechanism, let \(\mathbf q_i^\Omega(\mathbf z,t)\) denote its realized
target-domain response.  Since Galerkin restriction is linear, these responses
add with the same coefficients \(c_i\) as the target PDE.  Moreover,
\(\dot{\mathbf a}(t)=\dot{\mathbf a}_{\rm bc}(t)
+N_\Omega\dot{\mathbf z}(t)\).  Projecting this coefficient velocity onto the
boundary-adapted space gives
\[
    N_\Omega^\top M_\Omega
    \left[
        \dot{\mathbf a}_{\rm bc}(t)
        +
        N_\Omega\dot{\mathbf z}(t)
    \right]
    =
    \sum_{i=1}^{N_{\rm blk}}
    c_i\,\mathbf q_i^\Omega(\mathbf z(t),t),
\]
which yields the reduced
dynamics
\begin{equation}\label{eq:reduced-z-ode}
    \dot{\mathbf z}(t)
    =
    \sum_{i=1}^{N_{\rm blk}}
    c_i\,\mathbf q_i^\Omega(\mathbf z(t),t)
    -
    N_\Omega^\top M_\Omega
    \dot{\mathbf a}_{\rm bc}(t).
\end{equation}
Equation~\eqref{eq:reduced-z-ode} is the finite-dimensional system advanced by
Geometry-aware LegONet.  The boundary condition is enforced by the coordinates
throughout the rollout, rather than by a penalty or post-step correction.

Time integration is problem dependent.  The default LegONet strategy composes
the realized mechanisms as subflows, with Strang splitting as the standard
second-order symmetric composition
\citep{strang1968splitting}.  Explicit Runge--Kutta, implicit or
Crank--Nicolson treatment of stiff linear terms, and IMEX updates may be used
when required by the target dynamics.

Note that the same interface supports operator composition.  Higher-order, coupled and
auxiliary mechanisms are assembled from reusable responses
\(\mathbf q_i^\Omega\), rather than learned as new monolithic networks.  Frozen
blocks can therefore be selected, scaled, reassembled and composed on new
domains without geometry-specific neural retraining.

\begin{algorithm}[H]
\caption{Geometry-aware realization and rollout of frozen LegONet blocks}
\label{alg:ga-legonet}
\begin{algorithmic}[1]
\Require Ambient domain \(Q\), cutoff \(K\), mechanism expressions
\(\{L_i^Q\}_{i=1}^{N_{\rm lib}}\), target domain \(\Omega\), and initial condition on \(\Omega\).
\Ensure Boundary-admissible rollout \(u_h(\mathbf x,t)\) on \(\Omega\).

\Statex \textbf{Offline: train reusable mechanisms on \(Q\)}
\State Construct the ambient basis
\(\mathcal V_K(Q)=\operatorname{span}\{\phi_1,\ldots,\phi_M\}\).
\For{\(i=1,\ldots,N_{\rm lib}\)}
    \State Supervise the mechanism \(L_i^Q\) and train
    \(\mathbf F_i^\theta:\mathbb R^M\rightarrow\mathbb R^M\).
\EndFor
\State Freeze the ambient library
\(\mathcal B_\theta^Q=\{\mathbf F_i^\theta\}_{i=1}^{N_{\rm lib}}\).

\Statex \textbf{Online: realize the target boundary}
\State Assemble the boundary system in Eq.~\eqref{eq:boundary-system}.
\State Compute \(N_\Omega\) using
Eqs.~\eqref{eq:boundary-svd} and \eqref{eq:boundary-adapted-N}.
\State Compute \(\mathbf a_{\rm bc}(t)\) and represent
\(\mathbf a(t)=\mathbf a_{\rm bc}(t)+N_\Omega\mathbf z(t)\).

\Statex \textbf{Rollout: evolve only admissible coordinates}
\State Project or fit the initial condition to obtain \(\mathbf z(0)\).
\For{each time step}
    \State Evaluate the required frozen blocks
    \(\mathbf F_i^\theta(\mathbf a_{\rm bc}(t)+N_\Omega\mathbf z(t))\).
    \State Advance the reduced ODE in Eq.~\eqref{eq:reduced-z-ode}.
\EndFor
\State Reconstruct
\[
    u_h(\mathbf x,t)
    =
    \Phi_Q(\mathbf x)
    \left(
        \mathbf a_{\rm bc}(t)+N_\Omega\mathbf z(t)
    \right),
    \qquad
    \mathbf x\in\Omega .
\]
\end{algorithmic}
\end{algorithm}

Algorithm~\ref{alg:ga-legonet} separates ambient pretraining, target-domain
realization and boundary-admissible rollout.  A standard stability argument
decomposes the resulting error into boundary-adapted approximation,
mechanism-response mismatch and time-discretization terms; see Supplementary
Information.

\subsection*{Sparse law identification in boundary-adapted coordinates}

We now consider the inverse problem of identifying an unknown interior law on a
known target domain \(\Omega\) with a prescribed boundary condition.  Let
\(\{L_p^\Omega\}_{p=1}^{P}\) be a candidate set of physical mechanisms.  We
assume that the unknown dynamics can be approximated by
\[
    \partial_t u(\mathbf x,t)
    \approx
    \sum_{p=1}^{P}
    c_p L_p^\Omega(u)(\mathbf x,t),
    \qquad
    \mathbf x\in\Omega ,
\]
with a subset of the coefficients
\(\mathbf c=(c_1,\ldots,c_P)^\top\in\mathbb R^P\) active.  Sparse observations
are available as
\(y_{n,s}=u(\mathbf x_s,t_n)+\eta_{n,s}\), where
\(\{\mathbf x_s\}_{s=1}^{N_{\rm obs}}\subset\Omega\) are interior sensors,
\(\{t_n\}_{n=1}^{N_t}\subset[0,T_{\rm id}]\) is a short identification window,
and \(\eta_{n,s}\) denotes observation noise.  The goal is to recover the active
mechanisms and their coefficients, then use the identified law for prediction
beyond \(T_{\rm id}\).  

Let \(\mathcal D_\theta^Q\) denote the frozen ambient dictionary corresponding
to the \(P\) candidate mechanisms.  Each candidate is realized on
\(\Omega\) using the default coefficient response in
Eq.~\eqref{eq:default-reduced-response} or the mechanism-dependent procedure in
Remark~\ref{rem:mechanism-dependent-realization}.  This gives the
target-domain response dictionary
\[
    \mathcal D_\theta^\Omega
    =
    \left\{
        \mathbf q_p^\Omega:
        \mathbb R^r\times[0,T_{\rm id}]
        \rightarrow\mathbb R^r
    \right\}_{p=1}^{P}.
\]
The unknown law in boundary-adapted coordinates is therefore
\[
    \dot{\mathbf z}(t)
    +
    N_\Omega^\top M_\Omega
    \dot{\mathbf a}_{\rm bc}(t)
    \approx
    \sum_{p=1}^{P}
    c_p\,\mathbf q_p^\Omega(\mathbf z(t),t).
\]

To infer the reduced trajectory, define the sensor matrix by
\((\Phi_{\rm obs})_{sj}=\phi_j(\mathbf x_s)\), and let
\(\mathbf y_n=(y_{n,1},\ldots,y_{n,N_{\rm obs}})^\top\).  Each observation is
mapped into the boundary-adapted coordinates through
\begin{equation}\label{eq:sensor-to-z}
    \mathbf z_n
    =
    \arg\min_{\mathbf z\in\mathbb R^r}
    \left\|
        \Phi_{\rm obs}
        \left(
            \mathbf a_{\rm bc}(t_n)+N_\Omega\mathbf z
        \right)
        -
        \mathbf y_n
    \right\|_2^2 .
\end{equation}
Temporal derivatives \(\dot{\mathbf z}_n\) are estimated from the short reduced
trajectory. 
The mechanism coefficients are then identified by a regularized linear
least-squares problem,
\begin{equation}\label{eq:law-id-regression}
    \widehat{\mathbf c}
    =
    \arg\min_{\mathbf c\in\mathbb R^P}
    \sum_{n=1}^{N_t}
    \left\|
        \dot{\mathbf z}_n
        +
        N_\Omega^\top M_\Omega
        \dot{\mathbf a}_{\rm bc}(t_n)
        -
        \sum_{p=1}^{P}
        c_p\,\mathbf q_p^\Omega(\mathbf z_n,t_n)
    \right\|_2^2
    +
    \lambda_{\rm ridge}\mathcal R(\mathbf c).
\end{equation}
Here \(\lambda_{\rm ridge}\ge0\) controls the regularization penalty
\(\mathcal R(\mathbf c)\).  If the active mechanisms are prescribed, we solve
the least-squares problem on that subset.  If the active subset is unknown, we
use sequential thresholded ridge least squares, with column normalization and
threshold--refit iterations  \citep{brunton2016discovering,rudy2017data,tibshirani1996regression} specified in Supplementary Information.

The identified coefficients immediately define a predictive reduced model:
\begin{equation}\label{eq:discovered-rollout}
    \dot{\mathbf z}(t)
    =
    \sum_{p=1}^{P}
    \hat c_p\,
    \mathbf q_p^\Omega(\mathbf z(t),t)
    -
    N_\Omega^\top M_\Omega
    \dot{\mathbf a}_{\rm bc}(t).
\end{equation}
No separate geometry-specific discretization is required after regression.
Identification and rollout use the same boundary-adapted coordinates and the
same realized mechanism responses.

This formulation targets sparse, randomly distributed and noisy observations
over a short time interval.  Sensors are used only to infer
\(\mathbf z_n\); the spatial mechanism features are supplied by the frozen
dictionary.  By contrast, the pointwise finite-difference baseline requires
dense geometry-matched observations to construct stable derivative features.
The numerical experiments test this distinction directly.

\begin{remark}[Dictionary content]\label{rem:mechanism-dependent-realization}
The dictionary contains candidate physical mechanisms. Its entries may be frozen learned blocks or fixed analytic Galerkin maps.  All candidates enter the regression through the same
reduced response \(\mathbf q_p^\Omega\). The same ambient
dictionary can therefore be reused across target domains while only its
deterministic realization changes.
\end{remark}

\begin{algorithm}[H]
\caption{Sparse law identification in boundary-adapted coordinates}
\label{alg:boundary-law-id}
\begin{algorithmic}[1]
\Require Frozen ambient dictionary \(\mathcal D_\theta^Q\), precomputed target
interface \(N_\Omega,M_\Omega\), affine lift \(\mathbf a_{\rm bc}(t)\), observations \(\{\mathbf y_n\}_{n=1}^{N_t}\).
\Ensure Identified coefficients \(\hat{\mathbf c}\) and target-domain rollout.

\State Build the sensor matrix \(\Phi_{\rm obs}\).
\For{\(n=1,\ldots,N_t\)}
    \State Infer \(\mathbf z_n\) using Eq.~\eqref{eq:sensor-to-z}.
\EndFor
\State Estimate \(\dot{\mathbf z}_n\) from the short reduced trajectory.
\For{\(n=1,\ldots,N_t\) and \(p=1,\ldots,P\)}
    \State Assemble
    \(\mathbf q_p^\Omega(\mathbf z_n,t_n)\) through the target-domain response
    interface.
\EndFor
\State Solve Eq.~\eqref{eq:law-id-regression} for \(\hat{\mathbf c}\).
\State Roll out the identified law using Eq.~\eqref{eq:discovered-rollout}.
\State Reconstruct
\(u_h(\mathbf x,t)=
\Phi_Q(\mathbf x)
(\mathbf a_{\rm bc}(t)+N_\Omega\mathbf z(t))\) on \(\Omega\).
\end{algorithmic}
\end{algorithm}

Algorithm~\ref{alg:boundary-law-id} converts sparse short-time observations
into a predictive law using the same boundary-adapted mechanism responses as the
forward solver.

\subsection*{Relation to classical Galerkin and geometry-aware neural operators}

Having specified the forward and inverse formulations, we now clarify how the proposed framework relates to classical null-space Galerkin methods and modern geometry-aware neural operators. Geometry-aware LegONet combines a classical boundary-compatible coordinate construction with a reusable library of learned physical mechanisms. The null-space idea itself is not new. When the governing PDE is fully known, one may construct a boundary-compatible basis, assemble the required operators in that basis and solve the resulting Galerkin system. For a single prescribed problem, such a classical solver may be the most direct choice.

The difference lies in how the interior mechanisms are represented and reused. In a generic boundary-compatible basis, an energy-based Galerkin discretization commonly takes the form
\[
\dot{\mathbf z}=-\nabla_{\mathbf z} E(\mathbf z),
\]
where \(\nabla_{\mathbf z} E(\mathbf z)\) denotes the energy gradient in coefficient space. Its evaluation still requires deriving and assembling the variational response associated with the chosen basis. For a large, dense or poorly conditioned basis, this step can be costly or numerically sensitive. Classical methods can reduce this cost through orthogonalization or matrix factorization, but the operator construction remains tied to the specified equation and representation.

In gLegONet, the physical mechanisms are pretrained in a shared ambient coefficient space using generator-based blocks,
\[
\mathbf F_i^\theta(\mathbf a) =
-G_i\nabla_{\mathbf a} E_i^\theta(\mathbf a)
+
J_i\nabla_{\mathbf a} H_i^\theta(\mathbf a)
+
\mathbf R_i(\mathbf a).
\]
Their dissipative or conservative structure is encoded through the positive-semidefinite and skew-symmetric operators, while automatic differentiation evaluates the learned generator gradients. For each target geometry, the null-space basis is mass-orthonormalized so that
\[
N_\Omega^\top M_\Omega N_\Omega=I.
\]
After this one-time geometry-dependent construction, the frozen mechanisms produce their reduced responses directly. No repeated reduced mass-matrix solve or new geometry-specific derivation of the learned variational response is required. The main advantage is therefore not the null-space algebra alone, but its use as an interface through which learned mechanisms can be reused, composed and evaluated on new domains.

Modern geometry-aware neural operators address a complementary problem. They usually learn an end-to-end solution map conditioned on geometry or boundary information, using training data sampled from a prescribed distribution of domains. The present experiments instead test whether a fixed mechanism library can be realized on previously unseen domains without geometry-specific neural training. Introducing a geometry-aware neural-operator baseline would require defining and training on a geometry distribution, and would therefore test geometry-conditioned amortization rather than the zero-shot realization of a fixed mechanism library studied here. We do not claim that gLegONet is universally superior to these approaches. Geometry-aware neural operators amortize solution prediction over a learned family of geometries, whereas gLegONet amortizes the learning of physical mechanisms and handles each new geometry through deterministic boundary-adapted assembly.

The framework therefore draws complementary strengths from neural and classical computation. The Galerkin coordinate construction provides algebraic boundary enforcement and efficient reduced evolution. The learned blocks provide data-driven, structure-aware mechanism representations that can be reused when analytic derivation or repeated geometry-specific assembly is difficult. The same realized mechanisms can then support both forward simulation and physical-law identification.

\subsection*{Use of AI-assisted tools}

A large language model, ChatGPT, was used during manuscript preparation to assist with language editing, consistency checks and the identification of passages that could benefit from further clarification. All suggested revisions were critically reviewed by the authors. The authors take full responsibility for the scientific claims, mathematical derivations, numerical results and final content of the manuscript. The AI tool was not used to generate research data or perform the reported numerical experiments.

\section*{Acknowledgements}
We would like to thank the support of National Science Foundation (DMS-2533878, DMS-2053746, DMS-2134209, ECCS-2328241, CBET-2347401 and OAC-2311848), and U.S.~Department of Energy (DOE) Office of Science Advanced Scientific Computing Research program, under the "Uncertainty Quantification for Multifidelity Operator Learning (MOLUcQ)" project (Project No. 81739), DE-SC0023161, the SciDAC LEADS Institute, and DOE–Fusion Energy Science, under grant number: DE-SC0024583.

\section*{Code availability}
The code used in this study, including Python scripts for data analysis, is available at \href{https://github.com/Yooki-YueqiWang/gLegONet}{https://github.com/Yooki-YueqiWang/gLegONet}.

\section*{Data availability} Data used to generate the numerical results in this study are available in the same GitHub repository as the code: \url{https://github.com/Yooki-YueqiWang/gLegONet}.

\section*{Author contributions}

J.Z., Y.W. and G.L. conceived the study and developed the overall methodology. J.Z. contributed to the mathematical formulation, theoretical analysis, algorithmic development and improvement, development of the computational framework, and analysis and interpretation of the numerical results. Y.W. contributed to the mathematical formulation, theoretical analysis and methodological development, developed the numerical implementation, designed and conducted the computational experiments, analysed the results and prepared the figures. G.L. supervised the research, contributed to the methodological development and interpretation of the results, and acquired funding. J.Z. and Y.W. wrote the initial manuscript. All authors discussed the results, reviewed and revised the manuscript, and approved the final version. J.Z., Y.W. and G.L. contributed equally to this work.

\section*{Competing interests}
The authors declare no competing interests.

\section*{Additional information}
\textbf{Correspondence and requests for materials} should be addressed to Guang Lin
(guanglin@purdue.edu).

\bibliography{sn-bibliography}

\clearpage
\section*{Extended Data}

\setcounter{equation}{0}
\setcounter{figure}{0}
\setcounter{table}{0}

\renewcommand{\theequation}{\arabic{equation}}
\renewcommand{\thefigure}{\arabic{figure}}
\renewcommand{\thetable}{\arabic{table}}

\renewcommand{\figurename}{Extended Data Fig.}
\renewcommand{\tablename}{Extended Data Table}

\captionsetup{
  labelfont={bf,color=black},
  textfont=normalfont,
  labelsep=space
}

\begin{figure}[H]
\centering
\begin{subfigure}{0.48\textwidth}
    \centering
    \includegraphics[width=\linewidth]{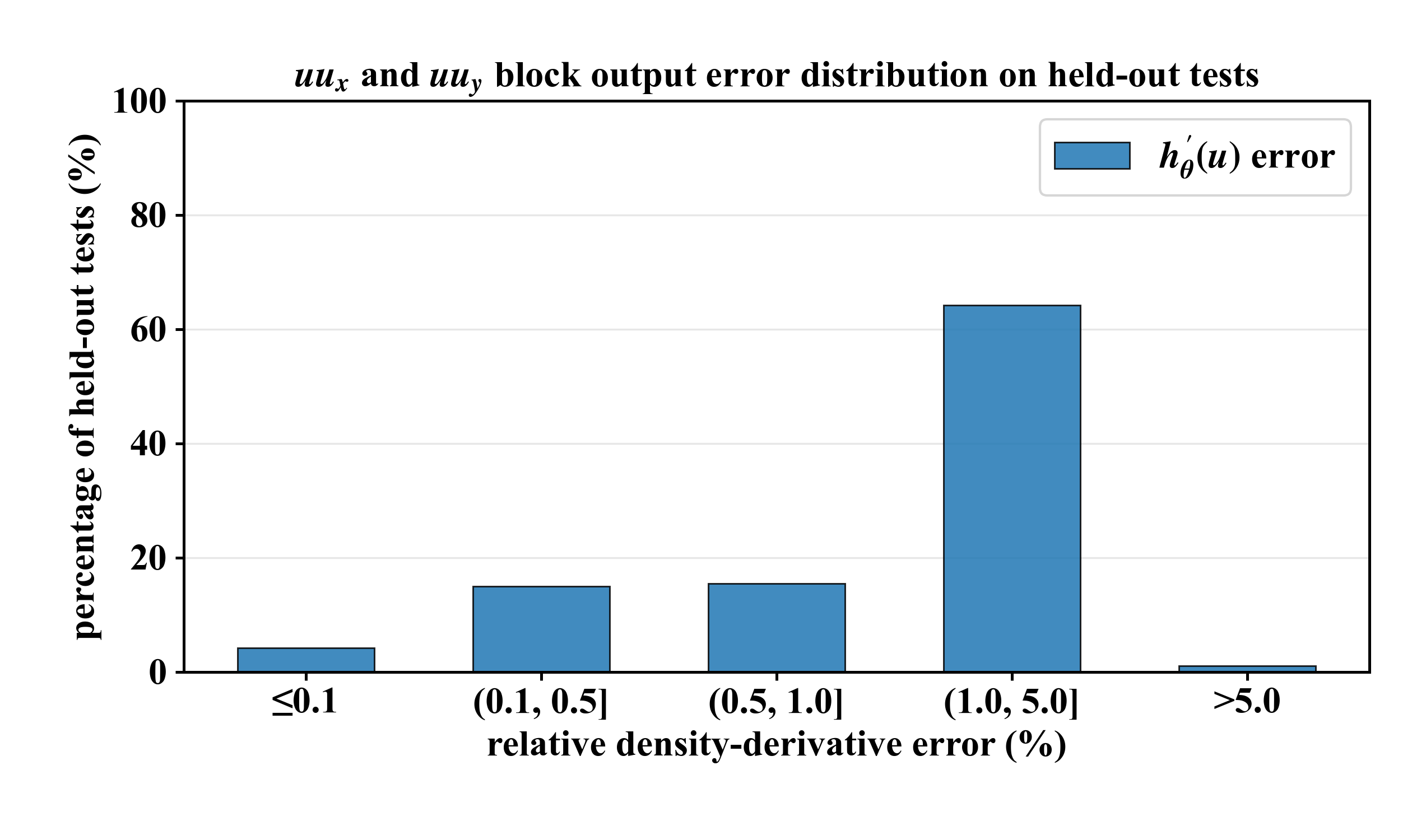}
    \caption{Transport primitive derivative.}
\end{subfigure}
\hfill
\begin{subfigure}{0.48\textwidth}
    \centering
    \includegraphics[width=\linewidth]{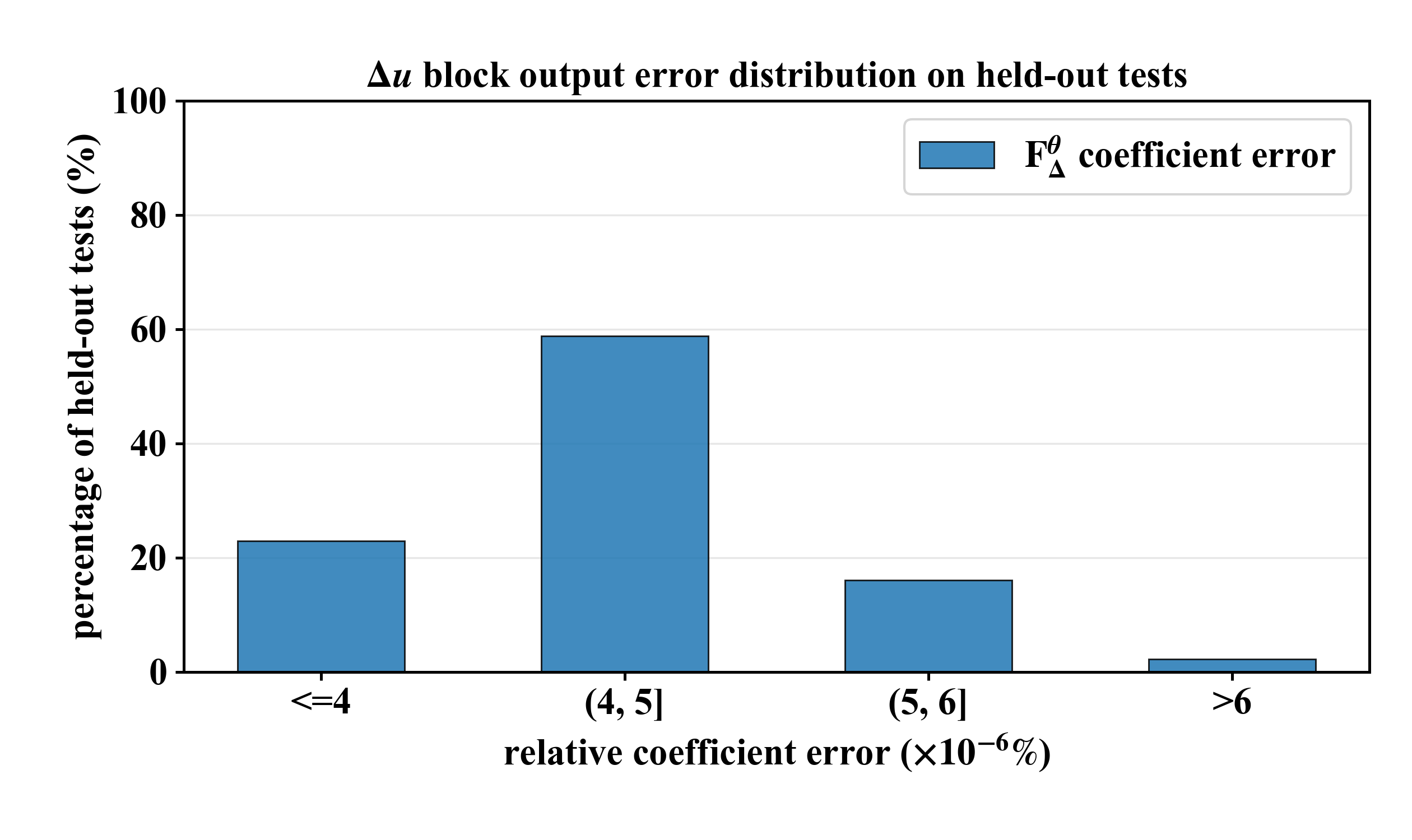}
    \caption{Laplacian coefficient output.}
\end{subfigure}
\caption{ Held-out diagnostics for the reusable ambient blocks. The left panel reports the relative error of the learned transport primitive derivative \(h_\theta'(u)\), and the right panel reports the coefficient-space output error of the Laplace block. Bars show the percentage of held-out samples in each error range. }
\label{fig:ambient_block_error_bins}
\end{figure}

\begin{figure}[H]
\centering
\begin{subfigure}{0.48\textwidth}
    \centering
    \includegraphics[width=\linewidth]{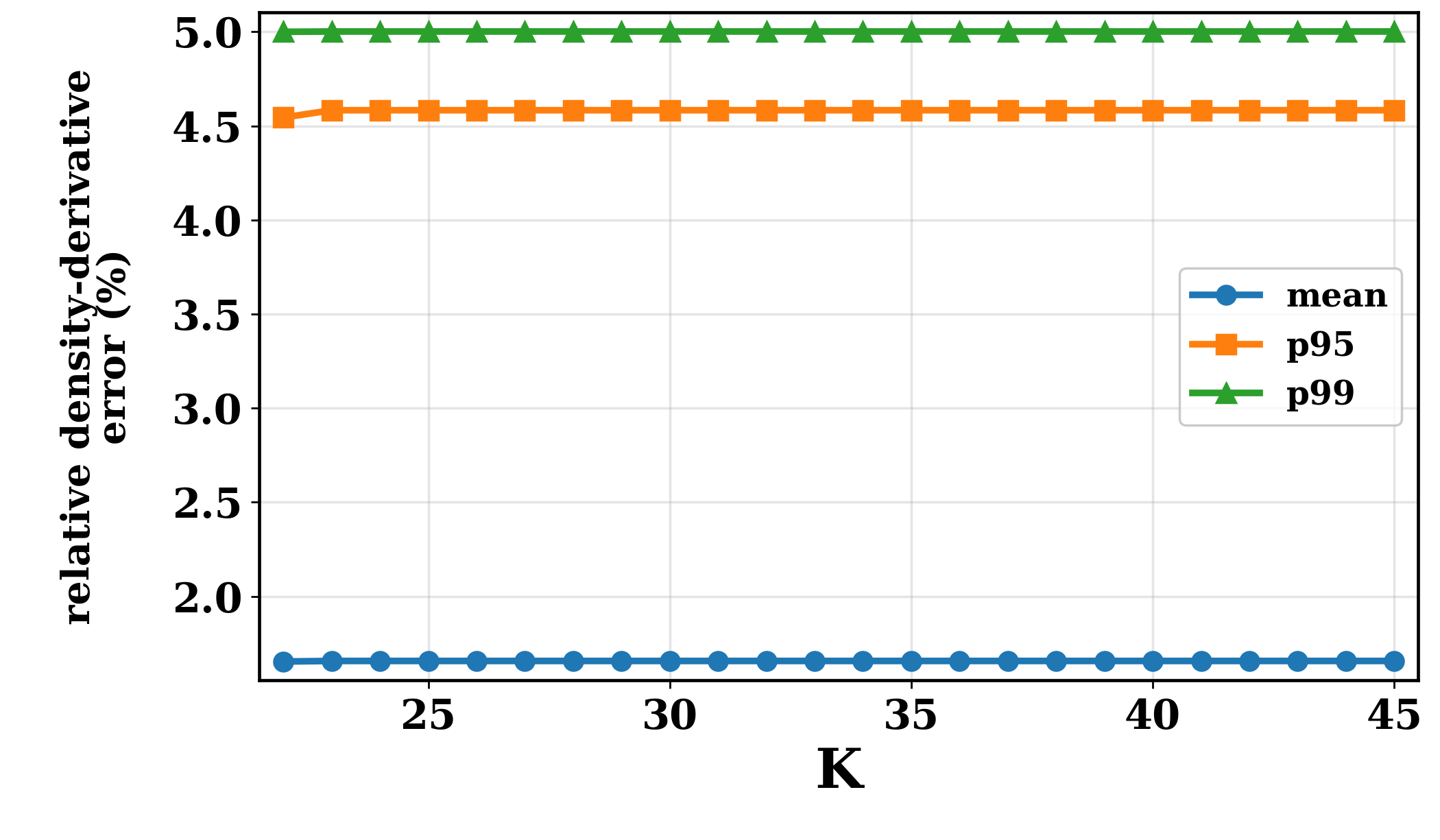}
    \caption{Transport primitive derivative.}
\end{subfigure}
\hfill
\begin{subfigure}{0.48\textwidth}
    \centering
    \includegraphics[width=\linewidth]{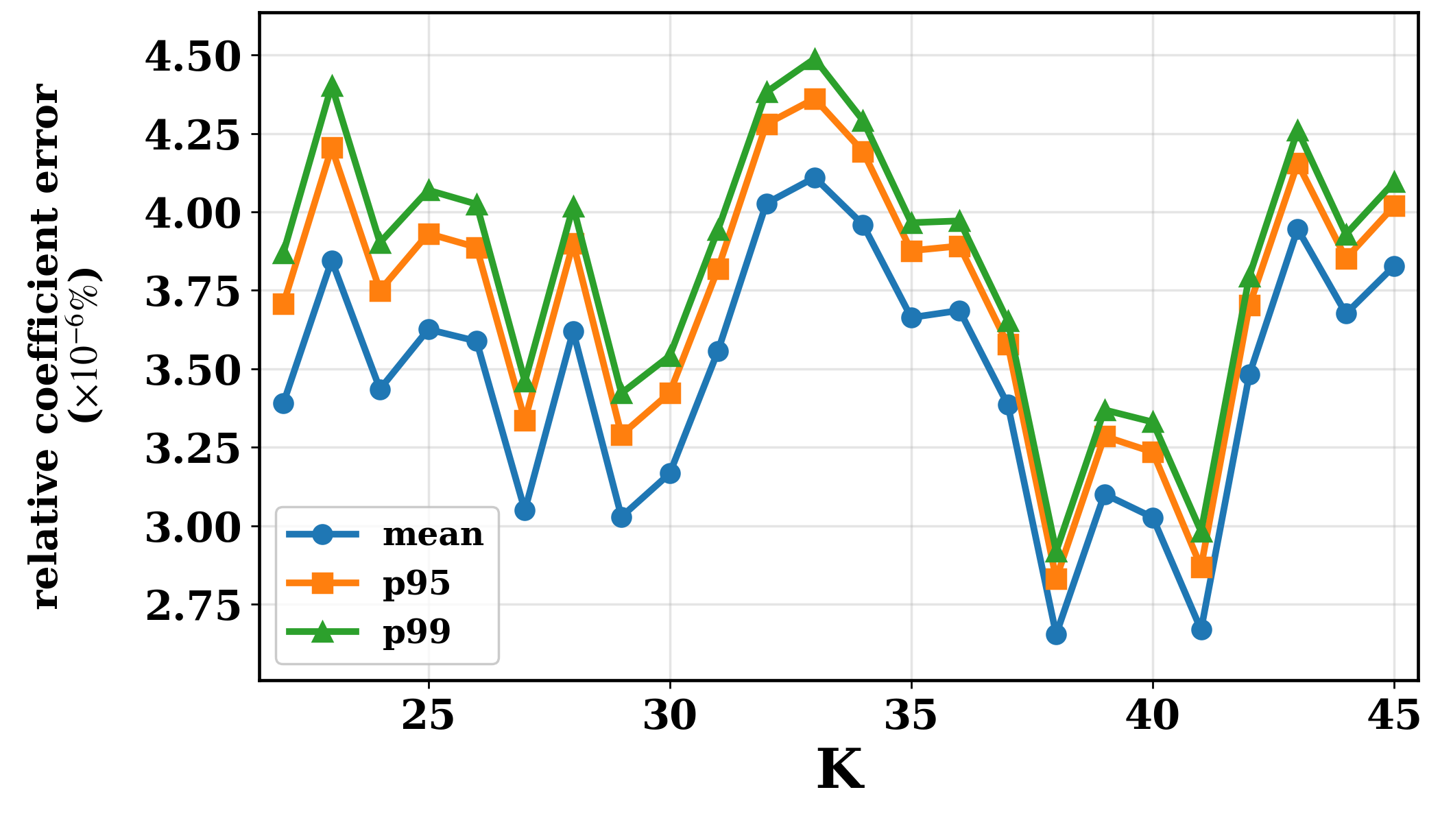}
    \caption{Laplacian coefficient output.}
\end{subfigure}
\caption{ Resolution sensitivity of ambient block pretraining. For each Fourier cutoff \(K\), a separate ambient library is trained and tested in the same space \(\mathcal V_K(Q)\). The diagnostics remain stable as the ambient resolution is increased, indicating that the learned blocks do not degrade over the tested range of \(K\). }
\label{fig:ambient_block_k_sensitivity}
\end{figure}

\begin{figure}[H]
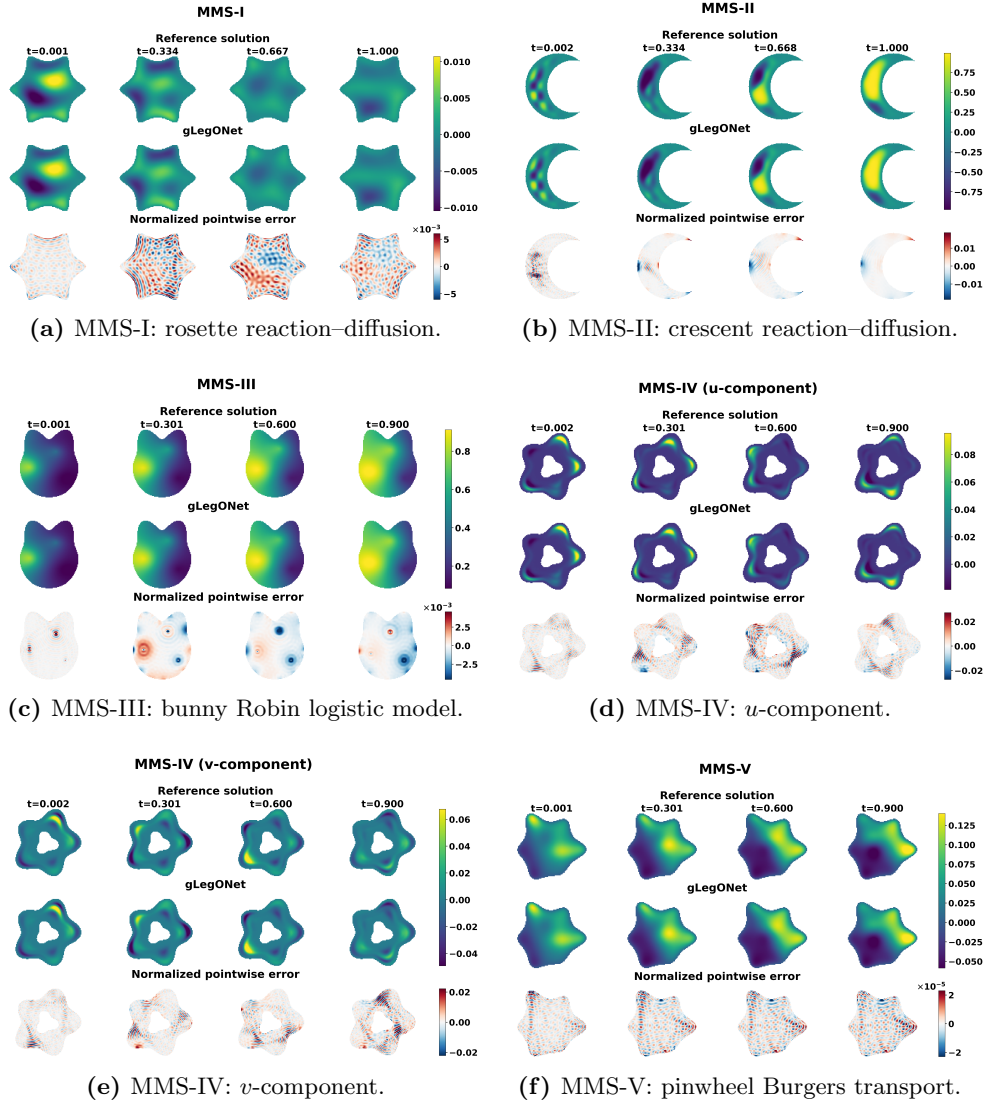

\centering
\scriptsize

\begin{subfigure}[t]{0.49\textwidth}
    \centering
    \includegraphics[width=\linewidth]{mms1_ours_ref_pred_error_3x4_gLegONet.png}
    \caption{MMS-I: rosette reaction--diffusion.}
\end{subfigure}
\hfill
\begin{subfigure}[t]{0.49\textwidth}
    \centering
    \includegraphics[width=\linewidth]{mms2_ours_ref_pred_error_3x4_gLegONet.png}
    \caption{MMS-II: crescent reaction--diffusion.}
\end{subfigure}

\vspace{2mm}

\begin{subfigure}[t]{0.49\textwidth}
    \centering
    \includegraphics[width=\linewidth]{mms3_ours_ref_pred_error_3x4_gLegONet.png}
    \caption{MMS-III: bunny Robin logistic model.}
\end{subfigure}
\hfill
\begin{subfigure}[t]{0.49\textwidth}
    \centering
    \includegraphics[width=\linewidth]{mms4_ours_ref_pred_error_3x4_u_gLegONet.png}
    \caption{MMS-IV: \(u\)-component.}
\end{subfigure}

\vspace{2mm}

\begin{subfigure}[t]{0.49\textwidth}
    \centering
    \includegraphics[width=\linewidth]{mms4_v_ours_ref_pred_error_3x4_v_gLegONet.png}
    \caption{MMS-IV: \(v\)-component.}
\end{subfigure}
\hfill
\begin{subfigure}[t]{0.49\textwidth}
    \centering
    \includegraphics[width=\linewidth]{mms5_ours_ref_pred_error_3x4_gLegONet.png}
    \caption{MMS-V: pinwheel Burgers transport.}
\end{subfigure}

\caption{Transferred-block rollouts on the manufactured-solution benchmarks. Each panel shows the reference solution, the Geometry-Aware LegONet prediction and the normalized pointwise error at representative times. The five tests cover Dirichlet, Neumann, Robin, mixed and transport-dominated boundary settings. }
\label{fig:mms_ours_rollout_3x4_all}
\end{figure}

\begin{figure}[p]
\centering
\includegraphics[width=0.58\textwidth]{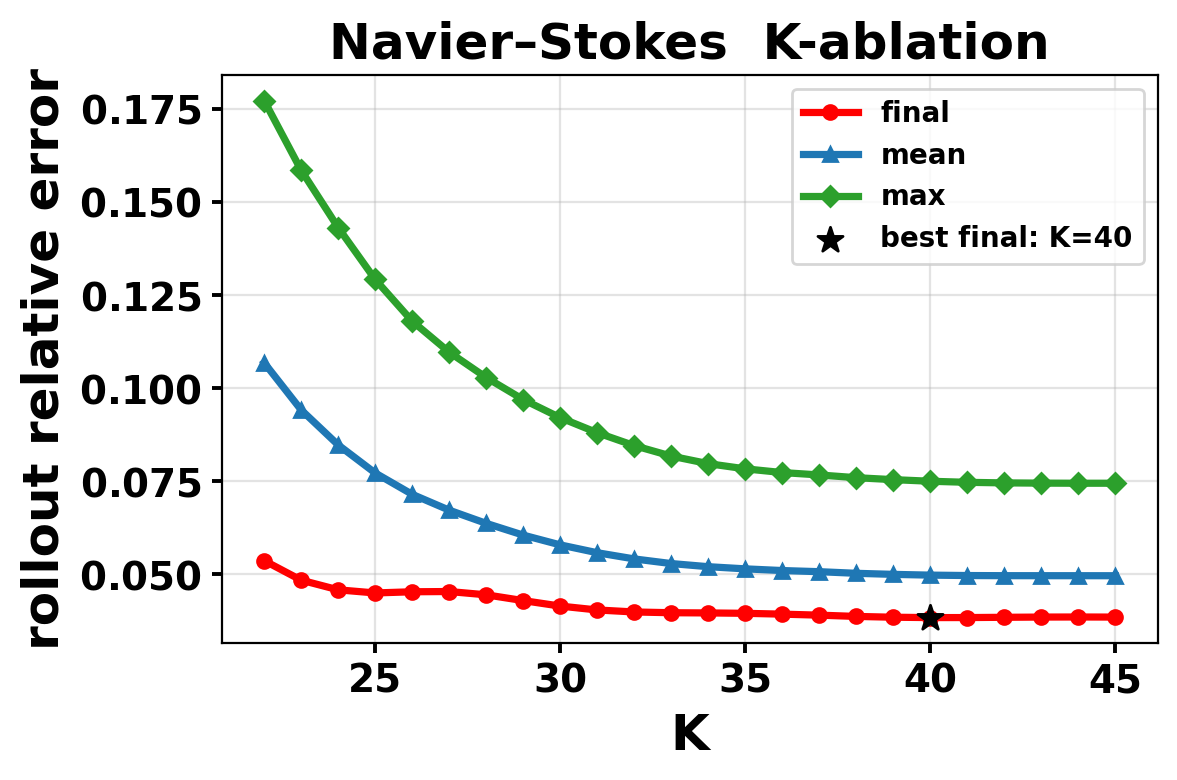}
\caption{ Resolution ablation for vorticity transport past an embedded cylinder. Increasing \(K\) enriches the ambient coefficient space and reduces the mean and maximum relative errors. The final-time error saturates near \(K=40\), showing that boundary enforcement remains stable while the interior rollout is limited by ambient resolution and retained admissible rank. }
\label{fig:ns_smooth_left_k_ablation}
\end{figure}

\begin{tableorg}[p]
\centering
\scriptsize
\setlength{\tabcolsep}{1.8pt}
\renewcommand{\arraystretch}{1.08}
\caption{ Baseline comparison on the five manufactured-solution benchmarks. The time column reports the target-stage cost in seconds: for PINN, FNO and UNO this is the target training time, whereas for gLegONet it is the target-domain assembly-and-rollout time, including boundary null-space construction, reduced operator assembly and time integration, with no target-domain neural training. Final, mean and max are relative \(L^2(\Omega)\) errors over the rollout. \(E_{\partial}^{\rm rms}\) is the RMS residual of the prescribed boundary operator, and \(E_{\partial}^{\rm sc}\) is the corresponding solution-scaled boundary residual. }
\label{tab:mms_baseline_comparison}
\resizebox{\linewidth}{!}{%
\begin{tabular}{llrrrrrr}
\toprule
Test & Method
& Time
& Final & Mean & Max &
\(E_{\partial}^{\rm rms}\) & \(E_{\partial}^{\rm sc}\) \\
& & (s)
& error & error & error & & \\
\midrule
MMS-I
& \textbf{gLegONet} & \(\mathbf{2.12}\)
& \(\mathbf{1.76\times10^{-3}}\) & \(\mathbf{2.12\times10^{-3}}\)
& \(\mathbf{3.16\times10^{-3}}\) & \(\mathbf{9.29\times10^{-11}}\)
& \(\mathbf{2.97\times10^{-8}}\) \\
& PINN & \(5.36\times10^3\)
& \(1.77\times10^{-1}\) & \(1.65\times10^{-1}\)
& \(5.37\times10^{-1}\) & \(8.43\times10^{-6}\)
& \(2.69\times10^{-3}\) \\
& FNO & \(6.20\times10^2\)
& \(1.59\times10^{-1}\) & \(9.20\times10^{-2}\)
& \(1.79\times10^{-1}\) & \(3.17\times10^{-4}\)
& \(1.01\times10^{-1}\) \\
& UNO & \(1.96\times10^3\)
& \(1.03\times10^{-2}\) & \(1.02\times10^{-2}\)
& \(2.42\times10^{-2}\) & \(2.07\times10^{-4}\)
& \(6.62\times10^{-2}\) \\
\midrule
MMS-II
& \textbf{gLegONet} & \(\mathbf{2.67}\)
& \(\mathbf{2.19\times10^{-3}}\) & \(\mathbf{2.98\times10^{-3}}\)
& \(\mathbf{4.78\times10^{-3}}\) & \(\mathbf{1.01\times10^{-7}}\)
& \(\mathbf{2.60\times10^{-7}}\) \\
& PINN & \(4.64\times10^3\)
& \(1.64\times10^{-2}\) & \(1.22\times10^{-2}\)
& \(1.73\times10^{-2}\) & \(1.05\times10^{-2}\)
& \(2.68\times10^{-2}\) \\
& FNO & \(1.49\times10^3\)
& \(3.65\times10^{-2}\) & \(5.39\times10^{-2}\)
& \(2.01\times10^{-1}\) & \(7.43\times10^{-1}\)
& \(1.91\times10^{0}\) \\
& UNO & \(2.19\times10^3\)
& \(7.42\times10^{-3}\) & \(6.27\times10^{-3}\)
& \(2.70\times10^{-2}\) & \(7.39\times10^{-1}\)
& \(1.90\times10^{0}\) \\
\midrule
MMS-III
& \textbf{gLegONet} & \(\mathbf{1.13\times10^{1}}\)
& \(\mathbf{1.11\times10^{-3}}\) & \(\mathbf{1.36\times10^{-3}}\)
& \(\mathbf{4.27\times10^{-3}}\) & \(\mathbf{5.33\times10^{-9}}\)
& \(\mathbf{1.04\times10^{-8}}\) \\
& PINN & \(6.55\times10^3\)
& \(1.54\times10^{-1}\) & \(4.10\times10^{-2}\)
& \(1.54\times10^{-1}\) & \(3.40\times10^{-2}\)
& \(6.63\times10^{-2}\) \\
& FNO & \(9.42\times10^{2}\)
& \(2.03\times10^{-3}\) & \(2.25\times10^{-3}\)
& \(9.59\times10^{-3}\) & \(1.63\times10^{-2}\)
& \(3.19\times10^{-2}\) \\
& UNO & \(9.44\times10^{2}\)
& \(6.95\times10^{-3}\) & \(7.67\times10^{-3}\)
& \(1.58\times10^{-2}\) & \(8.94\times10^{-2}\)
& \(1.75\times10^{-1}\) \\
\midrule
MMS-IV
& \textbf{gLegONet} & \(\mathbf{8.89}\)
& \(\mathbf{5.77\times10^{-3}}\) & \(\mathbf{6.12\times10^{-3}}\)
& \(\mathbf{7.70\times10^{-3}}\) & \(\mathbf{4.51\times10^{-10}}\)
& \(\mathbf{1.93\times10^{-8}}\) \\
& PINN & \(8.11\times10^3\)
& \(8.54\times10^{-2}\) & \(3.47\times10^{-2}\)
& \(8.54\times10^{-2}\) & \(1.41\times10^{-3}\)
& \(8.57\times10^{-2}\) \\
& FNO & \(4.52\times10^3\)
& \(1.79\times10^{-2}\) & \(9.60\times10^{-3}\)
& \(2.65\times10^{-2}\) & \(9.37\times10^{-4}\)
& \(9.73\times10^{-2}\) \\
& UNO & \(3.32\times10^3\)
& \(1.45\times10^{-2}\) & \(8.43\times10^{-3}\)
& \(1.46\times10^{-2}\) & \(1.00\times10^{-3}\)
& \(5.30\times10^{-2}\) \\
\midrule
MMS-V
& \textbf{gLegONet} & \(\mathbf{1.95\times10^{1}}\)
& \(\mathbf{5.60\times10^{-4}}\) & \(\mathbf{5.06\times10^{-4}}\)
& \(\mathbf{5.71\times10^{-4}}\) & \(\mathbf{3.91\times10^{-11}}\)
& \(\mathbf{8.17\times10^{-10}}\) \\
& PINN & \(6.56\times10^3\)
& \(1.50\times10^{-3}\) & \(1.61\times10^{-3}\)
& \(4.28\times10^{-3}\) & \(4.61\times10^{-5}\)
& \(7.92\times10^{-4}\) \\
& FNO & \(2.25\times10^2\)
& \(1.14\times10^{-2}\) & \(4.71\times10^{-3}\)
& \(1.14\times10^{-2}\) & \(3.85\times10^{-3}\)
& \(1.22\times10^{-1}\) \\
& UNO & \(2.26\times10^3\)
& \(1.39\times10^{-2}\) & \(6.65\times10^{-3}\)
& \(1.39\times10^{-2}\) & \(3.86\times10^{-3}\)
& \(1.22\times10^{-1}\) \\
\botrule
\end{tabular}%
}
\end{tableorg}

\begin{tableorg}[p]
\centering
\small
\caption{ Parameterized manufactured families used to train the supervised FNO and UNO baselines. Target parameters are the held-out values used in the five MMS benchmarks; training ranges define the surrounding amortized trajectory families. }
\label{tab:operator_family_ranges}
\setlength{\tabcolsep}{4pt}
\begin{tabular}{lll}
\toprule
Test & Target parameters & Training ranges \\
\midrule
MMS-I &
\(A=0.42,\ \omega=12\) &
\(A\in[0.30,0.55],\ \omega\in[9,15]\) \\
MMS-II &
\(A=1.0,\ \omega=7.5\) &
\(A\in[0.80,1.20],\ \omega\in[6.50,8.50]\) \\
MMS-III &
\(c_f=3.2,\ \alpha_s=3.0,\ A_{\rm rip}=0.05\) &
\(c_f\in[2.70,3.70],\ \alpha_s\in[2.40,3.60],\
A_{\rm rip}\in[0.02,0.08]\) \\
MMS-IV &
\(A_u=4.0,\ A_v=3.2,\ \omega=4.6\) &
\(A_u\in[3.40,4.60],\ A_v\in[2.70,3.70],\
\omega\in[3.80,5.40]\) \\
MMS-V &
\(A=2.6,\ \omega=4.0\) &
\(A\in[2.00,3.20],\ \omega\in[3.20,4.80]\) \\
\botrule
\end{tabular}
\end{tableorg}

\begin{tableorg}[p]
\centering
\scriptsize
\setlength{\tabcolsep}{3.6pt}
\renewcommand{\arraystretch}{1.08}
\caption{ Downstream baseline comparison on nonlinear target-domain rollouts. The time column reports the target-stage cost in seconds: for PINN, FNO and UNO this is the target training time, whereas for gLegONet it is the target-domain assembly-and-rollout time without target-domain neural training. Final, mean and max are relative \(L^2\) errors against the independent FEM reference for Allen--Cahn and the finite-volume reference for Burgers. \(E_{\partial}^{\rm rms}\) is the dense boundary residual. The diagnostic column reports the structure-related quantity: \(D_{\rm mass}^{\max}\) for mass conservation in Allen--Cahn and \(R_E^{\rm final}\) for the final kinetic-energy budget residual in Burgers. }
\label{tab:downstream_ac_burgers_baselines}
\resizebox{\linewidth}{!}{%
\begin{tabular}{lrrrrrr}
\toprule
Method
& Time
& Final
& Mean
& Max
& \(E_{\partial}^{\rm rms}\)
& Diagnostic \\
& (s)
& error
& error
& error
&
&  \\
\midrule
\multicolumn{7}{l}{\textit{Volume-constrained Allen--Cahn on a disk}} \\
\midrule
\textbf{gLegONet}
& \(\mathbf{3.05{\times}10^{1}}\)
& \(\mathbf{2.98{\times}10^{-2}}\)
& \(\mathbf{3.18{\times}10^{-2}}\)
& \(\mathbf{5.83{\times}10^{-2}}\)
& \(\mathbf{3.47{\times}10^{-9}}\)
& \(\mathbf{3.63{\times}10^{-5}}\) \\

PINN
& \(8.22{\times}10^{2}\)
& \(4.81{\times}10^{-2}\)
& \(4.94{\times}10^{-2}\)
& \(7.61{\times}10^{-2}\)
& \(2.49{\times}10^{-3}\)
& \(8.56{\times}10^{-3}\) \\

FNO
& \(1.24{\times}10^{4}\)
& \(1.19{\times}10^{-1}\)
& \(9.31{\times}10^{-2}\)
& \(1.62{\times}10^{-1}\)
& \(7.57{\times}10^{-3}\)
& \(1.37{\times}10^{-3}\) \\

UNO
& \(1.06{\times}10^{4}\)
& \(6.31{\times}10^{-2}\)
& \(5.39{\times}10^{-2}\)
& \(1.62{\times}10^{-1}\)
& \(3.14{\times}10^{-2}\)
& \(1.01{\times}10^{-2}\) \\

\midrule
\multicolumn{7}{l}{\textit{Inner-square vector Burgers}} \\
\midrule
\textbf{gLegONet}
& \(\mathbf{3.88{\times}10^{1}}\)
& \(\mathbf{1.83{\times}10^{-2}}\)
& \(\mathbf{1.17{\times}10^{-2}}\)
& \(\mathbf{1.83{\times}10^{-2}}\)
& \(\mathbf{1.91{\times}10^{-13}}\)
& \(\mathbf{3.63{\times}10^{-7}}\) \\

PINN
& \(1.31{\times}10^{3}\)
& \(5.58{\times}10^{-2}\)
& \(3.10{\times}10^{-2}\)
& \(5.58{\times}10^{-2}\)
& \(4.23{\times}10^{-3}\)
& \(4.21{\times}10^{-2}\) \\

FNO
& \(1.61{\times}10^{4}\)
& \(1.89{\times}10^{-1}\)
& \(1.34{\times}10^{-1}\)
& \(1.90{\times}10^{-1}\)
& \(3.88{\times}10^{-3}\)
& \(1.44{\times}10^{-2}\) \\

UNO
& \(1.19{\times}10^{4}\)
& \(3.62{\times}10^{-2}\)
& \(2.88{\times}10^{-2}\)
& \(3.62{\times}10^{-2}\)
& \(2.10{\times}10^{-3}\)
& \(2.55{\times}10^{-2}\) \\

Direct inner-square
& \(3.25{\times}10^{1}\)
& \(1.83{\times}10^{-2}\)
& \(1.17{\times}10^{-2}\)
& \(1.83{\times}10^{-2}\)
& \(1.17{\times}10^{-14}\)
& \(2.59{\times}10^{-8}\) \\
\botrule
\end{tabular}%
}
\end{tableorg}

\begin{tableorg}[p]
\centering
\scriptsize
\setlength{\tabcolsep}{4.0pt}
\renewcommand{\arraystretch}{1.08}
\caption{ Representative entries from the \(K\)-ablation for the Navier--Stokes cylinder-wake benchmark. Final, mean and max are relative \(L^2(\Omega)\) errors against the masked-grid reference. \(E_{\partial}^{\rm rms}\) is the RMS dense-boundary residual over saved rollout states, and Rank is the retained boundary-adapted dimension. }
\label{tab:ns_smooth_left_ablation}
\resizebox{\linewidth}{!}{%
\begin{tabular}{lrrrrr}
\toprule
Method & Final error & Mean error & Maximum error
& \(E_{\partial}^{\rm rms}\)
& Rank \\
\midrule
gLegONet (\(K=22\))
& \(5.36{\times}10^{-2}\)
& \(1.07{\times}10^{-1}\)
& \(1.77{\times}10^{-1}\)
& \(1.41{\times}10^{-11}\)
& 1375 \\

gLegONet (\(K=31\))
& \(4.04{\times}10^{-2}\)
& \(5.58{\times}10^{-2}\)
& \(8.80{\times}10^{-2}\)
& \(7.17{\times}10^{-12}\)
& 2811 \\

gLegONet (\(K=39\))
& \(3.84{\times}10^{-2}\)
& \(4.99{\times}10^{-2}\)
& \(7.54{\times}10^{-2}\)
& \(1.31{\times}10^{-12}\)
& 4545 \\

gLegONet (\(K=40\))
& \(3.83{\times}10^{-2}\)
& \(4.97{\times}10^{-2}\)
& \(7.50{\times}10^{-2}\)
& \(1.44{\times}10^{-12}\)
& 4789 \\

gLegONet (\(K=45\))
& \(3.85{\times}10^{-2}\)
& \(4.96{\times}10^{-2}\)
& \(7.44{\times}10^{-2}\)
& \(1.59{\times}10^{-12}\)
& 6099 \\
\botrule
\end{tabular}%
}
\end{tableorg}

\begin{tableorg}[p]
\centering
\scriptsize
\setlength{\tabcolsep}{5pt}
\renewcommand{\arraystretch}{1.08}
\caption{
Resolution ablation for the clamped Swift--Hohenberg benchmark on the
smooth heart-shaped domain.
}
\label{tab:heart_sh_k_ablation}
\begin{tabular}{crrrrr}
\toprule
Method 
& Final error
& Mean error
& Maximum error
& \(E_{\partial}^{\rm rms}\)
& Rank \\
\midrule
gLegONet (\(K=22\))
& \(1.47{\times}10^{-2}\)
& \(1.76{\times}10^{-2}\)
& \(4.62{\times}10^{-2}\)
& \(1.53{\times}10^{-8}\)
& \(445\) \\

gLegONet (\(K=24\))
& \(1.08{\times}10^{-2}\)
& \(1.33{\times}10^{-2}\)
& \(3.66{\times}10^{-2}\)
& \(8.87{\times}10^{-9}\)
& \(522\) \\

gLegONet (\(K=26\))
& \(9.43{\times}10^{-3}\)
& \(1.15{\times}10^{-2}\)
& \(2.98{\times}10^{-2}\)
& \(3.44{\times}10^{-9}\)
& \(610\) \\

gLegONet (\(K=28\))
& \(8.39{\times}10^{-3}\)
& \(9.87{\times}10^{-3}\)
& \(2.43{\times}10^{-2}\)
& \(4.06{\times}10^{-9}\)
& \(698\) \\

gLegONet (\(K=30\))
& \(7.01{\times}10^{-3}\)
& \(8.33{\times}10^{-3}\)
& \(2.08{\times}10^{-2}\)
& \(3.51{\times}10^{-9}\)
& \(795\) \\
\botrule
\end{tabular}
\end{tableorg}


\end{document}